\documentclass[11pt, a4paper]{article}
\usepackage[letterpaper,top=2.5cm,bottom=2.5cm,left=1.25cm,right=1.25cm,marginparwidth=1.75cm]{geometry}

\usepackage{amsmath, amssymb, amsthm} 
\usepackage{hyperref} 
\usepackage{float} 
\usepackage{booktabs, multirow} 
\usepackage{enumitem}
\usepackage{changepage} 
\usepackage{enumitem} 
\usepackage{dsfont}
\usepackage[table]{xcolor}
\usepackage{adjustbox}  
\usepackage{array}      
\usepackage{comment} 
\usepackage{subcaption} 
\usepackage{tikz}
\usepackage[numbers]{natbib}

\usepackage{changepage}
\usepackage{hyperref}
\usepackage{todonotes}
\usepackage{makecell} 

\usepackage{rotating}
\theoremstyle{plain}
\newtheorem{proposition}{Proposition}
\newtheorem{lemma}{Lemma}
\newtheorem{corollary}{Corollary}

\theoremstyle{definition}
\newtheorem{example}{Example}
\newtheorem{instance}{Instance}

\usepackage{siunitx}
\usepackage{graphicx}
\usepackage{tabularx}
\usepackage{placeins}

\usepackage[affil-it]{authblk}

\definecolor{resultgreen}{HTML}{A9D18E}
\definecolor{resultyellow}{HTML}{FFC000}
\definecolor{resultblue}{HTML}{8EA9DB}

\begin{document}
\title{Mixed-integer programming formulations for the quadratic bin packing problem: linearization, hybridization, and enhancements}
\author{Fatih Burak Ak\c{c}ay\thanks{Department of Econometrics and OR, Tilburg University, The Netherlands, email: F.B.Akcay@tilburguniversity.edu} \qquad
Maxence Delorme\thanks{Department of Econometrics and OR, Tilburg University, The Netherlands, email: M.Delorme@tilburguniversity.edu} \qquad
Renata Sotirov\thanks{Department of Econometrics and OR, Tilburg University, The Netherlands, email: r.sotirov@tilburguniversity.edu}
}

\date{ }
\noindent
\maketitle
\vspace*{-7ex}

\abstract{The quadratic bin packing problem (QBPP) extends the classical bin packing problem by incorporating pairwise dissimilarity values between items assigned to the same bin.
The objective is to minimize the total cost, composed of the cost associated with the number of bins used and the sum of pairwise dissimilarity values for items assigned to the same bin, subject to the bin-capacity constraints.
We study the quadratic programming formulation together with a broad family of equivalent mixed-integer linear programming (MILP) formulations for the QBPP based on classical linearization techniques.
We analyze the relationships among the linear programming relaxations of the MILP formulations and strengthen several of these relaxations by exploiting the structure of the QBPP.
We further introduce the first hybrid exact formulation for the QBPP that combines elements of two distinct modeling approaches. In particular, we show that the combination of explicit bin patterns and assignment-based representation performs particularly well.
For the hybrid formulation, we investigate the effect of item selection, and for both the hybrid and non-hybrid formulations, we study the impact of item ordering on formulation size, relaxation strength and computational performance.
Although a recently introduced tailored branch-and-price algorithm remains the state-of-the-art for the QBPP, we show that the proposed ordering-based non-hybrid and hybrid formulations substantially improve computational performance over the best MILP formulation known to the literature, in addition of solving 33 new benchmark instances to optimality.}


\maketitle
\section{Introduction}
Many allocation and partitioning problems require more than simply determining whether a collection of objects satisfies a given capacity restriction.
Pairwise interactions can therefore make certain combinations of objects desirable or undesirable, while separating them using additional resources may incur an extra cost.
Such trade-offs arise naturally in many practical settings.
In the advertising industry, for example, advertisements may be assigned to capacity-limited banners, where competing advertisements are undesirable on the same banner, whereas complementary advertisements may benefit from being displayed together.
Similarly, in transportation, products may be assigned to capacity-limited vehicles, where compatibility or handling requirements of the products can make certain combinations undesirable.
In both cases, this creates a trade-off between forming desirable combinations and limiting the number of resources used.
The \emph{quadratic bin packing problem} (QBPP) captures this trade-off.
In the QBPP, the objects to be assigned are represented as weighted items and the available resources as capacitated bins.
Every additional bin used incurs a fixed cost, while assigning two items to the same bin incurs a pairwise dissimilarity value, which may be positive (corresponding to penalty), negative (corresponding to profit) or zero.
The objective therefore balances the number of bins used against the composition of those bins.
This distinguishes the QBPP from the classical \emph{bin packing problem} (BPP), where only the number of bins used is optimized, and from clustering and partitioning problems in which the number of groups is typically fixed or not directly penalized, and only their composition is optimized.

More broadly, the QBPP lies at the intersection of several well-studied combinatorial optimization problems.  Depending on the assumptions on the capacity constraints, the number of bins and the dissimilarity values, several combinatorial problems can be recovered as special cases of the QBPP.
The main relations between selected problems and the QBPP are summarized in Table~\ref{tab:qbpp_related}, where we use bin packing terminology for consistency.
Here, bins represent groups such as partitions, clusters or knapsacks.
Capacity indicates whether a restriction is imposed on the total weight of each group.
The number of bins indicates whether the number of groups is fixed in advance or determined as part of the optimization problem.
Bin costs represent costs associated with forming these groups, while dissimilarity values represent pairwise interactions between items.
\begin{table}[!htbp]
  \centering
  \caption{Relationship of the QBPP to selected combinatorial problems.}
  \label{tab:qbpp_related}
  \scriptsize
  \setlength{\tabcolsep}{3.5pt}
  \begin{tabular}{lcccc}
  \toprule
  \textbf{Problem}
  & \textbf{capacity}
  & \textbf{\#bins}
  & \textbf{bin cost}
  & \textbf{dissimilarity values} \\
  \midrule

  Quadratic bin packing problem
  & $\checkmark$
  & variable
  & nonnegative
  & general \\

  \midrule\midrule

  Bin packing problem
  & $\checkmark$
  & variable
  & 1
  & 0 \\

  \midrule

  \makecell[l]{Bin packing problem with conflicts}
  & $\checkmark$
  & variable
  & 1
  & \makecell{$+\infty$ for conflicts, 0 otherwise} \\

  \midrule

  Graph coloring
  & $\times$
  & variable
  & 1
  & \makecell{$+\infty$ for edges, 0 otherwise} \\

  \midrule

  \makecell[l]{Node-capacitated  graph partitioning}
  & $\checkmark$
  & fixed
  & 0
  & general \\

  \midrule

  \makecell[l]{Minimum {$k$}-partition (max-$k$-cut)}
  & $\times$
  & fixed
  & 0
  & general \\

  \midrule

  \makecell[l]{Weighted correlation clustering}
  & $\times$
  & variable
  & 0
  & general \\
  \midrule

  \makecell[l]{Quadratick multiple knapsack problem}
  & $\checkmark$
  & fixed
  & ---
  & nonnegative \\
  \bottomrule
  \end{tabular}
\end{table}
The classical BPP, for which a broad range of mathematical models and exact algorithms has been developed (see, e.g.,~\cite{DIM16} for a review) is obtained by setting all dissimilarity values to zero and bin cost to one.
The BPP with conflicts extends this setting by forbidding certain pairs of items from sharing a bin.
In the QBPP, such conflicts can be represented by assigning a dissimilarity value of infinity (in practice, a prohibitively large value) to conflicting pairs and zero to all remaining pairs (see~\cite{MIMT10} for a detailed study).
Additionally, if the capacity restrictions are removed, each bin then corresponds to an independent set of the conflict graph, resulting in the graph coloring problem (see~\cite{L21} for a comprehensive overview of algorithms and applications).
Restricting the number of bins instead connects the QBPP to the node-capacitated graph partitioning problem, where vertices are partitioned into at most $k$ capacitated subsets while minimizing the cost of edges between subsets.
Since the total interaction value over all pairs is constant, minimizing between-subset interactions is equivalent to maximizing within-subset interactions. Hence, after reversing the signs of the interaction coefficients and setting the bin cost to zero, the problem can be described as a QBPP with at most $k$ bins (see~\cite{FMSWW98} for an extensive computational study).
Removing the capacity constraints and the bin-opening cost leads to the minimum $k$-partition problem, which is equivalently expressed as the max-$k$-cut problem when the objective is written in terms of edges between different subsets (for a recent comparison of formulations for max-$k$-cut, see~\cite{FVHTZ23}).
A further connection is with the weighted correlation clustering problem, where the number of clusters is not fixed and pairwise relations may favor either grouping or separating two objects (see~\cite{WH22} for a recent review).
Its objective can be expressed, up to an additive constant, as a sum of pairwise costs for objects assigned to the same cluster.
Hence, weighted correlation clustering corresponds to an uncapacitated QBPP without bin-opening costs.

The QBPP is also closely related to the \emph{quadratic multiple knapsack problem} (QMKP), which combines the capacity constraints of the multiple knapsack problem with individual and pairwise nonnegative profits (for a recent survey of the problem and its extensions, see~\cite{CILM22}).
The QMKP differs from the QBPP in several aspects: it (i) uses a fixed number of knapsacks, (ii) does not require all items to be assigned and (iii) maximizes individual and pairwise profits, which are typically assumed to be nonnegative.
In contrast, the QBPP requires every item to be assigned, allows pairwise dissimilarity values of either sign and includes a fixed cost for each bin used.
Although the QMKP is structurally closer to the QBPP, the \emph{quadratic knapsack problem} (QKP) is also relevant for us since several of its formulation and linearization techniques (which, to the best of our knowledge, have not yet been adapted to the QMKP) are also valid for the QBPP (see~\cite{GMT25} for a recent review).

In the present work, we study several classical linearization approaches for the QBPP, including those of~\cite{FP10}, \cite{GW74}, \cite{G75}, and \cite{AS86}.
The formulations based on the Fan--Pardalos and Glover--Woolsey linearizations have already been considered in the recent work of \cite{CLMI26} on the QBPP. We complement their results with a theoretical analysis of the corresponding \emph{linear programming} (LP) relaxations.
We show that, without symmetry-breaking constraints, the two QBPP formulations have the same optimal LP-relaxation value, whereas after symmetry-breaking constraints are incorporated, the Glover--Woolsey-based QBPP formulation dominates the Fan--Pardalos-based formulation and the dominance is strict for some instances.
To the best of our knowledge, we also provide the first exact QBPP formulation based on Glover's linearization and strengthen its LP relaxation by exploiting the specific structure of the QBPP.
Finally, we derive a linear integer program for the QBPP following the approach introduced by \cite{AS86} and establish its relationship with the other formulations.

Beyond these classical linearizations, we introduce to the best of our knowledge, the first hybrid \emph{mixed-integer linear programming} (MILP) formulation for the QBPP, combining elements of two distinct modeling structures. In particular, the hybrid formulation combines explicit bin patterns generated over a selected subset of items with an assignment-based representation that can accommodate all items.
We also study the role of item ordering in assignment-based formulations.
We show that the item ordering can substantially affect formulation size, LP-relaxation strength and computational performance. Furthermore, we demonstrate that an appropriate ordering improves both the non-hybrid and item-partitioned hybrid formulations.
Finally, we conduct an extensive computational study to compare the proposed formulations and evaluate how model strengthening techniques (such as valid inequalities), hybridization, item selection and item ordering affect the relaxation strength and computational performance, measured by the number of instances solved to optimality and the solution time.

The remainder of the paper is organized as follows.
Section~\ref{sec:litrev} reviews the literature related to the formulations and solution approaches considered in this paper.
Section~\ref{sec:probdef} formally defines the QBPP and introduces the notation used throughout the paper.
Section~\ref{sec:formulations} presents the MILP formulations based on classical linearization techniques, as well as the pattern-based and item-partitioned hybrid formulations, and analyzes the relationships among their LP relaxations.
Section~\ref{sec:valid} introduces additional valid inequalities derived from the underlying BPP.
Section~\ref{sec:comp_exp} reports the computational experiments, comparing the formulations and examining the effects of model strengthening techniques, item ordering and hybridization on relaxation strength and computational performance.

\section{Literature Review} \label{sec:litrev}

The problems related to the QBPP outlined in the previous section also provide useful connections in terms of modeling and solution methods.
In this section, however, we focus on a selection of these problems whose formulations and solution approaches are closely related to the methods studied in this paper.

The classical BPP has been studied through a broad range of mathematical formulations and exact algorithms.
While much of the recent literature on exact solution methods for the BPP has focused on \emph{branch-and-price} (B\&P) approaches (see, e.g., \cite{BCCF24}), the problem is also well known for the wide variety of \emph{integer linear programming} (ILP) formulations that have been proposed (see, e.g., \cite{DIM16}).
{Our formulations for the QBPP largely retain the assignment structure of the textbook BPP model (often attributed to \cite{K60})}, but differ in how the quadratic interactions are represented.
Similarly, both ILP and pattern-based formulations have been studied for the BPP with conflicts.
\cite{MIMT10} presented an assignment-based ILP formulation and a set-covering B\&P approach, while \cite{WLBL20} later developed a competitive branch-price-and-cut algorithm based on a set-partitioning formulation.
For the node-capacitated graph partitioning problem, \cite{FMSWW96} developed ILP formulations and valid inequalities, which were later incorporated into a branch-and-cut algorithm~\citep{FMSWW98}.

The QKP and QMKP are particularly relevant to the formulation techniques considered in this paper.
For the QKP, \cite{BS04} proposed two MILP formulations based on an adaptation of Glover's linearization~\citep{G75}.
Rather than introducing an auxiliary variable for every pair of items, they aggregate the quadratic contribution associated with each item, requiring only a linear number of additional variables and constraints.
The coefficients used in the linearization are tightened through auxiliary knapsack problems, with an enhanced variant additionally exploiting information from a feasible solution.
For the QMKP, \cite{GMRT21} investigated several polynomial-size formulations based on classical linearization techniques for \emph{binary quadratic programming} (BQP) formulations and analyzed the relationships among their continuous relaxations.
\cite{B20} proposed the first exact algorithm for the QMKP, based on an exponential-size formulation over feasible knapsack assignments solved by a B\&P algorithm.
The associated pricing problem is a QKP with additional constraints, for which a specialized reformulation and exact solution procedure were developed.
\cite{F22} subsequently proposed a branch-and-bound algorithm that branches directly on item-to-knapsack assignments.
Its main component is an upper bound obtained by distributing pairwise values among individual items, estimating their maximum potential contributions through fractional knapsack problems and solving a transportation problem.
The method outperformed the previously proposed exact algorithms on the benchmark instances considered.

The QBPP itself has only recently been introduced and studied as a standalone optimization problem by \cite{CLMI26}.
They proposed three MILP formulations together with enhanced versions having strengthened continuous relaxations.
They also introduced an exponential-size set-partitioning formulation and developed a tailored B\&P algorithm to solve it.
This approach is structurally similar to Bergman's B\&P method for the QMKP, as both rely on columns representing feasible subsets of items and lead to a quadratic knapsack-type pricing problem.
In the QBPP, however, the pricing problem may contain both positive and negative coefficients, for which \cite{CLMI26} developed a dedicated heuristic and exact solution procedures.
Their computational experiments show that the enhanced MILP formulations are effective on smaller instances, with the so-called ``Extended Fortet--Glover--Woolsey'' (EFGW) formulation performing best among them, whereas the B\&P approach performs best overall, particularly on larger and more challenging instances.
\section{Problem definition and notation} \label{sec:probdef}

In the QBPP, we are given a set of $n$ items, where each item $i \in \{1, \dots, n\}$ has a positive integer weight  $w_i$.
In addition, $b$ identical bins are available, each with positive integer capacity $W$ and a nonnegative integer usage cost $c$.
Let $D=(d_{ij})$ denote an $n \times n$ dissimilarity matrix, where $d_{ij}$ is the integer dissimilarity value associated with placing items $i$ and $j$ in the same bin.
Each item must be assigned to exactly one bin and the total weight of the items in any bin cannot exceed the bin capacity.
The objective is to minimize the total cost, defined as the sum of pairwise dissimilarity {values} between items packed in the same bin and the cost associated with the number of bins used.
Let $\mathcal{N} = \{1, 2, \ldots, n\}$ be the set of items and $\mathcal{B} = \{1, 2, \ldots, b\}$ the set of bins. 
We assume that:
\begin{itemize}
    \item each item can fit into a bin without exceeding its capacity: $0 < w_i \leq W$ ($ i \in \mathcal{N} $),
    \item the number of available bins is equal to the number of items: $b=n$ with $n\ge2$,
    \item the relative order of items in a bin does not affect the dissimilarity value, i.e., dissimilarity values are symmetric: $ d_{ij} = d_{ji}$ ($i, j \in \mathcal{N}$),
    \item an item has no self-dissimilarity value: $ d_{ii} = 0 $ ($ i \in \mathcal{N} $).
\end{itemize}

We introduce two small instances that will be used throughout the paper in the numerical examples.
\begin{instance} \label{inst1}
$n=4$, $c=8$, $W=24$, $(w_1,w_2,w_3,w_4)=(4,8,12,12)$ and
$(d_{12},d_{13},d_{14}, d_{23}, d_{24}, d_{34})= (-4,8,12,$ $-8,-12,8)$.
An optimal solution assigns items $1$, $2$ and $3$ to bin $1$ and item $4$ to bin $2$, yielding an optimal objective value of $12$.
\end{instance}
\begin{instance}
\label{inst2}
This instance has the same parameters as Instance~\ref{inst1}, except that $W=30$. The optimal solution described for Instance~\ref{inst1} remains optimal for Instance~\ref{inst2}, with the same optimal objective value of $12$.
\end{instance}

\section{Mixed-integer programming formulations} \label{sec:formulations}
In this section, we first present a BQP formulation of the QBPP.
We then introduce symmetry-breaking constraints and incorporate them into the formulations in the remainder of the section.

\paragraph{Binary quadratic programming formulation:}
Let the binary variable $x_{ik}$ denote whether item $i\ (i\in \mathcal{N}$) is packed in bin $k\ (k \in \mathcal{B})$ and let binary variable $y_k$ indicate whether bin $k\ (k \in \mathcal{B})$ is used. Then the BQP formulation for the QBPP is given by
\begin{align}
(\mathrm{F_{BQP}}) \quad \min \quad & \sum_{i=1}^{n-1} \sum_{j=i+1}^n \sum_{k \in \mathcal{B}} d_{ij}x_{ik}x_{jk} + c \sum_{k \in \mathcal{B}} y_k \label{BQP1}\\
\text{s.t.} \quad
& \sum_{k \in \mathcal{B}} x_{ik} = 1 &  i \in \mathcal{N} \label{BQP3}\\
& \sum_{i \in \mathcal{N}} w_i x_{ik} \leq W y_{k} &  k \in \mathcal{B} \label{BQP4} \\
& x_{ik} \leq y_k &   i \in \mathcal{N}, k \in \mathcal{B} \label{BQP2} \\
& y_k \in \{0,1\} & k \in \mathcal{B} \label{BQP6}\\
& x_{ik} \in \{0,1\} & i \in \mathcal{N}, k \in \mathcal{B}. \label{BQP5}
\end{align}
The objective function~\eqref{BQP1} minimizes the total dissimilarity among items assigned to the same bin and the cost of bin usage.
Constraints~\eqref{BQP3} require each item to be placed in exactly one bin, while constraints~\eqref{BQP4} impose the bin capacity restriction.
Although implied by the capacity constraints~\eqref{BQP4}, constraints~\eqref{BQP2} are included to strengthen the LP relaxation by ensuring that an item can be assigned to a bin only if the bin is used.

\begin{example}
\label{ex1}
Consider Instance \ref{inst1}.
An optimal solution of $\mathrm{F_{BQP}}$ is given by $x_{11}=x_{21}=x_{31}=x_{42}=1$ and $y_1=y_2=1,$
with all other variables equal to zero,
yielding an optimal objective value of $12$.
The same solution is also optimal for Instance~\ref{inst2}.
\end{example}

\paragraph{Symmetry-breaking constraints:}
Since all bins are identical in capacity and cost, the indices of the bins are interchangeable and have no impact on the objective value.
In other words, let $(\mathcal{N}_1, \mathcal{N}_2, \dots, \mathcal{N}_{\kappa})$, where $\kappa \le b$, be a partition of $\mathcal{N}$.
Suppose that a feasible solution assigns the items in subset $\mathcal{N}_{k}$ to bin $k$ for $k=1, 2, \dots, \kappa$.
Then, assigning the same subsets of items to any $\kappa$ bins yields an equivalent solution in terms of feasibility and objective value.
Hence, the solution is invariant under relabeling of the active bins.
Therefore, using $\kappa$ out of $b$ bins results in ${}_bP_\kappa = \frac{b!}{(b-\kappa)!}$ symmetric solutions.
To reduce symmetry in the model, we adopt the symmetry-breaking constraints from \cite{BNVM12}, originally proposed for the graph partitioning problem under capacity constraints:
\begin{align}
    x_{ik} &=0 &i\in \mathcal{N},k \in \mathcal{B}, i<k \label{SYM1}\\
    x_{ik} &\leq x_{kk} &i\in \mathcal{N},k \in \mathcal{B}, k<i. \label{SYM2}
\end{align}
Constraints~\eqref{SYM1} restrict item~$i$ to be packed only into the first $i$ bins and
constraints~\eqref{SYM2} ensure that
bin $k$ is nonempty if and only if item $k$ is packed into it.
Each variable $y_k$ can therefore be replaced by $x_{kk}$ and constraints \eqref{SYM2} can be used in place of constraints \eqref{BQP2}.
The capacity constraints~\eqref{BQP4} can then be written as
\begin{equation}
    \sum_{i\in\mathcal{N}} w_i x_{ik} \le W\, x_{kk} \qquad k \in \mathcal{B}. \label{CapSym}
\end{equation}
Accordingly, the summations in the assignment constraints~\eqref{BQP3} and capacity constraints~\eqref{CapSym} can be restricted to the relevant index sets, namely $k=1,\ldots,i$ and $i=k,\ldots,n$, respectively.
We denote the resulting BQP formulation for the QBPP by $\mathrm{F^S_{BQP}}$. \phantomsection \label{form:BQP}

For Instance~\ref{inst1} and~\ref{inst2}, the solution displayed in Example \ref{ex1} is infeasible for $\mathrm{F^S_{BQP}}$: item $4$ is assigned to bin $2$, whereas item $2$ is not, violating constraints~\eqref{SYM2}.
An optimal solution of $\mathrm{F^S_{BQP}}$ assigns items $1$, $2$, and $3$ to bin $1$ and item $4$ to bin $4$, leaving bins $2$ and $3$ empty.

\subsection{Mixed-integer linear programming formulations}
In quadratic optimization problems with binary variables, any polynomial expression can be reformulated as a linear one by introducing auxiliary variables together with additional linear constraints. For the QBPP, we derive formulations based on four linearization approaches: Fan--Pardalos, Glover--Woolsey, Adams--Sherali and Glover.
We also develop an item-partitioned hybrid formulation that combines a pattern-based representation with an assignment-based formulation.
In the remainder of this work, for a formulation $\mathrm{F_x}$, we denote its LP relaxation by $\mathrm{LP}(\mathrm{F_x})$. We use $\mathrm{opt}(\mathrm{F_x})$ and $\mathrm{opt}(\mathrm{LP}(\mathrm{F_x}))$ to denote the corresponding optimal objective values.
Numerical values reported in the examples are rounded to the nearest two decimal places.
We report the corresponding explicit optimal solutions in Appendix~\ref{appdx:num_examples}.

\paragraph{The Fan--Pardalos linearization approach:}
Following the  linearization framework of \cite{FP10}, we
introduce for each $i,j \in \mathcal{N}$ with $i<j$ an auxiliary variable $z_{ij}:=\sum_{k\in\mathcal{B}}x_{ik}x_{jk}$, which indicates whether items $i$ and $j$ are packed together.
We then replace the quadratic terms in the objective function with the auxiliary variables and introduce constraints
to ensure that $z_{ij}$ equals $1$ if items $i$ and $j$  are packed into the same bin, and 0 otherwise.
The resulting linear formulation of the QBPP is given as follows:
\begin{align}
(\mathrm{F_{FP}}) \quad \min \quad & \sum_{i=1}^{n-1} \sum_{j=i+1}^{n} d_{ij}z_{ij}+c\sum_{k \in \mathcal{B}}y_k \\
\text{s.t.}  \quad  & \eqref{BQP3}-\eqref{BQP5} \notag\\
    &x_{ik}+x_{jk}-1\leq z_{ij}   & i,j\in\mathcal{N}, k\in \mathcal{B}, i<j \label{FP1}\\
    &z_{ij} \leq 1-x_{ik}+x_{jk}  & i,j\in\mathcal{N}, k\in \mathcal{B}, i<j \label{FP2}\\
    &z_{ij} \leq 1+x_{ik}-x_{jk}  & i,j\in\mathcal{N}, k\in \mathcal{B}, i<j \label{FP3} \\
    &z_{ij} \ge 0  & i,j\in \mathcal{N}, i<j. \label{FP4}
\end{align}
The formulation without constraints~\eqref{BQP2} is presented by \cite{CLMI26}. The linearization constraints~\eqref{FP1}--\eqref{FP4} enforce that  $z_{ij}$  takes the correct value for each $i,j \in \mathcal{N}$ with $i<j$ (see \cite{M25}).
This linearization introduces $\mathcal{O}(n^2)$ auxiliary variables and $\mathcal{O}(n^2b)$ additional constraints. We denote the feasible solutions of $\mathrm{F_{FP}}$ with $(X,y, Z)$, where $X=(x_{ij})$, $y=(y_{k})$ and $Z=(z_{ij})$.
The following lemma establishes redundancy among the linearization constraints in $\mathrm{F_{FP}}$.
A similar result was obtained by~\cite{M25} for one of the constraint sets. For completeness, we include the full statement and proof.
\begin{lemma}
In $\mathrm{F_{FP}}$, constraint set \eqref{FP2} is redundant in the presence of \eqref{FP3}, and constraint set \eqref{FP3} is redundant in the presence of \eqref{FP2}.
\end{lemma}
\begin{proof}
    Consider a feasible solution $(X, y, Z)$ of $\mathrm{F_{FP}}$ without constraints~\eqref{FP3}.
    Suppose that for some $k \in \mathcal{B}$, we have $x_{ik} = 0$ and $x_{jk} = 1$.
    Then $z_{ij}$ is not directly constrained   in the absence of~\eqref{FP3}.
    However, by the assignment constraints \eqref{BQP3}, there exists a bin $k' \in \mathcal{B}$, with $k' \neq k$, such that $x_{ik'}=1$ and $x_{jk'}=0$.
    For this $k'$, constraints \eqref{FP2} and \eqref{FP4} imply that $z_{ij}=0$. Hence, constraints \eqref{FP3} are redundant.
    An analogous argument shows the redundancy of~\eqref{FP2} in the presence of constraints~\eqref{FP3}.
\end{proof}
Next, we show that $\mathrm{opt(LP(F_{FP}))}$ provides a weak bound.
\begin{proposition} \label{prop: LP_FP}
Let $(\bar X, \bar y, \bar Z)$ be defined as
\[
\bar x_{ik}:=\tfrac{1}{b}, \qquad
\bar y_k:=\max\!\Big\{\tfrac{1}{b}, \;\tfrac{\sum_{i\in\mathcal{N}} w_i}{W\,b}\Big\}, \qquad
\bar z_{ij}:=\begin{cases}1,& d_{ij}<0 \\[4pt]0,& d_{ij}\ge 0,  \end{cases} \quad i,j \in \mathcal{N},\, k\in\mathcal{B},\, i<j.
\]
Then $(\bar X, \bar y, \bar Z)$ is the optimal solution of $\mathrm{LP(F_{FP})}$
with objective value
$
\sum_{i=1}^{n-1}\sum_{\substack{j=i+1\\d_{ij}<0}}^{n} d_{ij} + \allowbreak c \max\!\left\{1, \;\tfrac{1}{W}\sum_{i\in\mathcal{N}} w_i\right\}.
$
\end{proposition}
\begin{proof}

It is straightforward to verify the feasibility of $(\bar X, \bar y, \bar Z)$.
To prove optimality, we derive a lower bound on the objective value of $\mathrm{LP(F_{FP})}$ and show that $(\bar X, \bar y, \bar Z)$ attains it.
We start with finding a lower bound on $\sum_{k\in\mathcal{B}}y_k$.
Constraints~\eqref{BQP3} and \eqref{BQP2} together imply that every feasible solution satisfies   $\sum_{k\in\mathcal{B}} y_k \ge 1$.
Summing the capacity constraints~\eqref{BQP4} over $k\in\mathcal{B}$ and using \eqref{BQP3} gives $\sum_{k\in\mathcal{B}}{y_k}\ge \frac{1}{W}\sum_{i\in\mathcal{N}}w_i$.
Combining these two bounds yields
\begin{equation}
   \sum_{k\in\mathcal{B}} y_k \ge \max\!\Big\{1, \frac{\sum_{i\in\mathcal{N}} w_i}{W}\Big\}. \label{doteq}
\end{equation}
Next, we find a lower bound on $\sum_{i=1}^{n-1}\sum_{j=i+1}^n d_{ij} z_{ij}$. Since $0\le z_{ij}\le1$ $(i,j \in \mathcal{N}, i<j)$, which follows trivially from \eqref{FP1}-\eqref{FP4}, we have
\begin{equation}
d_{ij}z_{ij}\;\ge\;
\begin{cases}
d_{ij},& d_{ij}<0 \\
0,& d_{ij}\ge0
\end{cases} \quad \implies \quad
\sum_{i=1}^{n-1}\sum_{j=i+1}^n d_{ij} z_{ij} \ge \sum_{i=1}^{n-1}\sum_{\substack{j=i+1\\d_{ij}<0}}^{n} d_{ij}. \label{threedots}
\end{equation}
Combining \eqref{doteq} and \eqref{threedots}, we obtain a lower bound on the objective value of $\mathrm{LP(F_{FP})}$:
\begin{equation}
 \sum_{i=1}^{n-1}\sum_{\substack{j=i+1\\d_{ij}<0}}^{n} d_{ij} \;+\; c\max\Big\{1, \;\tfrac{1}{W}\sum_{i\in\mathcal{N}} w_i\Big\}. \label{obj_bound}
\end{equation}
Evaluating the objective at $(\bar X, \bar y, \bar Z)$ yields exactly the value in \eqref{obj_bound}, thereby proving that $(\bar X, \bar y, \bar Z)$ is optimal.
\end{proof}

It follows from the previous proposition that the optimal solution $(\bar X, \bar y, \bar Z)$  of $\mathrm{LP(F_{FP}})$ involves only negative interaction values,
and that the variables $\bar X$ and $\bar y$ take averaged values.
Since the optimal solution  of the LP relaxation fails to capture the combinatorial structure of the original problem, the resulting  bound is weak.
One can strengthen an LP relaxation of so highly symmetric problem by e.g.,~adding symmetry-breaking constraints.
In the sequel, we add the symmetry-breaking constraints~\eqref{SYM1} and~\eqref{SYM2} to  $\mathrm{F_{FP}}$, replace each variable~$y_k$ with~$x_{kk}$, and substitute the capacity constraints~\eqref{BQP4} with~\eqref{CapSym}.
The resulting formulation of the QBPP is presented below and is also given by~\cite{CLMI26}.
\[
(\mathrm{F^S_{FP}}) \quad \min \sum_{i=1}^{n-1} \sum_{j=i+1}^{n} d_{ij} z_{ij}
+ c \sum_{k \in \mathcal{B}} x_{kk}
\quad \text{s.t.} \quad
\eqref{BQP3}, \eqref{BQP5}-\eqref{CapSym}, \eqref{FP1}-\eqref{FP4}.
\]
For completeness, the reduced formulation of $\mathrm{F^S_{FP}}$, obtained by eliminating the variables fixed to zero and by restricting the constraint index sets accordingly, is given in Appendix~\ref{appdx:FP}.

\paragraph{The Glover--Woolsey linearization approach:}
The linearization of products of binary variables dates back to \cite{F60} and was later refined by \cite{GW74}. We apply the Glover--Woolsey linearization approach to the QBPP by introducing variables  $g_{ijk} = x_{ik} x_{jk}$ for $i,j \in \mathcal{N}$, $k \in \mathcal{B}$, with $i<j$, where $g_{ijk}$ indicates whether items $i$ and $j$ are packed together in bin $k$.
Replacing the quadratic terms in the objective with these variables and adding linear constraints to enforce $g_{ijk} = x_{ik} x_{jk}$ yields the following linear formulation for the QBPP:
\begin{align}
(\mathrm{F_{GW}}) \quad \min \quad & \sum_{i=1}^{n-1} \sum_{j=i+1}^{n} \sum_{k \in \mathcal{B}} d_{ij}g_{ijk}+c\sum_{k \in \mathcal{B}}y_{k} \label{GWobj} \\
\text{s.t.}  \quad  & \eqref{BQP3}-\eqref{BQP5} \notag\\
&g_{ijk} \leq x_{ik} & i,j\in \mathcal{N}, k\in \mathcal{B}, i<j \label{GWL1}\\
&g_{ijk} \leq x_{jk} & i,j\in \mathcal{N}, k\in \mathcal{B}, i<j \label{GWL2} \\
&x_{ik} + x_{jk} - 1 \leq g_{ijk} & i,j\in \mathcal{N}, k\in \mathcal{B}, i<j \label{GWL3}\\
&g_{ijk} \ge 0 & i,j\in \mathcal{N}, k\in \mathcal{B}, i<j. \label{GWL4}
\end{align}
\cite{CLMI26} presents a similar formulation, where constraints \eqref{BQP2} are omitted.
The linearization constraints~\eqref{GWL1}--\eqref{GWL4} ensure that $g_{ijk}$ takes the correct value for each $i,j \in \mathcal{N}$, $k \in \mathcal{B}$ with $i<j$.
The linearization introduces $\mathcal{O}(n^2 b)$ auxiliary variables and $\mathcal{O}(n^2 b)$ additional constraints.
We denote feasible solutions of $\mathrm{F_{GW}}$ by $(X,y,G)$, where $X=(x_{ik})$, $y=(y_k)$ and $G=(g_{ijk})$.
Note that due to constraints~\eqref{GWL1}--\eqref{GWL4} and the integrality of $X$, the variables $g_{ijk}$ are implicitly binary.
Below, we show that $\mathrm{LP(F_{GW})}$ and $\mathrm{LP(F_{FP})}$ have the same optimal objective value (see Proposition~\ref{prop: LP_FP}).
\begin{proposition} \label{prop: LP_GW}
Let $(\bar X, \bar y, \bar G)$ be defined as
\[
\bar x_{ik}:=\tfrac{1}{b}, \qquad
\bar y_k:=\max\!\Big\{\tfrac{1}{b}, \;\tfrac{\sum_{i\in\mathcal{N}} w_i}{Wb}\Big\}, \qquad
\bar g_{ijk}:=\begin{cases}\tfrac{1}{b},& d_{ij}<0 \\[4pt]0,& d_{ij}\ge0,\end{cases} \quad i,j \in \mathcal{N}, k\in\mathcal{B}, i<j.
\]
Then $(\bar X, \bar y, \bar G)$ is the optimal solution of $\mathrm{LP(F_{GW})}$, with objective value
$
\sum_{i=1}^{n-1}\sum_{\substack{j=i+1\\d_{ij}<0}}^{n} d_{ij} +\allowbreak c \max\!\left\{1, \;\tfrac{1}{W}\sum_{i\in\mathcal{N}} w_i\right\}.
$
\end{proposition}
\begin{proof}
The feasibility of $(\bar X, \bar y, \bar G)$ can be easily verified.
We first derive a lower bound on the objective value of $\mathrm{LP(F_{GW})}$ and then show that $(\bar X, \bar y, \bar G)$ attains it.
From \eqref{BQP3}, \eqref{BQP4} and  \eqref{BQP2} it follows
\begin{align}\label{lby}
   \sum_{k\in\mathcal{B}} y_k \ge \max\!\Big\{1, \frac{\sum_{i\in\mathcal{N}} w_i}{W}\Big\},
\end{align}
see also derivation of the inequality \eqref{doteq}.

Next, we find a lower bound on the term $\sum_{i=1}^{n-1} \sum_{j=i+1}^{n} \sum_{k \in \mathcal{B}} d_{ij}g_{ijk}$. Since  $0\le g_{ijk}\le1$, which follows from \eqref{GWL1}-\eqref{GWL4}, we have
\begin{equation}
d_{ij}\sum_{k\in\mathcal{B}} g_{ijk}\;\ge\;
\begin{cases}
d_{ij},& d_{ij}<0 \\
0,& d_{ij}\ge0
\end{cases}
\quad \implies \quad
\sum_{i=1}^{n-1}\sum_{\substack{j=i+1\\d_{ij}<0}}^{n} d_{ij} \sum_{k\in\mathcal{B}} g_{ijk} \ge \sum_{i=1}^{n-1}\sum_{\substack{j=i+1\\d_{ij}<0}}^{n} d_{ij}. \label{twodoteq}
\end{equation}
By combining \eqref{lby} and \eqref{twodoteq}, we obtain
$\displaystyle\sum_{i=1}^{n-1}\sum_{\substack{j=i+1\\d_{ij}<0}}^{n} d_{ij} + c \max\!\Big\{1, \;\tfrac{1}{W}\sum_{i\in\mathcal{N}} w_i\Big\}$
as a lower bound on the objective value of $\mathrm{LP(F_{GW})}$.
Since the  objective value at $(\bar X, \bar y, \bar G)$  equals this lower bound, we conclude that   $(\bar X, \bar y, \bar G)$ is optimal.
\end{proof}

Next, we summarize the results of Proposition \ref{prop: LP_FP}  and Proposition \ref{prop: LP_GW}.
\begin{corollary}\label{prop: same_opt}
The optimal value of $\mathrm{LP(F_{GW})}$ is equal to the optimal value of $\mathrm{LP(F_{FP})}$.
\end{corollary}

It follows that both relaxations fail to capture the combinatorial structure of the QBPP.
As a result, the LP relaxation of $\mathrm{F_{GW}}$, like that of $\mathrm{F_{FP}}$, remains weak unless additional structure is imposed.
To obtain a formulation for the QBPP whose LP relaxation is stronger than that of $\mathrm{F_{GW}}$, we add the symmetry-breaking constraints \eqref{SYM1} and \eqref{SYM2}, replace each variable $y_k$ with $x_{kk}$, and substitute the capacity constraints \eqref{BQP4} with \eqref{CapSym}.
The resulting QBPP formulation is presented below and also provided by \cite{M25}.
\[
(\mathrm{F^S_{GW}}) \quad  \min \sum_{i=1}^{n-1} \sum_{j=i+1}^{n} \sum_{k \in \mathcal{B}} d_{ij}g_{ijk}+c\sum_{k \in \mathcal{B}}x_{kk}
\quad \text{s.t.} \quad
\eqref{BQP3}, \eqref{BQP5}-\eqref{CapSym}, \eqref{GWL1}-\eqref{GWL4} .
\]
The reduced formulation of $\mathrm{F^S_{GW}}$, obtained after removing the variables fixed to zero and reducing the associated index sets, is provided in Appendix~\ref{appdx:GW}.

Next, we show that $\mathrm{LP(F^S_{FP})}$ is dominated by $\mathrm{LP(F^S_{GW})}$.
\begin{proposition} \label{prop:LPrelstrong}
\leavevmode
\begin{enumerate}
    \item[(i)] $\mathrm{opt(LP(F^S_{GW}))} \;\ge\; \mathrm{opt(LP(F^S_{FP}))}$,
    \item[(ii)] the inequality in (i) is strict for some problem instances.
\end{enumerate}
\end{proposition}
\begin{proof}
(i)  Let $(X,G)$ be a feasible  solution of $\mathrm{LP(F^S_{GW})}$, and define
$z_{ij}:=\sum_{k\in\mathcal B} g_{ijk}$ $(i,j\in\mathcal N,\ i<j).$
We show that $(X,Z)$ is feasible for $\mathrm{LP(F^S_{FP})}$ and has the same objective value.
Constraints \eqref{BQP3}, \eqref{SYM1}, \eqref{SYM2}, and \eqref{CapSym}, together with the continuous relaxation of~\eqref{BQP5}, involve only $X$ and are identical in both formulations.
It remains to verify the linearization constraints~\eqref{FP1}-\eqref{FP4} from $\mathrm{LP(F^S_{FP})}$.
From constraints \eqref{GWL3} of $\mathrm{LP(F^S_{GW})}$, we have
\begin{equation}
x_{ik} + x_{jk} - 1 \le g_{ijk} \le \sum_{\ell\in\mathcal{B}} g_{ij\ell} = z_{ij} \quad i,j\in \mathcal{N}, k\in \mathcal{B}, i<j, \label{derivationFPconst}
\end{equation}
which corresponds to \eqref{FP1}.

For any $k'\in\mathcal{B}$, constraint~\eqref{GWL1} gives $g_{ijk'}\le x_{ik'}$ where $ i,j\in \mathcal{N}$ with $i<j$. Summing constraints \eqref{GWL2} over $k\in\mathcal{B}$ where $k\ne k'$ results in
$\sum_{\substack{k\in\mathcal{B}\\ k\ne k'}} g_{ijk}
\le \sum_{\substack{k\in\mathcal{B}\\ k\ne k'}} x_{jk}$. Combining these two inequalities and using $\sum_{k\in\mathcal{B}} x_{jk}=1$, we get
\begin{equation}
 \sum_{k\in\mathcal{B}} g_{ijk}
\le \sum_{\substack{k\in\mathcal{B}\\ k\ne k'}} x_{jk}+x_{ik'}
= (1 - x_{jk'})+x_{ik'} \quad\quad i,j\in \mathcal{N}, i<j.  \label{int_step}
\end{equation}
For each $k'\in\mathcal{B}$, replacing $\sum_{k\in\mathcal{B}} g_{ijk}$ with $z_{ij}$ results in constraints~\eqref{FP3}.
Interchanging $i$ and $j$ yields constraints \eqref{FP2}.
Finally, nonnegativity constraints~\eqref{FP4} follow directly from \eqref{GWL4}.
Therefore, $(X,Z)$ is a feasible solution for $\mathrm{LP(F^S_{FP})}$, and
by the definition of $z_{ij}$, both formulations yield the same objective value.

(ii) Consider Instance~\ref{inst1}.
For this instance, the optimal objective values of $\mathrm{LP(F^S_{GW})}$ and $\mathrm{LP(F^S_{FP})}$ are $4.67$ and $4$, respectively.
Hence, $\mathrm{opt}\bigl(\mathrm{LP(F^S_{GW})}\bigr) > \mathrm{opt}\bigl(\mathrm{LP(F^S_{FP})}\bigr)$, which establishes that the inequality in part~(i) can be strict.
Corresponding optimal solutions are reported in Appendix~\ref{appdx:LPrelstrong}.
\end{proof}

We now derive {valid equalities} for $\mathrm{F^S_{GW}}$ that follow from the symmetry-breaking constraints and can be used to strengthen its LP relaxation.

\begin{lemma} \label{lemma:gikk}
For $i,j \in \mathcal{N}$ with $i < j$, the equality $g_{iji} = x_{ji}$ is valid for $\mathrm{F^S_{GW}}$ and strengthens its LP relaxation.
\end{lemma}
\begin{proof}
If $x_{ji} = 0$, then \eqref{GWL2} and \eqref{GWL4} give $g_{iji}=0$.
If $x_{ji}=1$, then \eqref{SYM2} gives $x_{ii}=1$, and together with \eqref{GWL1} and \eqref{GWL3}, this implies $g_{iji}=1$.
Thus, $x_{ji} = g_{iji}$ holds for all integer feasible solutions.
Since the LP relaxation may admit fractional solutions with $g_{iji}<x_{ji}$, adding $g_{iji} = x_{ji}$ strengthens the relaxation.
\end{proof}

Below, we provide a numerical example illustrating Lemma~\ref{lemma:gikk}.
\begin{example}
\label{ex2}
Consider Instance \ref{inst1} and recall that the optimal objective value of $\mathrm{LP(F^S_{GW})}$ is $4.67$.
We denote by $\mathrm{F^{S'}_{GW}}$ the formulation obtained from $\mathrm{F^S_{GW}}$ by adding equalities from Lemma~\ref{lemma:gikk}.
The optimal objective value of $\mathrm{LP(F^{S'}_{GW})}$ is $5.33$.
The corresponding solution is provided in Appendix~\ref{appdx:ex2}.
\end{example}
We next introduce inequalities that strengthen the LP relaxations of both $\mathrm{F^S_{FP}}$ and $\mathrm{F^S_{GW}}$ by replacing the linearization constraints \eqref{FP1} and \eqref{GWL3}, respectively.

\begin{lemma}\label{lemma:stronger}
In the formulations \(\mathrm{F^S_{FP}}\) and \(\mathrm{F^S_{GW}}\), for
\(i,j \in \mathcal{N}\) and \(k \in \mathcal{B}\) with \(i<j\), the inequalities \eqref{FP1} and \eqref{GWL3} can be replaced by the valid inequalities
\begin{align}
x_{ik} + x_{jk} - x_{kk} &\le z_{ij} \label{eq:strong-z}\\
x_{ik} + x_{jk} - x_{kk} &\le g_{ijk} \label{eq:strong-g}
\end{align}
respectively, which tighten the corresponding LP relaxations.
\end{lemma}
\begin{proof}
We first show that \eqref{eq:strong-g} is valid for the integer feasible set of
$\mathrm{F^S_{GW}}$. The proof for $\mathrm{F^S_{FP}}$ is analogous.
We distinguish the five possible relative positions of $k$.

\emph{(i) $i<j<k$.} By \eqref{SYM1}, $x_{ik}=x_{jk}=0$, so the left-hand side of
\eqref{eq:strong-g} is nonpositive.

\emph{(ii) $i<j=k$.} Here $x_{jk}=x_{kk}$ and $x_{ik}=0$ by \eqref{SYM1}, so
\eqref{eq:strong-g} reduces to $g_{ijk}\ge 0$.

\emph{(iii) $i<k<j$.} By \eqref{SYM1} $x_{ik}=0$ and by \eqref{SYM2}
$x_{jk}\le x_{kk}$. Thus, the left-hand side of \eqref{eq:strong-g} is nonpositive.

\emph{(iv) $k=i<j$.} Since $x_{ik}=x_{kk}$, \eqref{eq:strong-g} reduces to
$x_{jk}\le g_{ijk}$. Moreover, $x_{jk}\le x_{kk}$ by \eqref{SYM2}. If
$x_{jk}=0$, the inequality follows from $g_{ijk}\ge 0$. If $x_{jk}=1$, then
$x_{kk}=x_{ik}=1$, and the linearization constraints imply $g_{ijk}=1$.

\emph{(v) $k<i<j$.} By \eqref{SYM2}, $x_{ik}\le x_{kk}$ and
$x_{jk}\le x_{kk}$. If $x_{kk}=0$, then $x_{ik}=x_{jk}=0$. If $x_{kk}=1$, then
\eqref{eq:strong-g} coincides with the original inequality
$x_{ik}+x_{jk}-1\le g_{ijk}$.

Thus \eqref{eq:strong-g} is valid in all cases.
To see that the replacement tightens the LP relaxation, note that $x_{kk}\le 1$
implies $x_{ik}+x_{jk}-x_{kk} \ge x_{ik}+x_{jk}-1$.
Hence, \eqref{eq:strong-g} dominates the original inequality \eqref{GWL3}.
The same argument applies to the corresponding inequality in $\mathrm{F^S_{FP}}$.
\end{proof}

We illustrate Lemma~\ref{lemma:stronger} with the following numerical example.
\begin{example}
\label{ex3}
Consider Instance~\ref{inst1}.
For this instance, the optimal objective values of $\mathrm{LP(F^S_{FP})}$ and $\mathrm{LP(F^S_{GW})}$ are $4$ and $4.67$, respectively.
After replacing \eqref{FP1} and \eqref{GWL3} by the strengthened inequalities \eqref{eq:strong-z} and \eqref{eq:strong-g}, respectively, the corresponding optimal LP-relaxation values increase to $5$ and $6.67$.
Thus, the inequalities in Lemma~\ref{lemma:stronger} strictly strengthen both LP relaxations for this instance.
Corresponding optimal solutions are reported in Appendix~\ref{appdx:ex3}.
\end{example}
\phantomsection\label{form:FP} \phantomsection\label{form:GW}
We denote by $\mathrm{F^{S+}_{FP}}$ and $\mathrm{F^{S+}_{GW}}$ the formulations obtained from $\mathrm{F^S_{FP}}$ and $\mathrm{F^S_{GW}}$ by replacing constraints~\eqref{FP1} and~\eqref{GWL3} with the inequalities~\eqref{eq:strong-z} and~\eqref{eq:strong-g}, respectively.
The reduced formulations of $\mathrm{F^{S+}_{FP}}$ and $\mathrm{F^{S+}_{GW}}$ are given in Appendices~\ref{appdx:FP} and~\ref{appdx:GW}, respectively.
These reductions exploit the symmetry-breaking constraints~\eqref{SYM1} and \eqref{SYM2} and the signs of the dissimilarity coefficients.
For $\mathrm{F^{S+}_{GW}}$, inequalities~\eqref{eq:strong-g} also permit the elimination of the variables $g_{iji}$ for all $i,j\in\mathcal N$ with $i<j$, using Lemma~\ref{lemma:gikk}.

\paragraph{The Adams--Sherali linearization approach:}
\cite{AS86} proposed the {\em reformulation linearization technique} (RLT) for BQP formulations.
The RLT provides a hierarchy of relaxations that are ranging between the continuous and convex hull representation for MILPs~\citep{SA90,SA13}. The RLT consists of the following two steps: a reformulation step in which additional nonlinear valid inequalities are generated, and a linearization step in which each product term is replaced by a continuous variable. The level of the hierarchy corresponds to the degree of polynomial terms produced during the reformulation stage.

To apply the RLT to the QBPP, we introduce binary variables $\hat{g}_{ijk} = x_{ik}x_{jk}$ for $i,j \in \mathcal{N}$ and $k \in \mathcal{B}$, indicating whether items $i$ and $j$ are assigned to bin $k$. These variables are similar to those used in the Glover--Woolsey linearization, but here all pairs $i,j \in \mathcal{N}$ are considered.
The resulting level-1 RLT formulation, with symmetry-breaking constraints, for the QBPP is
\phantomsection\label{form:RLT1}
\begin{align}
(\mathrm{F^S_{RLT1}}) \ \min \ & \frac{1}{2}\sum_{i\in\mathcal{N}} \sum_{\substack{j\in\mathcal{N}}} \sum_{k \in \mathcal{B}} d_{ij}\hat{g}_{ijk}+c\sum_{k \in \mathcal{B}}x_{kk} \label{RLobj} \\
\text{s.t.} \quad &
\eqref{BQP3}, \eqref{BQP5}-\eqref{CapSym} \notag\\
&\hat{g}_{ijk} \leq x_{ik} & i,j\in \mathcal{N}, k\in \mathcal{B} \label{R1}\\
&x_{ik} + x_{jk} - 1 \leq \hat{g}_{ijk} &i,j\in \mathcal{N}, k\in \mathcal{B} \label{R2}\\
& \hat{g}_{ijk} = \hat{g}_{jik} & i,j\in \mathcal{N}, k\in \mathcal{B}, i<j \label{R3}\\
&\sum_{i\in \mathcal{N}} w_i \hat{g}_{ijk}\leq W\, \hat{g}_{kjk} & j \in \mathcal{N}, k \in \mathcal{B} \label{RLTC1_1+} \\
&\sum_{i\in \mathcal{N}} w_i (x_{ik}-\hat{g}_{ijk}) \leq W (x_{kk} - \hat{g}_{kjk}) &   j \in \mathcal{N}, k \in \mathcal{B} \label{RLTC2_1-} \\
&\hat{g}_{ijk} \leq \hat{g}_{jkk} & i,j \in \mathcal{N}, k \in \mathcal{B}, k\!<\!i \label{RLTS_1+} \\
&x_{ik}-x_{kk} \le \hat{g}_{ijk} - \hat{g}_{jkk} & i,j\in \mathcal{N},k \in \mathcal{B}, k\!<\!i\label{RLTS_1-} \\
& x_{ik} = \hat{g}_{iik}  & i\in\mathcal{N}, k \in \mathcal{B} \label{Connector} \\
& \hat{g}_{ijk} \ge 0 & i,j\in \mathcal{N}, k\in \mathcal{B}. \label{R5}
\end{align}

Notice that formulation $\mathrm{F^S_{RLT1}}$ without constraints~\eqref{RLTC1_1+}–\eqref{Connector} reduces to $\mathrm{F^S_{GW}}$.
Constraints \eqref{GWL1}, \eqref{GWL3} and \eqref{GWL4} correspond to~\eqref{R1}, \eqref{R2}, and \eqref{R5}, respectively,
and constraints~\eqref{GWL2} follow from~\eqref{R1} and~\eqref{R3}.
Constraints \eqref{RLTC1_1+} and \eqref{RLTC2_1-} are RLT constraints derived from the capacity constraints \eqref{CapSym}, whereas \eqref{RLTS_1+} and \eqref{RLTS_1-} are RLT constraints obtained from the symmetry-breaking constraints \eqref{SYM2}.
Lastly, in constraints~\eqref{Connector}, we exploit the relation $x_{ik}^2=x_{ik}$. Although implied by the linearization constraints for integer solutions, they are included to strengthen the LP relaxation and used to organize the reprsentations of the RLT constraints in the formulation.

We denote feasible solutions of $\mathrm{F^S_{RLT1}}$ by $(X, \hat G)$, where $X = (x_{ij})$ and $\hat G = (\hat g_{ijk})$.
Similar to the $G$ variables in $\mathrm{F^{S}_{GW}}$, each $\hat g_{ijk}$ $(i,j \in \mathcal{N},\ k \in \mathcal{B})$ is implicitly binary.
This linearization introduces $\mathcal{O}(n^2 b)$ auxiliary variables and $\mathcal{O}(n^2 b)$ additional constraints.

We note that $\mathrm{F^S_{RLT1}}$ excludes the RLT constraints derived from the symmetry-breaking constraints~\eqref{SYM1} and the assignment constraints~\eqref{BQP3}. The former are omitted because the corresponding RLT constraints are already implied by constraints~\eqref{SYM1} and \eqref{R1}, while the latter are excluded as a modeling choice to restrict the formulation to variables with at most three indices. Including those RLT constraints
would indeed strengthen the LP relaxation, but at the cost of introducing four-indexed variables, since such constraints involve products across different bins.
We first identify equalities that are implied by the LP relaxation of $\mathrm{F^S_{RLT1}}$
and that will be used in the proof of Lemma~\ref{lem:redundancyR2}.
\begin{lemma}\label{lem:redundancy}
    Constraints $x_{ik}=\hat{g}_{ikk}$ for $i\in\mathcal{N}$ and $k\in\mathcal{B}$ are implied by $\mathrm{LP(F^S_{RLT1})}$.
\end{lemma}
\begin{proof}
    Consider the LP relaxation of $\mathrm{F^S_{RLT1}}$. For $j=k$, \eqref{RLTC2_1-} becomes $\sum_{i\in \mathcal{N}} w_i (x_{ik}-\hat{g}_{ikk}) \leq W (x_{kk} - \hat{g}_{kkk})$ for $k \in \mathcal{B}$, where $\hat{g}_{kkk}=x_{kk}$ by \eqref{Connector}.
    Constraints~\eqref{R1} imply $x_{ik}-\hat{g}_{ikk}\ge 0$ for $i\in\mathcal{N}$ and $k\in\mathcal{B}$.
    Since $w_i > 0$ for $i\in\mathcal{N}$, $x_{ik}-\hat{g}_{ikk}=0$ for $i\in\mathcal{N}$ and $k\in\mathcal{B}$.
\end{proof}

We next use this implication to show that the linearization constraints~\eqref{R2} are redundant for the LP relaxation of $\mathrm{F^S_{RLT1}}$.

\begin{lemma}\label{lem:redundancyR2}
    Constraints \eqref{R2} are redundant for $\mathrm{LP(F^S_{RLT1})}$.
\end{lemma}
\begin{proof}
We prove that each inequality in~\eqref{R2} is implied by the remaining constraints of $\mathrm{LP(F^S_{RLT1})}$.
We distinguish the position of $k$ relative to $i$. In the first two cases, the
conclusion follows for all $j\in\mathcal N$ without further case distinctions.

\emph{(i) $k<i$}: By \eqref{RLTS_1-} and Lemma~\ref{lem:redundancy}, $x_{ik}-x_{kk} \le \hat g_{ijk}-\hat g_{jkk} = \hat g_{ijk}-x_{jk}.$
Equivalently, $x_{ik}+x_{jk}-x_{kk}\le \hat g_{ijk}.$
Since $x_{kk}\le 1$, this implies $x_{ik}+x_{jk}-1 \le \hat{g}_{ijk}$.

\emph{(ii) $k>i$}: By \eqref{SYM1}, $x_{ik}=0$, hence \eqref{R2} reduces to $x_{jk}-1\le \hat g_{ijk}$, which follows from $x_{jk}\le 1$ and $\hat g_{ijk}\ge 0$.

\emph{(iii) $k=i$}: Then \eqref{R2} becomes $x_{ii}+x_{ji}-1\le \hat g_{iji}.$
If $j<i$, then $x_{ji}=0$ by \eqref{SYM1}, and the inequality follows from $x_{ii}\le 1$ and $\hat g_{iji}\ge 0$.
If $j=i$, then, by \eqref{Connector}, $\hat g_{iii}=x_{ii}$, and the inequality reduces to $x_{ii}\le 1$.
Finally, if $j>i$, then \eqref{R3} and Lemma~\ref{lem:redundancy} give $\hat g_{iji}=\hat g_{jii}=x_{ji}.$
Thus, the inequality reduces to $x_{ii}\le 1$.

Hence, constraints~\eqref{R2} are implied in all cases and therefore redundant in $\mathrm{LP(F^S_{RLT1})}$.
\end{proof}
Since constraints~\eqref{R2} are redundant for the LP relaxation, they can also be omitted from $\mathrm{F^{S}_{RLT1}}$.
Further reductions in formulation size are possible
by omitting the variables $x_{ik}$ with $i<k$ and $\hat g_{ijk}$ with $\min\{i,j\}<k$, whose values are necessarily zero by \eqref{SYM1}, \eqref{R1} and \eqref{R3}.
The constraints can then be imposed only over the corresponding index sets.
For completeness, Appendix~\ref{appdx:RLT} presents the formulation in which these reductions are explicitly incorporated.
Finally, since removing constraints~\eqref{RLTC1_1+}–\eqref{Connector} from $\mathrm{F^S_{RLT1}}$ yields a formulation equivalent to $\mathrm{F^S_{GW}}$,
the LP relaxation of $\mathrm{F^S_{RLT1}}$ dominates that of $\mathrm{F^S_{GW}}$.

\paragraph{The Glover linearization approach:}
For BQPs, \cite{G75} introduced a linearization approach that results in a formulation with fewer variables and constraints compared to the methods of \cite{FP10}, \cite{GW74}, and \cite{AS86}.
Following this approach, we write the quadratic part of the objective function in $\mathrm{F_{BQP}}$ as
\[
\sum_{k \in \mathcal{B}} \sum_{i=1}^{n-1} \sum_{j=i+1}^n d_{ij}x_{ik}x_{jk}
=
\sum_{k \in \mathcal{B}} \sum_{i=1}^{n-1} x_{ik} f_{ik}(x),
\]
where $f_{ik}(x)=\sum_{j=i+1}^{n} d_{ij} x_{jk}$.
We introduce the auxiliary variable $t_{ik}:= x_{ik} f_{ik}(x)$, which represents the total dissimilarity {value} between item $i$ and items $j=i+1,\dots, n$ packed in bin $k$.
Using lower and upper bounds $l_i$ and $u_i$ on $f_{ik}(x)$, which also
serve as valid bounds for $t_{ik}$,  we linearize the product $x_{ik}f_{ik}(x)$ and obtain the following MILP formulation for the QBPP:
\phantomsection\label{form:G}
\begin{align}
(\mathrm{F^S_{G}}) \quad \min \quad & {\sum_{i =1}^{n-1}} \sum_{k \in \mathcal{B}}  t_{ik} + c \sum_{k \in \mathcal{B}} x_{kk}\label{Gobj}\\
\text{s.t.} \quad &
\eqref{BQP3}, \eqref{BQP5} - \eqref{CapSym} \notag\\
& l_{i} x_{ik} \leq t_{ik} &&  i \in \mathcal{N}\setminus\!\{n\}, k \in \mathcal{B} \label{G1} \\
& \sum_{j=i+1}^{n} d_{ij} x_{jk} -u_{i} (1-x_{ik}) \leq t_{ik} &&  i \in \mathcal{N}\setminus\!\{n\}, k \in \mathcal{B} \label{G2} \\
& l_i \le t_{ik} \le u_i && i \in \mathcal{N}\setminus\!\{n\}, k \in \mathcal{B}. \label{G3}
\end{align}
The bounds $l_{i}$ and $u_{i}$ for $i\in\mathcal{N}\setminus\{n\}$ are defined by\\
\begin{minipage}{0.48\textwidth}
\begin{align*}
l_{i} := \sum_{\substack{j=i+1 \\ d_{ij} < 0}}^n d_{ij},
\end{align*}
\end{minipage}%
\hfill
\begin{minipage}{0.48\textwidth}
\begin{align*}
u_{i} := \sum_{\substack{j=i+1 \\ d_{ij} > 0}}^n d_{ij},
\end{align*}
\end{minipage}
and if the corresponding sum is empty, its value is taken to be zero.

Since $t_{ik}$ represents the product $x_{ik}f_{ik}(x)$, constraints~\eqref{G1} and \eqref{G2} ensure that $t_{ik}$ takes the correct value depending on whether item $i$ is assigned to bin $k$.
In particular, constraints~\eqref{G1} impose the lower bound $l_i$, while constraints~\eqref{G2} ensure $t_{ik}=f_{ik}(x)$ when $x_{ik}=1$ and allow $t_{ik}=0$ when $x_{ik}=0$.
The term $u_i(1-x_{ik})$ is required since $\sum_{j=i+1}^n d_{ij}$ may be positive, which would force $t_{ik}$ to remain bounded below by a positive value even when $x_{ik}=0$ (although $t_{ik}$ should then be zero).
We note that, if all $d_{ij}$ were nonpositive, the term $u_i(1-x_{ik})$ would not be necessary.
This approach introduces $\mathcal{O}(nb)$ auxiliary variables and $\mathcal{O}(nb)$ additional constraints.

\begin{lemma}
\label{lem:G_diag_lqp6}
For $i\in\mathcal{N}\setminus\{n\}$, the inequality
\begin{equation}
    \sum_{j=i+1}^{n} d_{ij}x_{ji} \le t_{ii}  \label{G_diag_lqp6_valid}
\end{equation}
is valid for all feasible solutions of $\mathrm{F^S_{G}}$.
Moreover, \eqref{G_diag_lqp6_valid} dominates \eqref{G2} with $k=i$ in the LP relaxation.
\end{lemma}

\begin{proof}
We distinguish two cases.
If $x_{ii}=1$, then constraints~\eqref{G2} with $k=i$ reduce to \eqref{G_diag_lqp6_valid}.
If $x_{ii}=0$, then by \eqref{SYM2}, we have $x_{ji}\le x_{ii}=0$ for $j>i$.
Therefore $\sum_{j=i+1}^{n} d_{ij}x_{ji}=0$.
Furthermore, constraints~\eqref{G1} with $k=i$ give $l_i x_{ii} \le t_{ii}$, and since $x_{ii}=0$, this implies $t_{ii}\ge 0$.
Thus, $\sum_{j=i+1}^{n} d_{ij}x_{ji}=0\le t_{ii},$ so~\eqref{G_diag_lqp6_valid} is valid.
Moreover, since $u_i(1-x_{ii})\ge 0$ for $0\le x_{ii}\le 1$,
\eqref{G_diag_lqp6_valid} implies \eqref{G2} with $k=i$.
Thus, \eqref{G_diag_lqp6_valid} dominates the original constraint in the LP relaxation.
\end{proof}

The following numerical example illustrates Lemma~\ref{lem:G_diag_lqp6}.
\begin{example}
Consider Instance~\ref{inst1}.
For this instance, the bounds on $f_{ik}(x)$ are $l=(-4,-20,0)$ and $u=(20,0,8)$.
Replacing constraints~\eqref{G2} for $k=i$ by the strengthened inequalities of Lemma~\ref{lem:G_diag_lqp6} increases the optimal LP-relaxation value from $3.6$ to $4.5$.
Corresponding optimal solutions are reported in Appendix~\ref{appdx:G1example}.
\end{example}

Because of the symmetry-breaking constraints~\eqref{SYM1}, the variables $x_{ik}$ and $t_{ik}$ with $i<k$ can be omitted.
Accordingly, constraints~\eqref{G1} are imposed only for $k\le i$.
Moreover, constraints~\eqref{G2} are treated separately for two cases: the original form is kept for $k<i$, while for $k=i$ it is replaced by the strengthening given in Lemma~\ref{lem:G_diag_lqp6}.

For the QKP, \cite{BS04} strengthened bounds $l_i$ by solving a series of knapsack subproblems.
Inspired by their approach, we strengthen the LP relaxation of $\mathrm{F^S_G}$ by exploiting the bin-capacity constraints to derive tighter bounds on the interaction terms. Our construction differs from that of \cite{BS04}, as the bounds are tailored to the residual capacity available in bin $k$.
Recall that the symmetry-breaking constraints \eqref{SYM1} and \eqref{SYM2} enforce that bin $k$ is used if and only if item $k$ is packed into bin $k$.
Exploiting this property, we {derive} item--bin-dependent lower and upper bounds $\tilde l_{ik}$ and $\tilde u_{ik}$ for $k\le i<n$ to be used in constraints  \eqref{G1} and \eqref{G3}, together with refined deactivation coefficients $\tilde M_{ik}$ for $k<i<n$ to be used in constraints \eqref{G2}.
No such quantities are required for $k>i$, because the corresponding assignment variables are eliminated due to symmetry breaking constraints.

Regarding $\tilde l_{ik}$ and $\tilde u_{ik}$ for $k\le i < n$ in constraints \eqref{G1} and \eqref{G3}, let us first observe that any values satisfying $\tilde l_{ik}\le 0\le\tilde u_{ik}$ are valid when $x_{ik}=0$, since, in that case, $t_{ik}$ is bounded below by 0 through \eqref{G1}.
We therefore focus on the case where $x_{ik}=1$, i.e., item $i$ is assigned to bin $k$.
When $k<i$, both items $i$ and $k$ must belong to bin~$k$, leaving residual capacity $W-w_i-w_k$ for items $j>i$.
When $k=i$, only item $i$ is fixed in the bin, leaving residual capacity $W-w_i$.
Using these residual capacities, we derive the following lower and upper bounds on $f_{ik}(x)$.
\vspace{-0.5cm}

\begin{adjustbox}{max width=\textwidth}
\begin{minipage}[t]{0.48\textwidth}
\begin{align*}
\tilde l_{ik}:=\min \quad
& \sum_{\substack{t=1\\ d_{i,i+t}<0}}^{n-i} d_{i,i+t}x_t \\
\text{s.t.}\quad
& \sum_{t=1}^{n-i} w_{i+t}x_t
\le
\begin{cases}
W-w_i-w_k, & k<i\\
W-w_i, & k=i
\end{cases}\\
& x_t\in\{0,1\}, \quad t=1,\ldots,n-i.
\end{align*}
\end{minipage}
\hfill
\begin{minipage}[t]{0.48\textwidth}
\begin{align*}
\tilde u_{ik}:=\max \quad
& \sum_{\substack{t=1\\ d_{i,i+t}>0}}^{n-i}
d_{i,i+t}x_t \\
\text{s.t.}\quad
& \sum_{t=1}^{n-i} w_{i+t}x_t
\le
\begin{cases}
W-w_i-w_k, & k<i,\\
W-w_i, & k=i,
\end{cases}\\
& x_t\in\{0,1\},
\quad t=1,\ldots,n-i.
\end{align*}
\end{minipage}
\end{adjustbox}
\vspace{1mm}

However, the bound $\tilde u_{ik}$ for $k<i<n$ is not necessarily valid in constraints \eqref{G2}, since it may be exceeded by the quantity $\sum_{j=i+1}^{n} d_{ij} x_{jk}$ when $x_{ik}=0$.
This is because, in this case, item $i$ does not occupy any capacity in bin $k$, while other items $j>i$ may still be assigned to that bin.
Therefore, a separate bound $\tilde M_{ik}$ is needed for these constraints.
Using reasoning similar to that described above, let us first observe that any
$0 \le \tilde M_{ik}$ value is valid when $x_{ik} = 1$ or when bin~$k$ is inactive.
We therefore focus on the case where bin $k$ is active and $x_{ik} = 0$, i.e., item $i$ is not assigned to bin $k$.
Since bin $k$ is active, it must contain item $k$, leaving at most $W-w_k$ capacity for items $j > i$.
Using this residual capacity, we derive the following upper bound for  $\sum_{j=i+1}^{n} d_{ij} x_{jk}$.
\begin{align*}
\tilde M_{ik}:=\max \quad
& \sum_{\substack{t=1\\ d_{i,i+t}>0}}^{n-i}
d_{i,i+t}x_t \\
\text{s.t.}\quad
& \sum_{t=1}^{n-i} w_{i+t}x_t\le W-w_k,\\
& x_t\in\{0,1\},
\quad t=1,\ldots,n-i.
\end{align*} \vspace{-0.5cm}

We point out that no coefficient $\tilde M_{ii}$ is needed because constraints~\eqref{G2} with $k=i$ are replaced by the strengthened inequalities~\eqref{G_diag_lqp6_valid}.
Throughout, we use the convention that if the objective function in any of these problems is an empty sum, its value is taken to be zero.
Although this approach requires solving a series of ($\mathcal{NP}$-hard) knapsack problems, these are empirically fast to solve for reasonably sized instances using dynamic programming.

\phantomsection\label{form:G2}

We denote by $\mathrm{F^{S+}_{G}}$ the formulation obtained from $\mathrm{F^S_{G}}$ by using the refined variable bounds $\tilde l_{ik}$ and $\tilde u_{ik}$, the deactivation coefficients $\tilde M_{ik}$, and replacing constraints~\eqref{G2} for $k=i$ with the strengthened inequalities \eqref{G_diag_lqp6_valid}.
We include the resulting reduced formulation of $\mathrm{F^{S+}_{G}}$ in Appendix~\ref{appdx:G}.
\begin{proposition}\label{prop:RLTvsG}
The relaxations $\mathrm{LP}(\mathrm{F^S_{RLT1}})$ and $\mathrm{LP}(\mathrm{F^{S+}_{G}})$
do not dominate each other.
\end{proposition}
\begin{proof}
Consider Instance~\ref{inst1} and~\ref{inst2}.
For Instance~\ref{inst1}, the optimal objective values of $\mathrm{LP}(\mathrm{F^S_{RLT1}})$ and $\mathrm{LP}(\mathrm{F^{S+}_{G}})$ are $12$ and $8.57$, respectively. Hence, $\mathrm{opt(LP(F^S_{RLT1}))}>$ $\mathrm{opt(LP(F^{S+}_{G}))}$.
For Instance~\ref{inst2}, the corresponding optimal objective values are $6.86$ and $6.91$, respectively.
Hence, $\mathrm{opt(LP(F^S_{RLT1}))}<$ $\mathrm{opt(LP(F^{S+}_{G}))}$.
Therefore, neither relaxation dominates the other. Corresponding optimal solutions are reported in Appendix~\ref{appdx:RLTvsG}.
\end{proof}

\paragraph{Pattern-based and item-partitioned hybrid formulations.}
We next consider a pattern-based {equivalent formulation} of the QBPP, in which the entire composition of each candidate bin is specified in advance.
This formulation is closely related to the classical set-covering formulation for the BPP, based on the pattern-based approach of \citet{GG61}, in which each feasible subset of items defines a bin pattern.
Here, we use a set-partitioning version in which the selected patterns form a partition of the item set. Let $\mathcal P$ denote the set of all nonempty feasible patterns: $\mathcal{P} = \left\{ p\subseteq\mathcal N: \sum_{i\in p}w_i\le W \right\}$.
For each $p\in\mathcal P$, let $\delta_{ip}=1$ if $i\in p$, and
$\delta_{ip}=0$ otherwise. Since the items contained in a pattern are known,
its cost can be computed in advance as $
\rho_p  = c+\sum_{i=1}^{n-1}\sum_{j=i+1}^{n} d_{ij}\delta_{ip}\delta_{jp}.
$
Introducing a binary variable $\lambda_p$ indicating whether pattern $p$ is
selected gives the following set-partitioning formulation of the QBPP:
\phantomsection\label{form:P}
\begin{align}
(\mathrm{F_P})\quad \min \quad
& \sum_{p\in\mathcal P}\rho_p\lambda_p \label{P_obj}\\
\text{s.t.}\quad
& \sum_{p\in\mathcal P}\delta_{ip}\lambda_p=1
&& i\in\mathcal N \label{P_cover}\\
& \lambda_p\in\{0,1\}
&& p\in\mathcal P. \label{P_bin}
\end{align}
Constraints~\eqref{P_cover} require each item to be contained in exactly one selected pattern.
When all feasible patterns can be enumerated, $\mathrm{F_P}$ offers several advantages.
Each variable represents a complete feasible bin, so the capacity requirement and all pairwise interaction costs are encoded directly in the pattern.
Consequently, the formulation does not require bin-activation, assignment or pairing variables and eliminates the symmetry caused by interchangeable bin indices.
Moreover, it contains only $\mathcal O(n)$ set-partitioning constraints and can be solved directly using a standard MILP solver.
Its LP relaxation coincides with the root relaxation of the B\&P formulation of \cite{CLMI26}.
For the instances for which $\mathrm{F_P}$ could be constructed and its LP relaxation solved, the resulting lower bound was never weaker than those of the alternative formulations and was strictly stronger in 95.8\% of cases.
The main drawback is the potentially exponential number of pattern variables. In the worst case, every nonempty subset of $\mathcal N$ is feasible, yielding $|\mathcal P|=2^n-1$. For example, when $n=25$, the number of feasible patterns can potentially exceed $33$ million.
When complete pattern enumeration results in too many variables to be handled directly by a MILP solver, the formulation can instead be solved within a B\&P framework, as done by \cite{CLMI26}.

Motivated by this limitation, we develop an item-partitioned hybrid formulation that combines two complementary modeling structures.
Although this particular hybridization is new for the QBPP, hybrid formulations that combine different modeling structures have also been proposed in other combinatorial settings (see, e.g., \cite{BDMM26,DMDPST16}).
To retain some of the advantages of explicit pattern enumeration while
controlling the size of the pattern set, we select a subset $\mathcal{A}\subseteq\mathcal{N}$ for pattern generation and only generate patterns corresponding to feasible subsets of $\mathcal{A}$.
We combine these patterns with the strengthened Glover--Woolsey formulation $\mathrm{F^{S+}_{GW}}$, which constitutes the \emph{assignment-based component} of the resulting hybrid formulation.
Restricting the patterns to subsets of $\mathcal{A}$ results in at most $2^{|\mathcal{A}|}-1$ patterns.
Items in $\mathcal{A}$ can be represented either by a pattern or through the assignment-based component, while items in $\mathcal {N}\setminus\mathcal{A}$ are represented exclusively through the assignment-based component.
Allowing both representations for items in $\mathcal{A}$ is necessary, since otherwise they could not share a bin with items outside $\mathcal{A}$.
This construction controls the size of both components: the pattern component contains at most $2^{|\mathcal A|}-1$ variables, while the assignment-based component requires at most $|\mathcal N\setminus\mathcal A|$ bin indices.
Thus, enlarging $|\mathcal{A}|$ increases the potential size of the pattern set while reducing the number of bin-indexed assignment and linearization variables.

Let $\bar{\mathcal{A}}=\mathcal N\setminus\mathcal{A}$ denote the set of items that must be represented through the assignment-based component.
We reorder the items such that the items in $\bar{\mathcal{A}}$ appear first and those in $\mathcal{A}$ afterwards.
Thus, letting $m=|\bar{\mathcal{A}}|$, we have $\bar{\mathcal{A}}=\{1,\ldots,m\}$, $\mathcal{A}=\{m+1,\ldots,n\}.$
Using this new ordering and the symmetry-breaking constraints \eqref{SYM1} and \eqref{SYM2}, bin $k$ is used if and only if $x_{kk}=1$ for some $k\in\bar{\mathcal{A}}$.
As in $\mathrm{F^S_{GW}}$, separate bin-activation variables $y_k$ are therefore unnecessary.
We next define the set of feasible patterns over $\mathcal{A}$ as $\mathcal P(\mathcal A)= \left\{p\subseteq \mathcal{A}:\allowbreak p\neq\emptyset,\; \sum_{i\in p} w_i\le W \right\}$. For each $p\in\mathcal
{P}(\mathcal{A})$, the parameters $\delta_{ip}$ and $\rho_p$, and the binary variable $\lambda_p$, are defined as in the pattern-based formulation above.
The resulting hybrid formulation of the QBPP is
\phantomsection\label{form:HI}
\begin{align}
\hspace{-0.6cm}
(\mathrm{F^{S+}_{HI}}) \ \min\
& \sum_{p\in\mathcal P(\mathcal A)}\rho_p\lambda_p
+ c\sum_{k\in\bar{\mathcal A}}x_{kk} + \sum_{i=1}^{n-1}\sum_{j=i+1}^{n}
\sum_{k\in\bar{\mathcal A}}
d_{ij}g_{ijk}  \phantom{\qquad\qquad\qquad\qquad\quad}
\label{PGW_obj}
\end{align}
\vspace{-0.9\baselineskip}
\begin{align}
\text{s.t.}\
& \sum_{p\in\mathcal P(\mathcal A)} \delta_{ip}\lambda_p + \sum_{k\in\bar{\mathcal A}} x_{ik} =1 && i\in\mathcal N \label{PGW_assign}\\
& x_{ik}=0 && i\in\mathcal N,\ k\in\bar{\mathcal A},\ i<k \label{PGW_SYM1}\\
& x_{ik}\le x_{kk} && i\in\mathcal N,\ k\in\bar{\mathcal A},\ k<i \label{PGW_SYM2}\\
& \sum_{i\in\mathcal{N}}w_i x_{ik} \le W x_{kk} && k\in\bar{\mathcal A} \label{PGW_cap}\\
& x_{ik}+x_{jk}-x_{kk} \le g_{ijk} && i,j\in\mathcal N,\ i<j,\ k\in\bar{\mathcal A} \label{PGW_gpos}\\
& g_{ijk}\le x_{ik} && i,j\in\mathcal N,\ i<j,\ k\in\bar{\mathcal A} \label{PGW_gneg1}\\
& g_{ijk}\le x_{jk} && i,j\in\mathcal N,\ i<j,\ k\in\bar{\mathcal A} \label{PGW_gneg2}\\
& g_{ijk}\ge 0 && i,j\in\mathcal N,\ i<j,\ k\in\bar{\mathcal A} \label{PGW_gbounds}\\
& x_{ik}\in\{0,1\} && i\in\mathcal N,\ k\in\bar{\mathcal A} \label{PGW_xbin}\\
& \lambda_p\in\{0,1\} &&p\in\mathcal P(\mathcal A). \label{PGW_lambin}
\end{align}

The objective function~\eqref{PGW_obj} combines the costs of the two representations: $\rho_p$ accounts for the costs of each pattern, while the remaining terms account for the fixed and pairwise interaction costs of bins represented by the assignment-based component.
Constraints~\eqref{PGW_assign} couple the two components by requiring each item to be represented exactly once either through a selected pattern or through an assignment variable.
The remaining constraints are analogous to those in $\mathrm{F^{S+}_{GW}}$, but restricted to items in $\bar{\mathcal A}$.
In particular, constraints~\eqref{PGW_SYM1}, \eqref{PGW_SYM2}, \eqref{PGW_cap} and \eqref{PGW_gpos} are the counterparts of \eqref{SYM1}, \eqref{SYM2}, \eqref{CapSym} and \eqref{eq:strong-g}, respectively.
Similarly, constraints~\eqref{PGW_gneg1}, \eqref{PGW_gneg2}, and \eqref{PGW_gbounds} correspond to the linearization constraints~\eqref{GWL1}, \eqref{GWL2} and \eqref{GWL4}, respectively, over the reduced set of bin indices.
The hybrid formulation includes the two underlying representations as limiting cases.
If $\mathcal{A}=\emptyset$, then no patterns are generated and the formulation reduces to $\mathrm{F^{S+}_{GW}}$.
If $\mathcal A=\mathcal{N}$, then all items are represented by patterns and the assignment-based component is unnecessary, yielding $\mathrm{F_P}$.

\section{Bin-count inequalities}\label{sec:valid}

We next show how the LP relaxations of the QBPP formulations developed in the previous section can be strengthened by exploiting lower bounds inherited from the classical BPP.
Every feasible QBPP solution must satisfy the same capacity and assignment constraints as the corresponding BPP.
Therefore, any lower bound on the number of bins required for the BPP is also valid for the QBPP and can be used to strengthen its respective LP relaxations.
In particular, consider the BPP defined on the same set of items $\mathcal{N}$ and bin capacity $W$, but without the pairwise interaction terms. Its objective is to minimize the number of bins required to pack all items.
Let $L_n$ denote the optimal number of bins for the corresponding BPP.
Hence, every feasible QBPP solution must use at least $L_n$ bins.
For formulations in which $x_{kk}$ indicates whether bin $k$ is active, the lower bound $L_n$ yields what we call the \emph{global bin-count} inequality:
\begin{equation}
    \sum_{k\in\bar{\mathcal{B}}} x_{kk}\ge L_n. \label{eq:bin_count}
\end{equation}
The symmetry-breaking constraints allow the global bin-count inequality~\eqref{eq:bin_count} to be strengthened through a family of cumulative bin-count inequalities.
Let $L_m$ denote the minimum number of bins required to pack the first $m$ items.
By~\eqref{SYM1}, item $i$ can only be assigned to bins with index $k\le i$.
Hence, the first $m$ items can only be assigned to bins $1,\ldots,m$.
Consequently, we obtain the following \emph{cumulative bin-count} inequalities:
\begin{equation}
    \sum_{k=1}^m x_{kk}\ge L_m \qquad m\in\mathcal{N}. \label{eq:cumulative_bin_count}
\end{equation}
For $m=n$, the cumulative bin-count inequality reduces to the global bin-count inequality~\eqref{eq:bin_count}.
It is not necessary to add~\eqref{eq:cumulative_bin_count} for every value of $m$.
If $L_m=L_{m-1}$, adding item $m$ does not increase the minimum number of required bins.
Hence, the cumulative inequality for $m$ is already enforced by the one for $m-1$ and can be omitted.
We therefore only need to add an inequality when $L_m>L_{m-1}$.
Although simpler lower bounds on $L_m$ could be used (see \cite{DIM16}), computing the exact values by solving the corresponding sequence of BPPs to optimality requires negligible computational effort in our experiments.

The same idea can be adapted to the item-partitioned hybrid formulation $\mathrm{F^{S+}_{HI}}$.
Since each used bin is represented either by the assignment-based component or by a selected pattern, the corresponding global bin-count inequality is
\begin{equation}
    \sum_{k\in\bar{\mathcal{A}}}x_{kk} + \sum_{p\in\mathcal P(\mathcal A)}\lambda_p \ge L_n. \label{eq:HI_bin_count}
\end{equation}
By construction, the items in $\bar{\mathcal A}$ occupy the first $|\bar{\mathcal A}|$ positions in the item ordering.
The ordering of the assignment-based component yields the following cumulative bin-count inequalities:
\begin{equation}
    \sum_{k=1}^m x_{kk}\ge L_m  \qquad m=1,\ldots,|\bar{\mathcal{A}}|, \label{eq:HI_cumulative_bin_count}
\end{equation}
where $L_m$ denotes the minimum number of bins required to pack the first $m$ items.
As above, it is sufficient to include these inequalities only for values of $m$ at which $L_m$ increases.
The global and cumulative bin-count inequalities also apply to the other formulations that use the same symmetry-breaking structure, including $\mathrm{F^{S+}_{FP}}$, $\mathrm{F^{S+}_{GW}}$, $\mathrm{F^S_{RLT1}}$ and $\mathrm{F^{S+}_{G}}$ formulations.
For the pattern-based formulation $\mathrm{F_P}$, the corresponding global bin-count inequality is $\sum_{p\in\mathcal P}\lambda_p\ge L_n$.

\section{Computational experiments} \label{sec:compEx} \label{sec:comp_exp}

In this section, we evaluate the effectiveness of the proposed models through extensive computational experiments on the benchmark instances introduced by \cite{CLMI26}.
The benchmark set consists of 675 instances characterized by four parameters: the number of items $n\in\{25,30,35,40,45\}$, a parameter $\mu\in\{0.6,1,2\}$ controlling the relative importance of bin usage and pairwise dissimilarity values, the density $\delta\in\{0.25,0.5,0.75\}$ of the dissimilarity matrix $D$ and $\sigma\in\{-,+,\pm\}$ representing the sign structure {of the dissimilarity values}. Five instances are generated for each parameter configuration.
The parameter $\mu$ affects both the bin capacity and the bin usage cost and therefore controls the trade-off between using fewer bins and grouping similar items. Larger values of $\mu$ decrease the bin capacity $W$ and increase the bin usage cost, favoring solutions with fewer bins despite the more restrictive capacity.
Conversely, smaller values of $\mu$ increase the bin capacity and decrease the bin usage cost, placing emphasis on the dissimilarity component.
For each fixed combination of $(n,\delta,\sigma)$, the instances corresponding to the three values of $\mu$ share the same item weights and pairwise dissimilarity values and differ only in the resulting bin capacity and bin cost.

Each of the tested approaches was implemented in {\tt C++}, and the code is available at \url{https://github.com/fatih-burak/QBPP}.
All computational experiments were performed using a single thread on a virtual machine with an {\tt AMD EPYC-Milan Processor} running at 2.25 GHz and 360 GiB of RAM under {\tt Ubuntu 24.04}.
The models were solved using \texttt{Gurobi 13.0.2} with the default parameter settings, unless stated otherwise. A time limit of 3600 seconds was imposed for each instance.
The formulations considered in the computational experiments are listed in the summary Table~\ref{tab:model_summary}.
Whenever a formulation is used computationally, we employ its reduced version incorporating all reductions described for that formulation.
The complete formulations, including these reductions, are provided in the Appendix~\ref{appendix:reduced_formulations}.
\begin{table*}[!htbp]
\centering
\caption{Summary of the formulations considered in the computational experiments.}
\label{tab:model_summary}
\scriptsize
\renewcommand{\arraystretch}{1.4}
\begin{adjustbox}{max width=\textwidth}
\begin{tabular}{llcl}
\toprule
\textbf{Formulation}& \textbf{Name}& \textbf{\makecell{Symmetry\\breaking}}& \textbf{Strengthening and refinement} \\
\midrule
\midrule
\midrule
\midrule
\midrule
\midrule
\midrule
\midrule
$\mathrm{F_{EFGW}}$ & \makecell[l]{Extended Fortet--Glover--Woolsey\\ (best MILP-formulation by \cite{CLMI26})} & \checkmark & \makecell[l]{$\mathrm{F^S_{RLT1}}$ without~\eqref{RLTS_1+} and~\eqref{RLTS_1-} \\[4pt] (or $\mathrm{F^{S}_{GW}}$ with~\eqref{RLTC1_1+}, \eqref{RLTC2_1-} and~\eqref{Connector})}\\
\bottomrule
\end{tabular}
\end{adjustbox}
\end{table*}

\subsection{Non-hybrid formulations}

We first evaluate the performance of the non-hybrid formulations without adding any valid inequalities introduced in Section~\ref{sec:valid}.
We report the results in Table~\ref{tab:results_nonhybrid_formulations}, where the instances are aggregated over the values of $\sigma$ and $\delta$.
The first three columns report the number of items $n$, the value of $\mu$ and the number of instances in each group, respectively.
For each formulation, we report the number of instances solved to proven optimality (column ``opt'') and the average CPU time in seconds (column ``T(s)''), including the runs terminated due to the time or memory limit (3600 seconds were used in both cases).
For comparison, we also include a reimplementation of EFGW, the best-performing MILP formulation reported by \cite{CLMI26}.
EFGW corresponds to $\mathrm{F^{S}_{GW}}$ augmented with the level-1 RLT constraints derived from the capacity constraints~\eqref{RLTC1_1+} and~\eqref{RLTC2_1-}.
Since \cite{CLMI26} used IBM CPLEX, we reimplemented EFGW and solved it with Gurobi under the same computational conditions as our other formulations.
For completeness, the formulation is provided in Appendix~\ref{appdx:EFGW}.
For the pattern-based formulation $\mathrm{F_P}$, we additionally set
\texttt{Presolve=0} and \texttt{Sifting=0}.
In our preliminary experiments, Gurobi's presolve often spent substantial time processing these large models while eliminating few or no rows and columns, after which sifting introduced additional root-processing time.
We therefore disable both procedures for $\mathrm{F_P}$.

For each group, among the formulations proposed in this paper, the formulation solving the largest number of instances is highlighted in yellow.
Among the formulations for which an LP-relaxation value is reported, the largest and second-largest LP bounds are highlighted in green and blue, respectively.
For $\mathrm{F_{P}}$, no value is reported for $z_{LP}$ in groups containing at least one instance for which the formulation could not be constructed due to the memory limit.

\begin{table*}[!htbp]
  \centering
  \caption{Computational performance of the formulations, aggregated over $\sigma$ and $\delta$}
  \label{tab:results_nonhybrid_formulations}
  \scriptsize
  \setlength{\tabcolsep}{1.8pt}
  \resizebox{\textwidth}{!}{%
    \begin{tabular}{l@{}cc @{\hspace{5pt}}| rrr @{\hspace{7pt}}| rr @{\hspace{7pt}} rrr @{\hspace{7pt}} rrr @{\hspace{7pt}} rrr @{\hspace{7pt}} rrr @{\hspace{7pt}} rrr}
    \toprule
    & & &
    \multicolumn{3}{c|}{EFGW} &
    \multicolumn{2}{c}{$\mathrm{F^{S}_{BQP}}$} &
    \multicolumn{3}{c}{$\mathrm{F^{S+}_{FP}}$} &
    \multicolumn{3}{c}{$\mathrm{F^{S+}_{GW}}$} &
    \multicolumn{3}{c}{$\mathrm{F^{S}_{RLT1}}$} &
    \multicolumn{3}{c}{$\mathrm{F^{S+}_{G}}$} &
    \multicolumn{3}{c}{$\mathrm{F_{P}}$} \\
    \cmidrule(lr){4-6} \cmidrule(lr){7-8} \cmidrule(lr){9-11} \cmidrule(lr){12-14} \cmidrule(lr){15-17} \cmidrule(lr){18-20} \cmidrule(lr){21-23}
    $n$ & $\mu$ & inst &
    opt & \multicolumn{1}{c}{T(s)} & \multicolumn{1}{c|}{$z_{\mathrm{LP}}$} &
    opt & \multicolumn{1}{c}{T(s)} &
    opt & \multicolumn{1}{c}{T(s)} & \multicolumn{1}{c}{$z_{\mathrm{LP}}$} &
    opt & \multicolumn{1}{c}{T(s)} & \multicolumn{1}{c}{$z_{\mathrm{LP}}$} &
    opt & \multicolumn{1}{c}{T(s)} & \multicolumn{1}{c}{$z_{\mathrm{LP}}$} &
    opt & \multicolumn{1}{c}{T(s)} & \multicolumn{1}{c}{$z_{\mathrm{LP}}$} &
    opt & \multicolumn{1}{c}{T(s)} & \multicolumn{1}{c}{$z_{\mathrm{LP}}$} \\
    \midrule
    \num{25} & \num{0.6} & \num{45} & \num{44} & \num{281.1} & \num{-2478.1} & \num{44} & \num{256.6} & \num{40} & \num{737.4} & \num{-2457.0} & \cellcolor{resultyellow}\num{45} & \num{182.2} & \num{-2398.7} & \cellcolor{resultyellow}\num{45} & \num{158.0} & \cellcolor{resultblue}\num{-1807.1} & \num{43} & \num{450.6} & \num{-2291.8} & \cellcolor{resultyellow}\num{45} & \num{58.2} & \cellcolor{resultgreen}\num{-1234.3} \\
    \num{30} & \num{0.6} & \num{45} & \num{36} & \num{862.2} & \num{-3766.2} & \num{37} & \num{788.8} & \num{29} & \num{1624.1} & \num{-3730.3} & \num{37} & \num{778.4} & \num{-3662.8} & \cellcolor{resultyellow}\num{42} & \num{832.2} & \cellcolor{resultgreen}\num{-2890.7} & \num{32} & \num{1275.3} & \num{-3546.5} & \num{2} & \num{3509.2} & -- \\
    \num{35} & \num{0.6} & \num{45} & \num{34} & \num{1367.9} & \num{-5345.4} & \cellcolor{resultyellow}\num{35} & \num{1013.6} & \num{19} & \num{2303.0} & \num{-5311.0} & \cellcolor{resultyellow}\num{35} & \num{1029.1} & \num{-5234.4} & \num{29} & \num{1918.6} & \cellcolor{resultgreen}\num{-4177.0} & \num{20} & \num{2204.6} & \num{-5034.9} & \num{0} & \num{3600} & -- \\
    \num{40} & \num{0.6} & \num{45} & \num{19} & \num{2527.2} & \num{-7208.8} & \cellcolor{resultyellow}\num{28} & \num{1939.1} & \num{12} & \num{2983.7} & \num{-7159.6} & \num{25} & \num{2111.1} & \num{-7081.7} & \num{12} & \num{2921.1} & \cellcolor{resultgreen}\num{-5716.1} & \num{14} & \num{2827.4} & \num{-6815.2} & \num{0} & \num{3600} & -- \\
    \num{45} & \num{0.6} & \num{45} & \num{11} & \num{2924.2} & \num{-9381.8} & \cellcolor{resultyellow}\num{16} & \num{2633.2} & \num{5} & \num{3215.6} & \num{-9341.4} & \num{15} & \num{2688.4} & \num{-9241.0} & \num{8} & \num{3330.7} & \cellcolor{resultgreen}\num{-7472.4} & \num{6} & \num{3211.1} & \num{-8840.1} & \num{0} & \num{3600} & -- \\
    \midrule avg. & \num{0.6} & \num{225} & \num{144} & \num{1592.5} & \num{-5636.0} & \cellcolor{resultyellow}\num{160} & \num{1326.3} & \num{105} & \num{2172.8} & \num{-5599.9} & \num{157} & \num{1357.8} & \num{-5523.7} & \num{136} & \num{1832.1} & \cellcolor{resultblue}\num{-4412.7} & \num{115} & \num{1993.8} & \num{-5305.7} & \num{47} & \num{2873.5} & \cellcolor{resultgreen}\num{-1183.6} \\
    \midrule
    \num{25} & \num{1} & \num{45} & \num{38} & \num{1133.5} & \num{-611.1} & \num{35} & \num{1093.5} & \num{30} & \num{1719.1} & \num{-767.2} & \num{42} & \num{754.1} & \num{-567.3} & \num{39} & \num{995.1} & \cellcolor{resultblue}\num{361.0} & \num{33} & \num{1306.9} & \num{-154.4} & \cellcolor{resultyellow}\num{45} & \num{73.5} & \cellcolor{resultgreen}\num{678.7} \\
    \num{30} & \num{1} & \num{45} & \num{20} & \num{2285.7} & \num{-1776.7} & \num{24} & \num{2033.4} & \num{13} & \num{2779.3} & \num{-1969.9} & \num{25} & \num{1938.7} & \num{-1715.7} & \num{20} & \num{2297.5} & \cellcolor{resultblue}\num{-330.6} & \num{15} & \num{2609.1} & \num{-1180.1} & \cellcolor{resultyellow}\num{35} & \num{1103.4} & \cellcolor{resultgreen}\num{234.5} \\
    \num{35} & \num{1} & \num{45} & \num{15} & \num{2695.9} & \num{-3224.7} & \num{19} & \num{2444.4} & \num{9} & \num{3062.1} & \num{-3456.5} & \num{19} & \num{2414.3} & \num{-3153.3} & \num{12} & \num{2785.7} & \num{-1210.8} & \num{8} & \num{3062.2} & \num{-2381.6} & \cellcolor{resultyellow}\num{21} & \num{2354.7} & -- \\
    \num{40} & \num{1} & \num{45} & \num{5} & \num{3237.6} & \num{-4944.3} & \cellcolor{resultyellow}\num{8} & \num{3137.9} & \num{3} & \num{3483.4} & \num{-5201.4} & \num{7} & \num{3128.3} & \num{-4862.3} & \num{5} & \num{3257.3} & \num{-2248.9} & \num{6} & \num{3476.3} & \num{-3845.7} & \num{0} & \num{3600} & -- \\
    \num{45} & \num{1} & \num{45} & \num{5} & \num{3434.3} & \num{-6982.1} & \num{5} & \num{3315.4} & \num{0} & \num{3600} & \num{-7286.1} & \cellcolor{resultyellow}\num{6} & \num{3315.3} & \num{-6888.4} & \num{3} & \num{3561.1} & \num{-3495.0} & \num{0} & \num{3599.6} & \num{-5517.6} & \num{0} & \num{3600} & -- \\
    \midrule avg. & \num{1} & \num{225} & \num{83} & \num{2557.4} & \num{-3507.8} & \num{91} & \num{2404.9} & \num{55} & \num{2928.8} & \num{-3736.2} & \cellcolor{resultyellow}\num{99} & \num{2310.2} & \num{-3437.4} & \num{79} & \num{2579.3} & \cellcolor{resultblue}\num{-1384.8} & \num{62} & \num{2810.8} & \num{-2615.9} & \num{101} & \num{2146.3} & \cellcolor{resultgreen}\num{241.6} \\
    \midrule
    \num{25} & \num{2} & \num{45} & \num{45} & \num{85.2} & \num{22033.7} & \num{43} & \num{177.5} & \num{43} & \num{198.9} & \num{21675.8} & \num{43} & \num{185.0} & \num{22053.8} & \num{43} & \num{219.8} & \cellcolor{resultblue}\num{23074.5} & \num{43} & \num{176.7} & \num{22891.7} & \num{43} & \num{164.6} & \cellcolor{resultgreen}\num{24329.9} \\
    \num{30} & \num{2} & \num{45} & \num{44} & \num{301.1} & \num{22025.6} & \num{43} & \num{347.8} & \num{33} & \num{1260.3} & \num{21483.9} & \cellcolor{resultyellow}\num{44} & \num{147.5} & \num{22057.4} & \cellcolor{resultyellow}\num{44} & \num{369.9} & \cellcolor{resultblue}\num{23461.7} & \num{43} & \num{307.7} & \num{23091.0} & \cellcolor{resultyellow}\num{44} & \num{86.3} & \cellcolor{resultgreen}\num{23922.7} \\
    \num{35} & \num{2} & \num{45} & \num{17} & \num{2519.4} & \num{21743.6} & \num{22} & \num{2178.6} & \num{9} & \num{3026.5} & \num{21041.8} & \num{26} & \num{1794.0} & \num{21785.3} & \num{19} & \num{2466.5} & \cellcolor{resultblue}\num{23852.8} & \num{16} & \num{2592.7} & \num{23190.5} & \cellcolor{resultyellow}\num{42} & \num{329.3} & \cellcolor{resultgreen}\num{24338.1} \\
    \num{40} & \num{2} & \num{45} & \num{7} & \num{3125.7} & \num{21212.0} & \num{11} & \num{2969.6} & \num{5} & \num{3275.2} & \num{20348.0} & \num{11} & \num{2829.9} & \num{21262.9} & \num{9} & \num{3032.1} & \cellcolor{resultblue}\num{24227.7} & \num{8} & \num{3031.7} & \num{23161.0} & \cellcolor{resultyellow}\num{43} & \num{361.2} & \cellcolor{resultgreen}\num{24826.7} \\
    \num{45} & \num{2} & \num{45} & \num{3} & \num{3491.0} & \num{20146.6} & \num{5} & \num{3330.4} & \num{4} & \num{3308.0} & \num{19137.4} & \num{9} & \num{2999.9} & \num{20208.3} & \num{4} & \num{3417.3} & \cellcolor{resultblue}\num{24304.7} & \num{3} & \num{3460.0} & \num{22765.9} & \cellcolor{resultyellow}\num{35} & \num{1568.6} & \cellcolor{resultgreen}\num{25047.3} \\
    \midrule avg. & \num{2} & \num{225} & \num{116} & \num{1904.5} & \num{21432.3} & \num{124} & \num{1800.8} & \num{94} & \num{2213.8} & \num{20737.4} & \num{133} & \num{1591.3} & \num{21473.5} & \num{119} & \num{1901.1} & \cellcolor{resultblue}\num{23784.3} & \num{113} & \num{1913.8} & \num{23020.0} & \cellcolor{resultyellow}\num{207} & \num{502.0} & \cellcolor{resultgreen}\num{24492.9} \\
    \midrule\midrule
    total & & \num{675} & \num{343} & \num{2018.2} & \num{4096.2} & \num{375} & \num{1844.0} & \num{254} & \num{2438.4} & \num{3800.4} & \cellcolor{resultyellow}\num{389} & \num{1753.1} & \num{4170.8} & \num{334} & \num{2104.2} & \cellcolor{resultblue}\num{5995.6} & \num{290} & \num{2239.5} & \num{5032.8} & \num{355} & \num{1840.6} & \cellcolor{resultgreen}\num{13549.2} \\
    \bottomrule
  \end{tabular}%
}
\end{table*}

Overall, among the non-hybrid formulations proposed in this paper, $\mathrm{F^{S+}_{GW}}$ achieved the best performance by solving 389 of the 675 instances to proven optimality.
Among the other formulations, interestingly, $\mathrm{F^{S}_{BQP}}$ is also competitive, solving 375 instances overall and ranking second among the proposed non-hybrid formulations.
In fact, for $\mu=0.6$, it solves slightly more instances than $\mathrm{F^{S+}_{GW}}$ (160 versus 157).
As expected, the performance of the pattern-based formulation $\mathrm{F_P}$ varies significantly with $\mu$.
For $\mu=2$, it solves 207 of the 225 instances, outperforming all other proposed formulations in this group.
In contrast, for $\mu=0.6$ and $\mu=1$, the formulation could not be constructed for 268 of the 450 instances due to the memory limit.
This behavior is consistent with the size of the enumerated pattern set. As $\mu$ increases, the bin capacity decreases, reducing the number of feasible item subsets and making complete pattern enumeration more manageable. Conversely, for smaller values of $\mu$, the larger bin capacity results in a much larger number of feasible patterns, resulting in memory issues.

Regarding the LP relaxations, $\mathrm{F_P}$ provides the strongest bounds whenever its full pattern set is enumerated, while $\mathrm{F^{S}_{RLT1}}$ generally provides the strongest bounds among the formulations, when $\mathrm{F_P}$ fails due to memory.
Despite its strong LP relaxation, the computational performance of $\mathrm{F^{S}_{RLT1}}$ remains limited, illustrating that relaxation strength alone does not compensate for formulation size and complexity.
The formulations $\mathrm{F^{S+}_{G}}$ and $\mathrm{F^{S+}_{FP}}$ are less competitive overall, solving 290 and 254 instances, respectively, with $\mathrm{F^{S+}_{G}}$ providing better computational performance.

Finally, we compare the proposed formulations with the best-performing MILP formulation, EFGW, reported by \cite{CLMI26}.
Among the MILP formulations, $\mathrm{F^{S+}_{GW}}$ clearly outperforms EFGW, solving 46 more instances to optimality.
This indicates that incorporating the strengthened inequalities~\eqref{eq:strong-g} alone is more effective than adding only the capacity-based RLT constraints used in EFGW.
In fact, $\mathrm{F^{S+}_{GW}}$ not only proves more instances to optimality but also provides, on average, a stronger LP relaxation than EFGW.
Note that in the experiments of \cite{CLMI26}, EFGW solved only 222 out of the 405 instances considered, which included only instances with 25, 30, and 35 items. Our reimplemantation solved 293 of the same instances suggesting that we do not underrepresent the performance of EFGW.
We note that although the strengthened inequalities~\eqref{eq:strong-g} improve the LP relaxation, as established in Lemma~\ref{lemma:stronger}, preliminary experiments showed only a limited effect on the number of instances solved and on average solution time.
Comparisons with the corresponding unstrengthened formulations indicated that much of the initial LP relaxation advantage is reduced after Gurobi's presolve.

To better understand why the stronger LP relaxation of $\mathrm{F^{S}_{RLT1}}$ does not translate into the best computational performance, we performed additional experiments on the instances with $n=25$, disabling cutting planes generated by Gurobi.
In these experiments, $\mathrm{F^{S}_{RLT1}}$ solved 123 of the 135 instances, compared to 119 for the $\mathrm{F^{S+}_{GW}}$, confirming that its stronger LP relaxation can be beneficial.
However, compared to $\mathrm{F^{S+}_{GW}}$, $\mathrm{F^{S}_{RLT1}}$ contained approximately four times as many variables, nine times as many constraints and more than ten times as many nonzero coefficients, and its average solution time was higher (595.4 versus 653.9 seconds).
This suggests that the stronger LP relaxation, which comes at the cost of a significant increase in formulation size, does not translate into improved computational performance, as its potential benefit seems to be outweighed by the solver's internal routines.

\subsection{Item-partitioned hybrid formulation}
We next evaluate the item-partitioned hybrid formulation $\mathrm{F^{S+}_{HI}}$ without including any of the valid inequalities introduced in Section~\ref{sec:valid}.
We investigate different rules for constructing the set $\mathcal{A}$ of items for which all feasible patterns are generated.
We consider five selection rules.
Under $R$, the set $\mathcal{A}$ consists of $|\mathcal{A}|$ randomly selected items.
Under $H$, $\mathcal A$ consists of the $|\mathcal{A}|$ heaviest items, whereas under $L$, it consists of the $|\mathcal{A}|$ lightest items.
The rule $I^{(1)}$ is interaction-based and constructs $\mathcal{A}$ greedily.
It starts with the pair of items having the smallest interaction coefficient $d_{ij}$ and subsequently adds, at each iteration, an item outside $\mathcal{A}$ having the smallest pairwise interaction with any item already in $\mathcal{A}$, until the desired cardinality is reached.
Finally, under $I^{(2)}$, each item $i$ is assigned the interaction-strength score
$s_i=\sum_{j\in\mathcal N\setminus\{i\}} |d_{ij}|$ and $\mathcal{A}$ consists of the $|\mathcal{A}|$ items with the smallest values of $s_i$.
We indicate the cardinality of set $\mathcal{A}$ with a subscript under the selection rule in Table~\ref{tab:results_hybrid}, where $|\mathcal{A}|\in\{10,20\}$.

\begin{table*}[ht]
  \centering
  \caption{Computational performance of the formulations, aggregated over $\mu$ and $\delta$}
  \label{tab:results_hybrid}
  \scriptsize
  \setlength{\tabcolsep}{2.2pt}
  \resizebox{\textwidth}{!}{%
  \begin{tabular}{llc rr rr rr rr rr rr rr rr rr rr}
    \toprule
     &  &  & \multicolumn{2}{c}{$R_{10}$} & \multicolumn{2}{c}{$R_{20}$} & \multicolumn{2}{c}{$H_{10}$} & \multicolumn{2}{c}{$H_{20}$} & \multicolumn{2}{c}{$L_{10}$} & \multicolumn{2}{c}{$L_{20}$} & \multicolumn{2}{c}{$I^{(1)}_{10}$} & \multicolumn{2}{c}{$I^{(1)}_{20}$} & \multicolumn{2}{c}{$I^{(2)}_{10}$} & \multicolumn{2}{c}{$I^{(2)}_{20}$} \\
    \cmidrule(lr){4-5} \cmidrule(lr){6-7} \cmidrule(lr){8-9} \cmidrule(lr){10-11} \cmidrule(lr){12-13} \cmidrule(lr){14-15} \cmidrule(lr){16-17} \cmidrule(lr){18-19} \cmidrule(lr){20-21} \cmidrule(lr){22-23}
    $n$ & $\sigma$ & inst & opt & \multicolumn{1}{c}{T(s)} & opt & \multicolumn{1}{c}{T(s)} & opt & \multicolumn{1}{c}{T(s)} & opt & \multicolumn{1}{c}{T(s)} & opt & \multicolumn{1}{c}{T(s)} & opt & \multicolumn{1}{c}{T(s)} & opt & \multicolumn{1}{c}{T(s)} & opt & \multicolumn{1}{c}{T(s)} & opt & \multicolumn{1}{c}{T(s)} & opt & \multicolumn{1}{c}{T(s)} \\
    \midrule
    \num{25} & $+$ & \num{45} & \num{42} & \num{421.3} & \num{44} & \num{303.5} & \num{42} & \num{420.4} & \num{44} & \num{397.8} & \cellcolor{resultgreen}\num{45} & \num{190.1} & \cellcolor{resultgreen}\num{45} & \num{199.4} & \cellcolor{resultgreen}\num{45} & \num{297.8} & \num{44} & \num{253.4} & \cellcolor{resultgreen}\num{45} & \num{263.0} & \num{44} & \num{230.3} \\
    \num{30} & $+$ & \num{45} & \num{37} & \num{816.5} & \num{36} & \num{854.5} & \num{36} & \num{888.3} & \num{37} & \num{877.0} & \cellcolor{resultgreen}\num{38} & \num{706.4} & \cellcolor{resultgreen}\num{38} & \num{795.7} & \cellcolor{resultgreen}\num{38} & \num{783.3} & \num{37} & \num{818.0} & \cellcolor{resultgreen}\num{38} & \num{758.7} & \num{37} & \num{765.4} \\
    \num{35} & $+$ & \num{45} & \num{26} & \num{1832.2} & \num{27} & \num{1928.3} & \num{25} & \num{1849.8} & \num{27} & \num{1805.6} & \cellcolor{resultgreen}\num{33} & \num{1434.1} & \cellcolor{resultgreen}\num{33} & \num{1432.5} & \num{29} & \num{1636.2} & \num{28} & \num{1761.7} & \num{28} & \num{1642.7} & \num{27} & \num{1722.5} \\
    \num{40} & $+$ & \num{45} & \num{18} & \num{2490.1} & \num{15} & \num{2678.8} & \num{16} & \num{2575.4} & \num{15} & \num{2624.4} & \num{16} & \num{2568.7} & \num{17} & \num{2547.2} & \num{18} & \num{2546.7} & \num{12} & \num{2876.2} & \cellcolor{resultgreen}\num{20} & \num{2302.5} & \num{17} & \num{2629.7} \\
    \num{45} & $+$ & \num{45} & \num{12} & \num{2870.9} & \num{7} & \num{3142.2} & \num{8} & \num{3083.6} & \num{6} & \num{3134.5} & \num{12} & \num{2848.9} & \num{11} & \num{2959.3} & \cellcolor{resultgreen}\num{13} & \num{2895.4} & \num{10} & \num{2908.5} & \num{12} & \num{2948.6} & \num{10} & \num{3005.7} \\\midrule
    avg. & $+$ & \num{225} & \num{135} & \num{1686.2} & \num{129} & \num{1781.5} & \num{127} & \num{1763.5} & \num{129} & \num{1767.9} & \cellcolor{resultgreen}\num{144} & \num{1549.6} & \cellcolor{resultgreen}\num{144} & \num{1586.8} & \num{143} & \num{1631.9} & \num{131} & \num{1723.6} & \num{143} & \num{1583.1} & \num{135} & \num{1670.7} \\
    \midrule
    \num{25} & $-$ & \num{45} & \num{43} & \num{746.8} & \num{40} & \num{1011.9} & \cellcolor{resultgreen}\num{44} & \num{721.0} & \num{42} & \num{801.3} & \num{42} & \num{669.1} & \num{36} & \num{1220.0} & \num{43} & \num{778.8} & \num{39} & \num{1104.8} & \num{42} & \num{636.1} & \num{40} & \num{980.9} \\
    \num{30} & $-$ & \num{45} & \num{25} & \num{1956.5} & \num{19} & \num{2404.0} & \num{25} & \num{1929.1} & \num{22} & \num{2219.7} & \num{24} & \num{1911.3} & \num{14} & \num{2492.6} & \num{24} & \num{1952.2} & \num{14} & \num{2511.1} & \cellcolor{resultgreen}\num{26} & \num{1723.0} & \num{22} & \num{2235.5} \\
    \num{35} & $-$ & \num{45} & \num{14} & \num{2915.2} & \num{9} & \num{3205.6} & \num{11} & \num{2986.0} & \num{9} & \num{3153.8} & \num{11} & \num{2940.9} & \num{7} & \num{3227.6} & \num{11} & \num{3039.5} & \num{8} & \num{3233.9} & \cellcolor{resultgreen}\num{15} & \num{2785.2} & \num{12} & \num{3038.4} \\
    \num{40} & $-$ & \num{45} & \num{2} & \num{3550.6} & \num{0} & \num{3600} & \num{2} & \num{3574.5} & \num{0} & \num{3600} & \num{2} & \num{3549.5} & \num{0} & \num{3600} & \num{1} & \num{3575.9} & \num{0} & \num{3600} & \cellcolor{resultgreen}\num{4} & \num{3445.3} & \num{1} & \num{3585.5} \\
    \num{45} & $-$ & \num{45} & \num{0} & \num{3600} & \num{0} & \num{3600} & \num{0} & \num{3600} & \num{0} & \num{3600} & \num{0} & \num{3600} & \num{0} & \num{3600} & \num{0} & \num{3600} & \num{0} & \num{3600} & \num{0} & \num{3600} & \num{0} & \num{3600} \\\midrule
    avg. & $-$ & \num{225} & \num{84} & \num{2553.8} & \num{68} & \num{2764.8} & \num{82} & \num{2562.1} & \num{73} & \num{2675.0} & \num{79} & \num{2534.1} & \num{57} & \num{2829.0} & \num{79} & \num{2589.3} & \num{61} & \num{2810.5} & \cellcolor{resultgreen}\num{87} & \num{2437.9} & \num{75} & \num{2688.3} \\
    \midrule
    \num{25} & $+\!/\!-$ & \num{45} & \cellcolor{resultgreen}\num{45} & \num{5.5} & \cellcolor{resultgreen}\num{45} & \num{23.4} & \cellcolor{resultgreen}\num{45} & \num{6.5} & \cellcolor{resultgreen}\num{45} & \num{13.9} & \cellcolor{resultgreen}\num{45} & \num{5.5} & \cellcolor{resultgreen}\num{45} & \num{36.2} & \cellcolor{resultgreen}\num{45} & \num{7.4} & \cellcolor{resultgreen}\num{45} & \num{22.9} & \cellcolor{resultgreen}\num{45} & \num{4.7} & \cellcolor{resultgreen}\num{45} & \num{18.4} \\
    \num{30} & $+\!/\!-$ & \num{45} & \cellcolor{resultgreen}\num{45} & \num{125.4} & \cellcolor{resultgreen}\num{45} & \num{238.8} & \cellcolor{resultgreen}\num{45} & \num{138.9} & \cellcolor{resultgreen}\num{45} & \num{157.0} & \cellcolor{resultgreen}\num{45} & \num{99.3} & \cellcolor{resultgreen}\num{45} & \num{270.8} & \cellcolor{resultgreen}\num{45} & \num{126.8} & \cellcolor{resultgreen}\num{45} & \num{260.0} & \cellcolor{resultgreen}\num{45} & \num{100.3} & \cellcolor{resultgreen}\num{45} & \num{162.1} \\
    \num{35} & $+\!/\!-$ & \num{45} & \num{43} & \num{473.6} & \num{42} & \num{849.4} & \num{43} & \num{572.1} & \num{42} & \num{748.2} & \cellcolor{resultgreen}\num{45} & \num{388.6} & \num{43} & \num{939.7} & \num{43} & \num{539.8} & \num{41} & \num{1006.9} & \cellcolor{resultgreen}\num{45} & \num{353.0} & \num{44} & \num{670.4} \\
    \num{40} & $+\!/\!-$ & \num{45} & \num{25} & \num{1936.6} & \num{21} & \num{2270.5} & \num{26} & \num{1952.1} & \num{21} & \num{2245.3} & \num{25} & \num{1826.4} & \num{25} & \num{2109.8} & \num{25} & \num{1965.3} & \num{19} & \num{2416.8} & \cellcolor{resultgreen}\num{27} & \num{1807.6} & \num{21} & \num{2155.6} \\
    \num{45} & $+\!/\!-$ & \num{45} & \num{16} & \num{2642.7} & \num{11} & \num{3122.8} & \num{17} & \num{2604.8} & \num{13} & \num{2981.5} & \num{17} & \num{2584.2} & \num{11} & \num{3080.4} & \num{16} & \num{2747.2} & \num{7} & \num{3276.9} & \cellcolor{resultgreen}\num{18} & \num{2493.4} & \num{12} & \num{3056.1} \\\midrule
    avg. & $+\!/\!-$ & \num{225} & \num{174} & \num{1036.8} & \num{164} & \num{1301.0} & \num{176} & \num{1054.9} & \num{166} & \num{1229.2} & \num{177} & \num{980.8} & \num{169} & \num{1287.4} & \num{174} & \num{1077.3} & \num{157} & \num{1396.7} & \cellcolor{resultgreen}\num{180} & \num{951.8} & \num{167} & \num{1212.5} \\
    \midrule\midrule
    total & & \num{675} & \num{393} & \num{1758.9} & \num{361} & \num{1949.1} & \num{385} & \num{1793.5} & \num{368} & \num{1890.7} & \num{400} & \num{1688.2} & \num{370} & \num{1901.1} & \num{396} & \num{1766.2} & \num{349} & \num{1976.9} & \cellcolor{resultgreen}\num{410} & \num{1657.6} & \num{377} & \num{1857.2} \\
    \bottomrule
  \end{tabular}%
  }
\end{table*}
Overall, the configurations with $|\mathcal{A}|=10$ provide the best trade-off.
Although increasing $|\mathcal{A}|$ generally strengthens the LP relaxation, it also increases the number of pattern variables and the resulting increase in formulation size outweighs the benefit of the stronger relaxation.
Among the item-selection rules, $I^{(2)}_{10}$ provides the best overall performance, solving 410 instances with an average solution time of 1657.6 seconds.
This corresponds to 21 additional instances compared to $\mathrm{F^{S+}_{GW}}$, the best-performing non-hybrid formulation, which solves 389 instances.
Its advantage is particularly pronounced for instances with only {negative} dissimilarity values, where it solves 87 instances compared to 84 for $R_{10}$, 82 for $H_{10}$, 79 for $L_{10}$ and 79 for $I^{(1)}_{10}$.
It also performs strongly for mixed-sign instances, solving 180 of the 225 instances.

{One possible explanation for the above-mentioned differences is solver randomness: it was observed in the literature (see, e.g., \cite{AD25}) that small performance difference could be attributed solely to chance.}
To measure the impact of solver randomness, we further run $\mathrm{F^{S+}_{GW}}$ and $I^{(2)}_{10}$ by varying Gurobi's \texttt{Seed} parameter over 5 different values.
Changing the seed parameter of the solver introduces small perturbations in the solver's decisions and results in different solution paths.
For $\mathrm{F^{S+}_{GW}}$, the number of instances solved to optimality ranges from 385 to 390, with average solution time ranging from 1751.4 to 1769.5 seconds. Additionally, 376 instances are solved to optimality under every seed value, whereas 400 are solved by at least one seed.
For $I^{(2)}_{10}$, the range of number of instances solved to optimality is between 406 and 413, with average solution time ranging from 1646.2 to 1667.1 seconds.
Under every seed value, 395 instances are solved to optimality, while 426 are solved to optimality for at least one seed. We observe that the results are relatively stable across different seed values and, overall, reiterate that $I^{(2)}_{10}$ outperforms $\mathrm{F^{S+}_{GW}}$.

To better understand these differences, we report in Table~\ref{tab:formulation_characteristics} the average LP-relaxation value $z_{\mathrm{LP}}$, the number of non-pattern variables ``$\mathrm{var}_{NP}$", the number of pattern variables ``$\mathrm{var}_{P}$" and the number of constraints ``const'' for the configurations where $|\mathcal{A}|=10$.
For each sign class, the strongest LP relaxation and the smallest values of the formulation-size statistics are highlighted in green.

\begin{table*}[ht]
  \centering
  \caption{Average formulation characteristics by sign structure.}
  \label{tab:formulation_characteristics}
  \scriptsize
  \setlength{\tabcolsep}{2.2pt}
  \resizebox{\textwidth}{!}{%
  \begin{tabular}{l rrrr rrrr rrrr rrrr rrrr}
    \toprule
     & \multicolumn{4}{c}{$R_{10}$} & \multicolumn{4}{c}{$H_{10}$} & \multicolumn{4}{c}{$L_{10}$} & \multicolumn{4}{c}{$I^{(1)}_{10}$} & \multicolumn{4}{c}{$I^{(2)}_{10}$} \\
    \cmidrule(lr){2-5} \cmidrule(lr){6-9} \cmidrule(lr){10-13} \cmidrule(lr){14-17} \cmidrule(lr){18-21}
     $\sigma$ & $z_{\mathrm{LP}}$ & $\mathrm{var}_{NP}$ & $\mathrm{var}_{P}$ & const & $z_{\mathrm{LP}}$ & $\mathrm{var}_{NP}$ & $\mathrm{var}_{P}$ & const & $z_{\mathrm{LP}}$ & $\mathrm{var}_{NP}$ & $\mathrm{var}_{P}$ & const & $z_{\mathrm{LP}}$ & $\mathrm{var}_{NP}$ & $\mathrm{var}_{P}$ & const & $z_{\mathrm{LP}}$ & $\mathrm{var}_{NP}$ & $\mathrm{var}_{P}$ & const \\ \midrule
     $+$ &  \num{13637} &  \num{4259} &  \num{687} &  \num{4294} &  \num{13638} &  \num{4262} & \cellcolor{resultgreen} \num{405} &  \num{4297} &  \num{13637} &  \num{4252} &  \num{1000} &  \num{4287} &  \num{13647} &  \num{4174} &  \num{681} &  \num{4209} & \cellcolor{resultgreen} \num{13655} & \cellcolor{resultgreen} \num{4085} &  \num{691} & \cellcolor{resultgreen} \num{4120} \\
    $-$ &  \num{1435} &  \num{4242} &  \num{689} &  \num{7919} &  \num{1435} &  \num{4240} & \cellcolor{resultgreen} \num{404} &  \num{7914} &  \num{1460} &  \num{4248} &  \num{1000} &  \num{7931} &  \num{1325} &  \num{4307} &  \num{688} &  \num{8049} & \cellcolor{resultgreen} \num{1775} & \cellcolor{resultgreen} \num{4078} &  \num{700} & \cellcolor{resultgreen} \num{7590} \\
    $+\!/\!-$ &  \num{-2563} &  \num{4240} &  \num{716} &  \num{6093} &  \num{-2568} &  \num{4243} & \cellcolor{resultgreen} \num{403} &  \num{6102} &  \num{-2549} &  \num{4231} &  \num{999} &  \num{6079} &  \num{-2614} &  \num{4283} &  \num{690} &  \num{6214} & \cellcolor{resultgreen} \num{-2486} & \cellcolor{resultgreen} \num{4071} &  \num{705} & \cellcolor{resultgreen} \num{5845} \\
    \midrule\midrule
     all&  \num{4169} &  \num{4247} &  \num{697} &  \num{6102} &  \num{4168} &  \num{4248} & \cellcolor{resultgreen} \num{404} &  \num{6104} &  \num{4183} &  \num{4244} &  \num{1000} &  \num{6099} &  \num{4119} &  \num{4255} &  \num{686} &  \num{6158} & \cellcolor{resultgreen} \num{4315} & \cellcolor{resultgreen} \num{4078} &  \num{699} & \cellcolor{resultgreen} \num{5852} \\
    \bottomrule
  \end{tabular}%
  }
\end{table*}

Since fewer combinations of heavy items satisfy the capacity constraints, $H_{10}$ yields, as expected, the smallest pattern component, with an average of 404 pattern variables.
In contrast, $L_{10}$ generates the largest pattern component, with an average of 1000 pattern variables, since a larger number of combinations of light items are capacity-feasible.
The number of pattern variables does not appear to determine LP relaxation strength directly. Although $L_{10}$ generates the most patterns and $H_{10}$ the fewest, neither yields the strongest relaxation.
In contrast, $I^{(2)}_{10}$ achieves the largest average LP value, $z_{\mathrm{LP}}=4315$, while also producing the smallest non-pattern component, with 4078 variables and 5852 constraints on average.
Hence, the superior performance of $I^{(2)}_{10}$ appears to result from the combination of two effects: a stronger relaxation and a smaller model size related to the assignment-based component.

\subsection{Effect of item ordering}

Unlike $I^{(1)}$, $I^{(2)}$ does not explicitly seek attractive patterns. Instead, it places in $\mathcal{A}$ the items with the smallest total absolute interaction, leaving the more interaction-intensive items in $\bar{\mathcal A}$ to be handled by the assignment-based component of the formulation.
Since the assignment-based component is order-dependent, this selection also affects the size and structure of the formulation.

This observation motivates a closer examination of the internal ordering of the items in $\bar{\mathcal A}$.
Item ordering affects the assignment-based formulation in two distinct ways.
First, it affects the LP relaxation because the symmetry-breaking constraints~\eqref{SYM1} and \eqref{SYM2} are index dependent.
Consequently, reordering the items changes which item weights and interaction values are subject to these restrictions in the capacity and linearization constraints.
Although the set of feasible integer packings is unchanged, the corresponding fractional feasible region and hence the LP relaxation value, can change.
Second, the ordering affects the size of the formulation.
For a pair of items having indices $i<j$ with $d_{ij}\neq 0$, the formulation contains variables $g_{ijk}$ only for $k<i<j$.
Hence, this interaction requires $i-1$ $g$-variables.
Consequently, different permutations of the same items can lead to different numbers of $g$-variables and, in turn, different numbers of associated constraints and nonzero coefficients.
For example, consider four items $a,b,c,d$ with nonzero interactions only between item $a$ and each of items $b$, $c$ and $d$.
Under the ordering $(a,b,c,d)$, these interactions correspond to index pairs $(1,2)$, $(1,3)$ and $(1,4)$.
Since the first index is $i=1$ in each case, there is no index $k<i$, and hence no $g$-variables are introduced.
Under the ordering $(b,c,d,a)$, the same three interactions correspond to index pairs $(1,4)$, $(2,4)$ and $(3,4)$.
The pair $(1,4)$ generates no $g$-variables, the pair $(2,4)$ generates $g_{241}$, and the pair $(3,4)$ generates $g_{341}$ and $g_{342}$.
Thus, the second ordering generates three $g$-variables, whereas the first generates none.
This example illustrates why placing highly connected items earlier in the ordering can reduce the size of the formulation.

We therefore define the interaction degree of item $i$ as $\deg(i)=|\{j\in\mathcal N\setminus\{i\}:$ $ d_{ij}\neq 0\}|$, that is, the number of items with which $i$ has a nonzero interaction.
We use decreasing interaction degree as a simple way of placing highly connected items earlier in the ordering.
To isolate the ordering effect from the selection effect, we fix the set $\mathcal A$ generated by $I^{(2)}_{10}$ and modify only the ordering of the items in $\bar{\mathcal A}$ before constructing the assignment-based component.
We consider three orderings: decreasing interaction degree, increasing interaction degree and a random ordering.
The decreasing-degree ordering places highly connected items first and is intended to reduce the size of the assignment-based component, while the increasing-degree ordering provides the opposite case and the random ordering serves as a control.
Since $\mathcal{A}$ is fixed, all three variants generate the same pattern variables and differ only in the ordering of the assignment-based representation.
We emphasize that
since all items in $\mathcal A$ are placed after those in $\bar{\mathcal A}$, permuting the items within $\mathcal{A}$ does not change the numbers of pattern variables ($\lambda_p$), assignment-based variables ($x_{ik}$ and $g_{ijk}$), constraints or nonzero coefficients.

Table~\ref{tab:ordering_analysis} reports the results for the decreasing-degree, increasing-degree and random ordering variants $I^{(2)}_{10}\text{-}D^\downarrow$, $I^{(2)}_{10}\text{-}D^\uparrow$ and
$I^{(2)}_{10}\text{-}R$.
The results are shown by the interaction density $\delta$, since $\delta$ determines the sparsity of the dissimilarity matrix and therefore directly affects the number of $g$-variables that can be eliminated through the ordering of the assignment-based component.
Since the set $\mathcal{A}$ is fixed, the column reporting the number of pattern variables is identical across the $I^{(2)}_{10}$ variants.
Yellow cells indicate the best computational performance, determined first by the highest number of instances solved to proven optimality and, in case of a tie, by the lowest average solution time.
Green cells indicate the strongest LP-relaxation value, while blue cells indicate the smallest numbers of assignment-based variables and constraints within each instance group.

\begin{table*}[!htbp]
  \centering
  \caption{Computational performance of the $I^{(2)}_{10}$ variants.}
  \label{tab:ordering_analysis}
  \scriptsize
  \setlength{\tabcolsep}{1.15pt}
  \resizebox{\textwidth}{!}{%
  \begin{tabular}{ccc @{\hspace{5pt}} rr@{\hspace{4pt}}rrrr @{\hspace{7pt}} rrrrrr @{\hspace{7pt}} rr@{\hspace{4pt}}rrrr @{\hspace{7pt}} rr@{\hspace{4pt}}rrrr}
    \toprule
    & & & \multicolumn{6}{c}{$I^{(2)}_{10}$} & \multicolumn{6}{c}{$I^{(2)}_{10}\text{-}D^\downarrow$} & \multicolumn{6}{c}{$I^{(2)}_{10}\text{-}R$} & \multicolumn{6}{c}{$I^{(2)}_{10}\text{-}D^\uparrow$} \\
    \cmidrule(lr){4-9} \cmidrule(lr){10-15} \cmidrule(lr){16-21} \cmidrule(lr){22-27}
    $n$ & $\delta$ & inst & opt & T(s) & \multicolumn{1}{c}{$z_{\mathrm{LP}}$} & $\mathrm{var}_{NP}$ & $\mathrm{var}_{P}$ & const & opt & T(s) & \multicolumn{1}{c}{$z_{\mathrm{LP}}$} & $\mathrm{var}_{NP}$ & $\mathrm{var}_{P}$ & const & opt & T(s) & \multicolumn{1}{c}{$z_{\mathrm{LP}}$} & $\mathrm{var}_{NP}$ & $\mathrm{var}_{P}$ & const & opt & T(s) &  \multicolumn{1}{c}{$z_{\mathrm{LP}}$} & $\mathrm{var}_{NP}$ & $\mathrm{var}_{P}$ & const \\
    \midrule
    \num{25} & 0.25 & \num{45} & \num{45} & \num{31.8} & \num{3434.0} & \num{744} & \num{526} & \num{997} & \cellcolor{resultyellow}\num{45} & \cellcolor{resultyellow}\num{21.7} & \cellcolor{resultgreen}\num{3486.7} & \cellcolor{resultblue}\num{663} & \cellcolor{resultblue}\num{526} & \cellcolor{resultblue}\num{878} & \num{45} & \num{35.5} & \num{3435.2} & \num{738} & \num{526} & \num{987} & \num{45} & \num{46.1} & \num{3386.0} & \num{813} & \num{526} & \num{1100}  \\
    \num{30} & 0.25 & \num{45} & \num{43} & \num{253.1} & \num{2829.6} & \num{1288} & \num{656} & \num{1757} & \cellcolor{resultyellow}\num{44} & \cellcolor{resultyellow}\num{182.9} & \cellcolor{resultgreen}\num{2912.5} & \cellcolor{resultblue}\num{1155} & \cellcolor{resultblue}\num{656} & \cellcolor{resultblue}\num{1560} & \num{44} & \num{234.3} & \num{2819.0} & \num{1276} & \num{656} & \num{1739} & \num{43} & \num{367.5} & \num{2753.4} & \num{1418} & \num{656} & \num{1949}  \\
    \num{35} & 0.25 & \num{45} & \num{39} & \num{722.7} & \num{2320.0} & \num{2034} & \num{715} & \num{2807} & \cellcolor{resultyellow}\num{43} & \cellcolor{resultyellow}\num{502.1} & \cellcolor{resultgreen}\num{2461.0} & \cellcolor{resultblue}\num{1825} & \cellcolor{resultblue}\num{715} & \cellcolor{resultblue}\num{2489} & \num{40} & \num{733.7} & \num{2277.5} & \num{2050} & \num{715} & \num{2838} & \num{39} & \num{840.8} & \num{2182.3} & \num{2262} & \num{715} & \num{3158}  \\
    \num{40} & 0.25 & \num{45} & \num{29} & \num{1614.3} & \num{1536.1} & \num{3003} & \num{756} & \num{4153} & \cellcolor{resultyellow}\num{36} & \cellcolor{resultyellow}\num{1106.0} & \cellcolor{resultgreen}\num{1747.8} & \cellcolor{resultblue}\num{2655} & \cellcolor{resultblue}\num{756} & \cellcolor{resultblue}\num{3642} & \num{30} & \num{1638.9} & \num{1570.3} & \num{2994} & \num{756} & \num{4130} & \num{24} & \num{1831.8} & \num{1379.6} & \num{3378} & \num{756} & \num{4720}  \\
    \num{45} & 0.25 & \num{45} & \num{22} & \num{2097.7} & \num{717.5} & \num{4282} & \num{819} & \num{5973} & \cellcolor{resultyellow}\num{26} & \cellcolor{resultyellow}\num{1752.2} & \cellcolor{resultgreen}\num{987.1} & \cellcolor{resultblue}\num{3764} & \cellcolor{resultblue}\num{819} & \cellcolor{resultblue}\num{5182} & \num{24} & \num{2021.9} & \num{723.5} & \num{4259} & \num{819} & \num{5928} & \num{22} & \num{2222.4} & \num{535.9} & \num{4800} & \num{819} & \num{6756}  \\ \midrule
     \multicolumn{2}{c}{avg.} &  \num{225} &  \num{178} &  \num{943.9} &  \num{2167.4} &  \num{2270} &  \num{694} &  \num{3137} &  \cellcolor{resultyellow}\num{194} &  \cellcolor{resultyellow}\num{713.0} &  \cellcolor{resultgreen}\num{2319.0} &  \cellcolor{resultblue}\num{2013} &  \cellcolor{resultblue}\num{694} &  \cellcolor{resultblue}\num{2750} &  \num{183} &  \num{932.8} &  \num{2165.1} &  \num{2263} &  \num{694} &  \num{3124} &  \num{173} &  \num{1061.7} &  \num{2047.4} &  \num{2534} &  \num{694} &  \num{3537}  \\
    \midrule
    \num{25} & 0.50 & \num{45} & \num{44} & \num{221.7} & \num{6285.8} & \num{1276} & \num{544} & \num{1821} & \cellcolor{resultyellow}\num{45} & \cellcolor{resultyellow}\num{168.2} & \cellcolor{resultgreen}\num{6377.2} & \cellcolor{resultblue}\num{1186} & \cellcolor{resultblue}\num{544} & \cellcolor{resultblue}\num{1686} & \num{45} & \num{196.0} & \num{6269.7} & \num{1278} & \num{544} & \num{1835} & \num{44} & \num{253.9} & \num{6229.6} & \num{1364} & \num{544} & \num{1957}  \\
    \num{30} & 0.50 & \num{45} & \num{37} & \num{913.2} & \num{5700.0} & \num{2259} & \num{634} & \num{3191} & \cellcolor{resultyellow}\num{37} & \cellcolor{resultyellow}\num{829.1} & \cellcolor{resultgreen}\num{5814.4} & \cellcolor{resultblue}\num{2106} & \cellcolor{resultblue}\num{634} & \cellcolor{resultblue}\num{2974} & \num{36} & \num{937.1} & \num{5678.7} & \num{2265} & \num{634} & \num{3208} & \num{36} & \num{1043.1} & \num{5615.2} & \num{2403} & \num{634} & \num{3404}  \\
    \num{35} & 0.50 & \num{45} & \num{25} & \num{1907.2} & \num{4550.7} & \num{3646} & \num{704} & \num{5188} & \cellcolor{resultyellow}\num{26} & \cellcolor{resultyellow}\num{1825.8} & \cellcolor{resultgreen}\num{4785.3} & \cellcolor{resultblue}\num{3368} & \cellcolor{resultblue}\num{704} & \cellcolor{resultblue}\num{4761} & \num{24} & \num{1934.8} & \num{4566.1} & \num{3643} & \num{704} & \num{5185} & \num{22} & \num{2043.7} & \num{4401.9} & \num{3917} & \num{704} & \num{5620}  \\
    \num{40} & 0.50 & \num{45} & \num{14} & \num{2732.7} & \num{3326.3} & \num{5327} & \num{761} & \num{7680} & \cellcolor{resultyellow}\num{16} & \cellcolor{resultyellow}\num{2587.3} & \cellcolor{resultgreen}\num{3546.5} & \cellcolor{resultblue}\num{4927} & \cellcolor{resultblue}\num{761} & \cellcolor{resultblue}\num{7066} & \num{14} & \num{2689.5} & \num{3304.5} & \num{5320} & \num{761} & \num{7682} & \num{14} & \num{2870.6} & \num{3094.9} & \num{5741} & \num{761} & \num{8286}  \\
    \num{45} & 0.50 & \num{45} & \num{8} & \num{3344.5} & \num{1574.5} & \num{7684} & \num{833} & \num{11091} & \cellcolor{resultyellow}\num{13} & \cellcolor{resultyellow}\num{3072.8} & \cellcolor{resultgreen}\num{1885.8} & \cellcolor{resultblue}\num{7082} & \cellcolor{resultblue}\num{833} & \cellcolor{resultblue}\num{10153} & \num{9} & \num{3290.5} & \num{1657.4} & \num{7657} & \num{833} & \num{11006} & \num{3} & \num{3565.0} & \num{1375.7} & \num{8258} & \num{833} & \num{11936}  \\ \midrule
     \multicolumn{2}{c}{avg.}&  \num{225} &  \num{128} &  \num{1823.9} &  \num{4287.5} &  \num{4039} &  \num{695} &  \num{5794} &  \cellcolor{resultyellow}\num{137} &  \cellcolor{resultyellow}\num{1696.6} &  \cellcolor{resultgreen}\num{4481.8} &  \cellcolor{resultblue}\num{3734} &  \cellcolor{resultblue}\num{695} &  \cellcolor{resultblue}\num{5328} &  \num{128} &  \num{1809.6} &  \num{4295.3} &  \num{4033} &  \num{695} &  \num{5783} &  \num{119} &  \num{1955.2} &  \num{4143.5} &  \num{4337} &  \num{695} &  \num{6240}  \\
    \midrule
    \num{25} & 0.75 & \num{45} & \cellcolor{resultyellow}\num{43} & \cellcolor{resultyellow}\num{650.4} & \num{9769.1} & \num{1855} & \num{583} & \num{2665} & \num{43} & \num{706.8} & \cellcolor{resultgreen}\num{9836.1} & \cellcolor{resultblue}\num{1766} & \cellcolor{resultblue}\num{583} & \cellcolor{resultblue}\num{2540} & \num{40} & \num{787.8} & \num{9787.0} & \num{1843} & \num{583} & \num{2637} & \num{43} & \num{663.6} & \num{9730.5} & \num{1923} & \num{583} & \num{2760}  \\
    \num{30} & 0.75 & \num{45} & \num{29} & \num{1415.7} & \num{8583.8} & \num{3287} & \num{623} & \num{4776} & \cellcolor{resultyellow}\num{29} & \cellcolor{resultyellow}\num{1382.7} & \cellcolor{resultgreen}\num{8710.5} & \cellcolor{resultblue}\num{3145} & \cellcolor{resultblue}\num{623} & \cellcolor{resultblue}\num{4549} & \num{29} & \num{1394.7} & \num{8605.7} & \num{3285} & \num{623} & \num{4743} & \num{29} & \num{1411.7} & \num{8498.2} & \num{3445} & \num{623} & \num{5004}  \\
    \num{35} & 0.75 & \num{45} & \num{24} & \num{2151.0} & \num{6998.2} & \num{5257} & \num{711} & \num{7625} & \cellcolor{resultyellow}\num{26} & \cellcolor{resultyellow}\num{1948.0} & \cellcolor{resultgreen}\num{7146.8} & \cellcolor{resultblue}\num{4986} & \cellcolor{resultblue}\num{711} & \cellcolor{resultblue}\num{7234} & \num{25} & \num{2089.8} & \num{6985.1} & \num{5230} & \num{711} & \num{7599} & \num{21} & \num{2344.6} & \num{6897.4} & \num{5468} & \num{711} & \num{7967}  \\
    \num{40} & 0.75 & \num{45} & \cellcolor{resultyellow}\num{8} & \cellcolor{resultyellow}\num{3208.5} & \num{4912.6} & \num{7938} & \num{787} & \num{11530} & \num{7} & \num{3155.3} & \cellcolor{resultgreen}\num{5161.4} & \cellcolor{resultblue}\num{7530} & \cellcolor{resultblue}\num{787} & \cellcolor{resultblue}\num{10943} & \num{6} & \num{3272.1} & \num{4945.6} & \num{7919} & \num{787} & \num{11547} & \num{4} & \num{3429.9} & \num{4763.7} & \num{8302} & \num{787} & \num{12120}  \\
    \num{45} & 0.75 & \num{45} & \num{0} & \num{3600} & \num{2183.0} & \num{11285} & \num{834} & \num{16518} & \cellcolor{resultyellow}\num{3} & \cellcolor{resultyellow}\num{3529.4} & \cellcolor{resultgreen}\num{2352.3} & \cellcolor{resultblue}\num{10781} & \cellcolor{resultblue}\num{834} & \cellcolor{resultblue}\num{15797} & \num{0} & \num{3600} & \num{2189.9} & \num{11306} & \num{834} & \num{16588} & \num{0} & \num{3600} & \num{2039.2} & \num{11796} & \num{834} & \num{17309}  \\ \midrule
     \multicolumn{2}{c}{avg.} &  \num{225} &  \num{104} &  \num{2205.1} &  \num{6489.3} &  \num{5924} &  \num{708} &  \num{8623} &  \cellcolor{resultyellow}\num{108} &  \cellcolor{resultyellow}\num{2144.4} &  \cellcolor{resultgreen}\num{6641.4} &  \cellcolor{resultblue}\num{5642} &  \cellcolor{resultblue}\num{708} &  \cellcolor{resultblue}\num{8213} &  \num{100} &  \num{2228.8} &  \num{6502.7} &  \num{5917} &  \num{708} &  \num{8623} &  \num{97} &  \num{2290.0} &  \num{6385.8} &  \num{6187} &  \num{708} &  \num{9032}  \\
    \midrule\midrule
     \multicolumn{2}{c}{total}  &  \num{675} &  \num{410} &  \num{1657.6} &  \num{4314.7} &  \num{4078} &  \num{699} &  \num{5852} &  \cellcolor{resultyellow}\num{439} &  \cellcolor{resultyellow}\num{1518.0} &  \cellcolor{resultgreen}\num{4480.8} &  \cellcolor{resultblue}\num{3796} &  \cellcolor{resultblue}\num{699} &  \cellcolor{resultblue}\num{5430} &  \num{411} &  \num{1657.1} &  \num{4321.0} &  \num{4071} &  \num{699} &  \num{5843} &  \num{389} &  \num{1769.0} &  \num{4192.2} &  \num{4353} &  \num{699} &  \num{6270}  \\
    \bottomrule
  \end{tabular}%
  }
\end{table*}
The results show a clear effect of the ordering of the items in $\bar{\mathcal{A}}$.
Overall, the decreasing-degree ordering reduces the number of assignment-based variables by 6.9\% and the number of constraints by 7.2\%, while increasing the average LP-relaxation value by 3.8\%.
These improvements translate into 29 additional instances solved to optimality and an 8.4\% reduction in average solution time.
The opposite behavior is observed for the increasing-degree ordering.
It produces the largest assignment-based component and the weakest average LP relaxation among the considered orderings and solves only 389 instances.
The random ordering remains close to the original ordering in both formulation size and LP relaxation strength.
These results indicate that the improvement obtained by the decreasing-degree ordering is not only an effect of relabeling.
Placing highly connected items earlier in the assignment-based representation simultaneously reduces the number of required variables and appears to restrict the fractional representation of their interactions, resulting in a stronger relaxation.

The effect of the ordering is particularly pronounced for sparse instances.
For $\delta=0.25$, the decreasing-degree ordering reduces the average number of non-pattern variables by approximately 11\%, whereas the reductions for $\delta=0.50$ and $\delta=0.75$ are approximately 8\% and 5\%, respectively.
A similar pattern is observed for the number of constraints.
The relative benefit of reordering is larger for sparse dissimilarity matrices.
The larger structural reduction for instances with sparser dissimilarity matrices is also accompanied by a larger computational improvement: the decreasing-degree ordering solves 16 additional instances for $\delta=0.25$, compared with 9 for $\delta=0.50$ and 4 for $\delta=0.75$.

The results show that ordering the items in $\bar{\mathcal{A}}$ by decreasing interaction degree improves both the size and computational performance of the item-partitioned hybrid formulation.
Although this ordering was introduced in the hybrid setting, the same idea can also be applied to the best-performing non-hybrid formulation, $\mathrm{F^{S+}_{GW}}$.
Since no subset $\mathcal{A}$ is selected for pattern generation, the ordering is applied directly to the complete item set $\mathcal{N}$.
Following the stronger effect observed for sparser instances in the hybrid formulation, we report the non-hybrid results aggregated by interaction density $\delta$.
Table~\ref{tab:gw_ordering_analysis} compares the formulation $\mathrm{F^{S+}_{GW}}$ with the decreasing-degree variant $\mathrm{F^{S+,D^\downarrow}_{GW}}$.

\begin{table*}[!htbp]
  \centering
  \caption{Effect of the decreasing interaction-degree ordering on $\mathrm{F^{S+}_{GW}}$ by interaction density.}
  \label{tab:gw_ordering_analysis}
  \scriptsize
    \setlength{\tabcolsep}{2pt}
  \begin{tabular}{@{}cc @{\hspace{8pt}} rrrrrr @{\hspace{10pt}} rrrrrr@{}}
    \toprule
    & &
    \multicolumn{6}{c}{$\mathrm{F^{S+}_{GW}}$} &
    \multicolumn{6}{c}{$\mathrm{F^{S+,D^\downarrow}_{GW}}$} \\
    \cmidrule(lr){3-8} \cmidrule(lr){9-14}
    $\delta$ & inst &
    opt & T(s) & $z_{\mathrm{LP}}$ & var & const & nonzero &
    opt & T(s) & $z_{\mathrm{LP}}$ & var & const & nonzero \\
    \midrule
    \num{0.25} & \num{225} &
    \num{172} & \num{1039.3} & \num{2067.7} &
    \num{2504} & \num{3458} & \num{9947} &
    \num{191} & \num{744.7} & \num{2310.7} &
    \num{2072} & \num{2813} & \num{8219} \\

    \num{0.50} & \num{225} &
    \num{122} & \num{1890.3} & \num{4121.0} &
    \num{4361} & \num{6253} & \num{17376} &
    \num{137} & \num{1678.9} & \num{4471.8} &
    \num{3808} & \num{5412} & \num{15163} \\

    \num{0.75} & \num{225} &
    \num{95} & \num{2329.7} & \num{6323.7} &
    \num{6215} & \num{9023} & \num{24792} &
    \num{108} & \num{2153.5} & \num{6607.0} &
    \num{5732} & \num{8319} & \num{22860} \\

    \midrule\midrule
    all & \num{675} &
    \num{389} & \num{1753.1} & \num{4170.8} &
    \num{4360} & \num{6245} & \num{17371} &
    \num{436} & \num{1525.7} & \num{4463.2} &
    \num{3871} & \num{5515} & \num{15414} \\
    \bottomrule
  \end{tabular}
\end{table*}

The decreasing-degree ordering leads to a significant improvement in the non-hybrid formulation.
Across all 675 instances, $\mathrm{F^{S+,D^\downarrow}_{GW}}$ solves 47 additional instances to proven optimality and reduces the average solution time by approximately 13\% compared to the original ordering.
It also reduces the average numbers of variables, constraints and nonzero coefficients by approximately 11\%, 12\% and 11\%, respectively, while improving the average value of $z_{\mathrm{LP}}$ by approximately 7\%.
Thus, the improvement is not solely attributable to the smaller formulation: the ordering also produces a stronger relaxation.

The same density-dependent pattern is observed in the non-hybrid formulation.
For $\delta=0.25$, the decreasing-degree ordering reduces the average numbers of variables, constraints and nonzero coefficients by approximately 17\%, 19\% and 17\%, respectively, and reduces the average solution time by approximately 28\%, while solving 19 additional instances.
The structural reductions decrease to approximately 13\% for $\delta=0.50$ and 8\% for $\delta=0.75$, while the numbers of additional instances solved are 15 and 13, respectively.

Overall, these results show that both item selection and internal ordering affect the item-partitioned hybrid formulation.
The choice of $\mathcal{A}$ determines the structure of the remaining assignment-based component, while, once $\mathcal{A}$ is fixed, ordering the items in $\bar{\mathcal{A}}$ by decreasing interaction degree reduces the formulation size and strengthens its LP relaxation.
The improvement observed for $\mathrm{F^{S+}_{GW}}$ further indicates that the benefit of this ordering is not specific to the hybrid formulation.
We therefore use the decreasing-degree ordering for both the non-hybrid and hybrid formulations in the subsequent experiments.

\subsection{Effect of the bin-count inequalities}

We next evaluate the computational effect of the valid inequalities introduced in Section~\ref{sec:valid}.
We consider the best-performing non-hybrid formulation, $\mathrm{F^{S+,D^\downarrow}_{GW}}$, and the best-performing item-partitioned hybrid formulation, $\mathrm{F^{S+}_{HI}}$ with selection rule $I^{(2)}_{10}$ and decreasing interaction-degree ordering.
We compare three variants of each formulation.
The ``base'' variant does not include any additional inequalities.
The ``global'' variant adds the global bin-count inequality~\eqref{eq:bin_count} for $\mathrm{F^{S+,D^\downarrow}_{GW}}$ and~\eqref{eq:HI_bin_count} for $\mathrm{F^{S+}_{HI}}$.
The ``global+cumulative'' variant additionally exploits the item ordering through the cumulative bin-count inequalities.
For $\mathrm{F^{S+,D^\downarrow}_{GW}}$, this corresponds to adding~\eqref{eq:cumulative_bin_count}, with the inequality for $m=n$ coinciding with the global bin-count inequality~\eqref{eq:bin_count}.
For $\mathrm{F^{S+}_{HI}}$, it corresponds to adding both~\eqref{eq:HI_bin_count} and~\eqref{eq:HI_cumulative_bin_count}.
As described in Section~\ref{sec:valid}, cumulative bin-count inequalities are included only for values of $m$ at which $L_m$ increases.

We report the results in Table~\ref{tab:global_cumulative_analysis}.
For each formulation and instance class, green highlights the strongest LP-relaxation value, while yellow indicates the best computational performance, determined first by the number of instances solved to optimality and then by the average solution time in case of a tie.
For reference, Table~\ref{tab:global_cumulative_analysis} also reports the computational performance of BP(5,10), the B\&P algorithm proposed by \cite{CLMI26} and currently the best-performing exact method on these benchmark instances.
The B\&P method is a specialized exact algorithm, whereas the formulations considered here are solved directly with a general-purpose MILP solver.

\begin{table*}[!htbp]
  \centering
  \caption{Effect of the global and cumulative inequalities on the best non-hybrid and hybrid formulations.}
  \label{tab:global_cumulative_analysis}
  \scriptsize
  \setlength{\tabcolsep}{1.3pt}
  \resizebox{\textwidth}{!}{%
  \begin{tabular}{lll @{\hspace{5pt}} rrr @{\hspace{7pt}} rrr @{\hspace{5pt}} rrr @{\hspace{5pt}} rrr @{\hspace{7pt}} rrr @{\hspace{5pt}} rrr @{\hspace{5pt}} rrr}
    \toprule
    & & & \multicolumn{3}{c}{\cite{CLMI26}} & \multicolumn{9}{c}{$\mathrm{F^{S+,D^\downarrow}_{GW}}$} & \multicolumn{9}{c}{$I^{(2)}_{10}\text{-}D^\downarrow$} \\
    \cmidrule(lr){4-6} \cmidrule(lr){7-15} \cmidrule(lr){16-24}
    & & & \multicolumn{3}{c}{BP(5,10)} & \multicolumn{3}{c}{base} & \multicolumn{3}{c}{global} & \multicolumn{3}{c}{global+cumulative} & \multicolumn{3}{c}{base} & \multicolumn{3}{c}{global} & \multicolumn{3}{c}{global+cumulative} \\
    \cmidrule(lr){4-6} \cmidrule(lr){7-9} \cmidrule(lr){10-12} \cmidrule(lr){13-15} \cmidrule(lr){16-18} \cmidrule(lr){19-21} \cmidrule(lr){22-24}
    $n$ & $\sigma$ & inst & opt & T(s) & $\mathrm{root}_{\mathrm{LP}}$ & opt & T(s) & \multicolumn{1}{c}{$z_{\mathrm{LP}}$} & opt & T(s) & \multicolumn{1}{c}{$z_{\mathrm{LP}}$} & opt & T(s) & \multicolumn{1}{c}{$z_{\mathrm{LP}}$} & opt & T(s) & \multicolumn{1}{c}{$z_{\mathrm{LP}}$} & opt & T(s) & \multicolumn{1}{c}{$z_{\mathrm{LP}}$} & opt & T(s) & \multicolumn{1}{c}{$z_{\mathrm{LP}}$} \\
    \midrule
    \num{25} & $+$ & \num{45} & \num{41} & \num{334.1} & \num{13123.3} & \num{45} & \num{188.1} & \num{12442.8} & \num{45} & \num{133.7} & \cellcolor{resultgreen}\num{13279.3} & \cellcolor{resultyellow}\num{45} & \cellcolor{resultyellow}\num{126.2} & \cellcolor{resultgreen}\num{13279.3} & \num{44} & \num{265.3} & \num{12443.9} & \num{45} & \num{156.1} & \cellcolor{resultgreen}\num{13280.2} & \cellcolor{resultyellow}\num{45} & \cellcolor{resultyellow}\num{152.4} & \cellcolor{resultgreen}\num{13280.2} \\
    \num{30} & $+$ & \num{45} & \num{42} & \num{352.1} & \num{13514.3} & \num{38} & \num{649.5} & \num{13215.0} & \cellcolor{resultyellow}\num{38} & \cellcolor{resultyellow}\num{632.8} & \cellcolor{resultgreen}\num{13332.2} & \num{38} & \num{647.6} & \cellcolor{resultgreen}\num{13332.2} & \num{38} & \num{667.2} & \num{13215.1} & \num{38} & \num{635.5} & \cellcolor{resultgreen}\num{13332.3} & \cellcolor{resultyellow}\num{38} & \cellcolor{resultyellow}\num{634.6} & \cellcolor{resultgreen}\num{13332.3} \\
    \num{35} & $+$ & \num{45} & \num{41} & \num{723.5} & \num{14054.5} & \num{33} & \num{1457.1} & \num{13735.2} & \cellcolor{resultyellow}\num{38} & \cellcolor{resultyellow}\num{920.5} & \cellcolor{resultgreen}\num{13806.5} & \num{37} & \num{981.2} & \cellcolor{resultgreen}\num{13806.5} & \num{33} & \num{1308.2} & \num{13738.6} & \num{36} & \num{1024.6} & \cellcolor{resultgreen}\num{13809.9} & \cellcolor{resultyellow}\num{37} & \cellcolor{resultyellow}\num{973.3} & \cellcolor{resultgreen}\num{13809.9} \\
    \num{40} & $+$ & \num{45} & \num{36} & \num{1213.4} & \num{14791.3} & \num{23} & \num{2002.3} & \num{14359.8} & \cellcolor{resultyellow}\num{26} & \cellcolor{resultyellow}\num{1826.5} & \cellcolor{resultgreen}\num{14418.2} & \num{24} & \num{1877.3} & \cellcolor{resultgreen}\num{14418.2} & \num{25} & \num{1899.2} & \num{14360.2} & \cellcolor{resultyellow}\num{27} & \cellcolor{resultyellow}\num{1795.2} & \cellcolor{resultgreen}\num{14418.5} & \num{26} & \num{1816.1} & \cellcolor{resultgreen}\num{14418.5} \\
    \num{45} & $+$ & \num{45} & \num{30} & \num{1577.6} & \num{15263.4} & \num{20} & \num{2426.2} & \num{14720.7} & \cellcolor{resultyellow}\num{22} & \cellcolor{resultyellow}\num{2169.5} & \cellcolor{resultgreen}\num{14767.0} & \num{22} & \num{2300.2} & \cellcolor{resultgreen}\num{14767.0} & \num{19} & \num{2475.6} & \num{14721.0} & \cellcolor{resultyellow}\num{22} & \cellcolor{resultyellow}\num{2270.4} & \cellcolor{resultgreen}\num{14767.3} & \num{21} & \num{2223.4} & \cellcolor{resultgreen}\num{14767.3} \\
    \midrule avg. &  & \num{225} & \num{190} & \num{840.1} & \num{14149.4} & \num{159} & \num{1344.6} & \num{13694.7} & \cellcolor{resultyellow}\num{169} & \cellcolor{resultyellow}\num{1136.6} & \cellcolor{resultgreen}\num{13920.6} & \num{166} & \num{1186.5} & \cellcolor{resultgreen}\num{13920.6} & \num{159} & \num{1323.1} & \num{13695.7} & \cellcolor{resultyellow}\num{168} & \cellcolor{resultyellow}\num{1176.4} & \cellcolor{resultgreen}\num{13921.6} & \num{167} & \num{1160.0} & \cellcolor{resultgreen}\num{13921.6} \\
    \midrule
    \num{25} & $-$ & \num{45} & \num{41} & \num{333.0} & \num{10730.8} & \num{44} & \num{564.1} & \num{7746.5} & \num{44} & \num{420.8} & \num{8800.9} & \cellcolor{resultyellow}\num{45} & \cellcolor{resultyellow}\num{481.7} & \cellcolor{resultgreen}\num{8806.2} & \num{44} & \num{627.1} & \num{7808.7} & \cellcolor{resultyellow}\num{45} & \cellcolor{resultyellow}\num{454.4} & \num{8832.0} & \num{45} & \num{508.7} & \cellcolor{resultgreen}\num{8835.1} \\
    \num{30} & $-$ & \num{45} & \num{45} & \num{194.1} & \num{9390.1} & \num{29} & \num{1655.9} & \num{5558.8} & \cellcolor{resultyellow}\num{29} & \cellcolor{resultyellow}\num{1618.2} & \num{5675.5} & \num{28} & \num{1621.6} & \cellcolor{resultgreen}\num{5676.9} & \num{27} & \num{1661.0} & \num{5595.4} & \num{26} & \num{1669.8} & \num{5705.9} & \cellcolor{resultyellow}\num{27} & \cellcolor{resultyellow}\num{1655.4} & \cellcolor{resultgreen}\num{5706.4} \\
    \num{35} & $-$ & \num{45} & \num{44} & \num{546.2} & \num{8528.7} & \num{16} & \num{2746.8} & \num{2820.5} & \cellcolor{resultyellow}\num{18} & \cellcolor{resultyellow}\num{2592.7} & \num{2877.2} & \num{18} & \num{2631.8} & \cellcolor{resultgreen}\num{2878.4} & \num{17} & \num{2670.7} & \num{2890.1} & \num{17} & \num{2581.1} & \num{2944.2} & \cellcolor{resultyellow}\num{19} & \cellcolor{resultyellow}\num{2573.1} & \cellcolor{resultgreen}\num{2944.9} \\
    \num{40} & $-$ & \num{45} & \num{30} & \num{1842.1} & \num{7380.1} & \num{6} & \num{3268.0} & \num{-618.6} & \cellcolor{resultyellow}\num{7} & \cellcolor{resultyellow}\num{3208.8} & \num{-560.0} & \num{7} & \num{3236.6} & \cellcolor{resultgreen}\num{-559.3} & \num{7} & \num{3267.3} & \num{-598.0} & \num{8} & \num{3235.2} & \num{-539.9} & \cellcolor{resultyellow}\num{8} & \cellcolor{resultyellow}\num{3199.5} & \cellcolor{resultgreen}\num{-539.7} \\
    \num{45} & $-$ & \num{45} & \num{18} & \num{2484.7} & \num{2200.6} & \num{0} & \num{3600} & \num{-5042.9} & \num{0} & \num{3600} & \num{-4983.7} & \num{0} & \num{3600} & \cellcolor{resultgreen}\num{-4982.8} & \num{0} & \num{3600} & \num{-5013.7} & \num{0} & \num{3600} & \num{-4955.1} & \num{0} & \num{3600} & \cellcolor{resultgreen}\num{-4954.8} \\
    \midrule avg. &  & \num{225} & \num{178} & \num{1080.0} & \num{7646.1} & \num{95} & \num{2366.9} & \num{2092.9} & \cellcolor{resultyellow}\num{98} & \cellcolor{resultyellow}\num{2288.1} & \num{2362.0} & \num{98} & \num{2314.6} & \cellcolor{resultgreen}\num{2363.9} & \num{95} & \num{2365.2} & \num{2136.5} & \num{96} & \num{2308.2} & \num{2397.4} & \cellcolor{resultyellow}\num{99} & \cellcolor{resultyellow}\num{2307.5} & \cellcolor{resultgreen}\num{2398.4} \\
    \midrule
    \num{25} & $+\!/\!-$ & \num{45} & \num{45} & \num{1.2} & \num{-79.8} & \num{45} & \num{4.0} & \num{-564.8} & \num{45} & \num{4.0} & \num{-503.7} & \num{45} & \num{4.4} & \cellcolor{resultgreen}\num{-501.4} & \cellcolor{resultyellow}\num{45} & \cellcolor{resultyellow}\num{4.3} & \num{-552.6} & \cellcolor{resultyellow}\num{45} & \cellcolor{resultyellow}\num{4.3} & \num{-496.5} & \num{45} & \num{4.5} & \cellcolor{resultgreen}\num{-495.7} \\
    \num{30} & $+\!/\!-$ & \num{45} & \num{45} & \num{21.5} & \num{-600.4} & \num{45} & \num{67.4} & \num{-1383.2} & \num{45} & \num{66.6} & \num{-1380.7} & \cellcolor{resultyellow}\num{45} & \cellcolor{resultyellow}\num{63.9} & \cellcolor{resultgreen}\num{-1380.6} & \num{45} & \num{66.5} & \num{-1373.0} & \num{45} & \num{69.0} & \num{-1371.1} & \cellcolor{resultyellow}\num{45} & \cellcolor{resultyellow}\num{70.4} & \cellcolor{resultgreen}\num{-1371.1} \\
    \num{35} & $+\!/\!-$ & \num{45} & \num{45} & \num{26.7} & \num{-1139.4} & \num{45} & \num{328.2} & \num{-2241.2} & \num{45} & \num{333.7} & \num{-2239.9} & \cellcolor{resultyellow}\num{45} & \cellcolor{resultyellow}\num{305.2} & \cellcolor{resultgreen}\num{-2239.6} & \cellcolor{resultyellow}\num{45} & \cellcolor{resultyellow}\num{296.9} & \num{-2235.5} & \num{45} & \num{312.1} & \num{-2234.4} & \num{45} & \num{323.2} & \cellcolor{resultgreen}\num{-2234.3} \\
    \num{40} & $+\!/\!-$ & \num{45} & \num{45} & \num{198.0} & \num{-1590.1} & \num{27} & \num{1662.8} & \num{-3311.5} & \num{28} & \num{1657.2} & \num{-3310.7} & \cellcolor{resultyellow}\num{28} & \cellcolor{resultyellow}\num{1643.3} & \cellcolor{resultgreen}\num{-3310.3} & \num{27} & \num{1682.1} & \num{-3306.6} & \cellcolor{resultyellow}\num{28} & \cellcolor{resultyellow}\num{1655.1} & \num{-3305.7} & \num{28} & \num{1665.8} & \cellcolor{resultgreen}\num{-3305.5} \\
    \num{45} & $+\!/\!-$ & \num{45} & \num{40} & \num{786.7} & \num{-2123.9} & \num{20} & \num{2265.4} & \num{-4489.6} & \cellcolor{resultyellow}\num{22} & \cellcolor{resultyellow}\num{2306.9} & \num{-4488.6} & \num{21} & \num{2321.6} & \cellcolor{resultgreen}\num{-4488.5} & \cellcolor{resultyellow}\num{23} & \cellcolor{resultyellow}\num{2279.1} & \num{-4482.2} & \num{21} & \num{2285.1} & \num{-4481.3} & \num{22} & \num{2262.3} & \cellcolor{resultgreen}\num{-4481.3} \\
    \midrule avg. &  & \num{225} & \num{220} & \num{206.8} & \num{-1106.7} & \num{182} & \num{865.6} & \num{-2398.1} & \cellcolor{resultyellow}\num{185} & \cellcolor{resultyellow}\num{873.7} & \num{-2384.7} & \num{184} & \num{867.7} & \cellcolor{resultgreen}\num{-2384.1} & \num{185} & \num{865.8} & \num{-2390.0} & \num{184} & \num{865.1} & \num{-2377.8} & \cellcolor{resultyellow}\num{185} & \cellcolor{resultyellow}\num{865.2} & \cellcolor{resultgreen}\num{-2377.6} \\
    \midrule\midrule
    total &  & \num{675} & \num{588} & \num{709.0} & \num{6896.2} & \num{436} & \num{1525.7} & \num{4463.2} & \cellcolor{resultyellow}\num{452} & \cellcolor{resultyellow}\num{1432.8} & \num{4632.6} & \num{448} & \num{1456.3} & \cellcolor{resultgreen}\num{4633.5} & \num{439} & \num{1518.0} & \num{4480.8} & \num{448} & \num{1449.9} & \num{4647.1} & \cellcolor{resultyellow}\num{451} & \cellcolor{resultyellow}\num{1444.2} & \cellcolor{resultgreen}\num{4647.5} \\
    \bottomrule
  \end{tabular}%
  }
\end{table*}
For $\mathrm{F^{S+,D^\downarrow}_{GW}}$, adding only the global inequality solves 16 additional instances and reduces the average solution time by approximately 6\%.
Adding the cumulative inequalities provides only a small additional improvement in the LP relaxation and solves four fewer instances than the global variant.
Hence, for the non-hybrid formulation, the cumulative inequalities do not provide a computational advantage over the global inequality alone.

A similar benefit from the global inequality is observed for the hybrid formulation.
The global variant solves nine additional instances compared with the base variant, while the global+cumulative variant solves 12 additional instances and reduces the average solution time by approximately 5\%.
The LP bounds of the global and global+cumulative variants are almost identical, indicating that most of the relaxation improvement is already obtained from the global inequality.

Interestingly, the cumulative inequalities do not further improve the LP-relaxation value for instances with only positive dissimilarity values.
For negative and mixed instances, they provide only a small additional strengthening beyond the global inequality.
Overall, however, the results indicate that most of the improvement in the LP relaxation is already obtained from the global inequality.

Finally, we compare the best non-hybrid and hybrid formulations with the B\&P results of \cite{CLMI26}.
The B\&P algorithm remains stronger overall, solving 588 of the 675 instances, compared with 452 and 451 instances for the best non-hybrid and hybrid variants, respectively.
Its advantage is particularly pronounced for the instance classes with nonpositive- and mixed-sign dissimilarity values.
Nevertheless, the approaches exhibit complementary behavior.
For $n=25$, several of the strengthened MILP variants solve all 135 instances, compared to 127 instances solved by BP(5,10).
Both of the MILP formulations solve 33 instances that BP(5,10) does not solve within the time limit, whereas BP(5,10) solves 165 instances that remain unsolved by both formulations.
The 33 instances uniquely solved by our approaches are concentrated mainly among instances with {positive} dissimilarity values and lower interaction density: 28 have {positive}, 5 have {negative} and none have mixed dissimilarity values, while 24 of the 33 have the lowest interaction density.
These results indicate that, although the specialized B\&P algorithm is more effective overall, the proposed MILP-based approaches remain competitive on particular classes of instances and provide complementary strengths.
Moreover, they are considerably easier to implement, requiring only a standard MILP solver and a simple item-ordering procedure rather than a tailored B\&P framework.

\section{Conclusion}\label{sec:conc}
In this paper, we studied a broad range of MILP formulations for the QBPP.
Starting from classical linearization techniques, we analyzed the relationships among their LP relaxations and developed several strengthened formulations.
We also introduced an item-partitioned hybrid formulation that combines explicitly generated bin patterns with assignment-based representation.
The theoretical analysis shows that formulations based on different linearization techniques can behave quite differently once the symmetry-breaking structure of the QBPP is taken into account.
In particular, although the Fan--Pardalos-based and Glover--Woolsey-based QBPP  formulations have the same LP-relaxation value without symmetry-breaking constraints, the latter dominates after the symmetry-breaking constraints are incorporated.

A central finding of the computational experiments is the importance of item ordering.
Simply reordering the items by decreasing interaction degree significantly improves both the non-hybrid and hybrid formulations {of the QBPP}.
For the non-hybrid formulation, this change alone increases the number of instances solved to optimality by 47.
At the same time, the ordering reduces formulation size and strengthens the LP relaxation.
These two effects are distinct, since an ordering that minimizes formulation size does not necessarily yield the strongest LP relaxation.
The computational experiments also show that item selection plays an important role in the hybrid formulation by determining the balance between the pattern-based and assignment-based components.
Among the additional model strengthening techniques, the global bin-count inequality derived from the underlying BPP provides a clear computational benefit, while the cumulative inequalities yield only limited improvement.

The specialized B\&P algorithm of \cite{CLMI26} remains more effective overall on the benchmark set.
Nevertheless, the computational experiments reveal complementary behavior with the formulations proposed in this paper.
The best non-hybrid and hybrid formulations solve 33 out of 87 instances that are not solved by B\&P within the time limit, primarily those with positive dissimilarities and lower dissimilarity matrix densities.
This suggests that carefully designed MILP formulations can remain effective for particular instance classes.
Moreover, they require significantly less implementation effort than a tailored B\&P algorithm, relying only on a standard MILP solver and a simple item-ordering procedure.

Future work could further exploit this complementarity by incorporating the ordering and strengthening ideas developed here into dedicated exact algorithms.
Another promising direction is to investigate adaptive strategies for the item-partitioned hybrid formulation, including how to choose the size of the set used for pattern generation and which items should be included in this set.
Finally, item ordering itself provides several directions for further study.
The ordering affects the formulation through two distinct mechanisms: it changes both the size of the resulting formulation and the strength of its LP relaxation.
For the non-hybrid formulation, the former objective is equivalent, up to an additive constant, to the minimum sum vertex cover problem on the interaction graph, which minimizes the number of generated pairing variables.
However, this problem is $\mathcal{NP}$-hard, as shown by \cite{FLT04}, and also computationally difficult to solve in practice.
Moreover, even an ordering that is optimal for minimal formulation size need not provide the strongest LP relaxation.
In this paper, we used decreasing interaction degree as a simple proxy that performed well with respect to both effects.
Developing ordering methods that explicitly account for formulation size and LP-relaxation strength, possibly tailored to the formulation and instance structure, is therefore an interesting direction for future research.

\bigskip

\bibliography{sn-bibliography}

\newpage

\section*{APPENDIX} \label{appdx}

\section{Numerical solutions of examples} \label{appdx:num_examples}

\subsection{Optimal solutions for Proposition~\ref{prop:LPrelstrong}}
\label{appdx:LPrelstrong}

For Instance~\ref{inst1}, let $(X_{GW}^*, G^*)$ and $(X_{FP}^*, Z^*)$ denote optimal solutions of
$\mathrm{LP(F^S_{GW})}$ and $\mathrm{LP(F^S_{FP})}$, respectively.
For each $k \in \mathcal{B}$, let $G_k^* = (g_{ijk}^*)$.
An optimal solution of $\mathrm{LP(F^S_{GW})}$ attaining objective value $4.67$ is given by
{\small
\setlength{\arraycolsep}{1.8pt}
\renewcommand{\arraystretch}{0.8} 
\[
X^*_{GW}\!=\!
\begin{bmatrix}
1 & 0 & 0 & 0 \\
0.50 & 0.50 & 0 & 0 \\
0.50 & 0.33 & 0.17 & 0 \\
0.50 & 0.33 & 0.17 & 0
\end{bmatrix},\
G^*_{1}\!=\!
\begin{bmatrix}
- & 0.50 & 0.50 & 0.50 \\
- & - & 0.50 & 0.50 \\
- & - & - & 0 \\
- & - & - & -
\end{bmatrix},\
G^*_{2}\!=\!
\begin{bmatrix}
- & 0 & 0 & 0 \\
- & - & 0.33 & 0.33 \\
- & - & - & 0 \\
- & - & - & -
\end{bmatrix},\
G^*_{3}\!=\!G^*_{4}\!=\!
\begin{bmatrix}
- & 0 & 0 & 0 \\
- & - & 0 & 0 \\
- & - & - & 0 \\
- & - & - & -
\end{bmatrix}.
\]
}
\noindent An optimal solution of $\mathrm{LP(F^S_{FP})}$ attaining objective value $4$ is given by
{\small
\setlength{\arraycolsep}{1.8pt}
\renewcommand{\arraystretch}{0.8} 
\[
X^*_{FP}\!=\!
\begin{bmatrix}
1    & 0    & 0 & 0 \\
0.25 & 0.75 & 0 & 0 \\
0.25 & 0.50 & 0.25 & 0 \\
0    & 0.50 & 0.25 & 0.25
\end{bmatrix}, \quad
Z^*\!=\!
\begin{bmatrix}
- & 0.25 & 0.25 & 0 \\
- & - & 0.75 & 0.75 \\
- & - & - & 0 \\
- & - & - & -
\end{bmatrix}.
\]
}
\subsection{Optimal solution for Example~\ref{ex2}}
\label{appdx:ex2}

\noindent For Instance~\ref{inst1}, let $(X',G')$ denote an optimal solution of
$\mathrm{LP(F^{S'}_{GW})}$, and for $k\in\mathcal B$, let
$G_k'=(g_{ijk}')$.
An optimal solution attaining the objective value $5.33$ is given by
{\fontsize{9pt}{10pt}\selectfont
\setlength{\arraycolsep}{1.8pt}
\renewcommand{\arraystretch}{0.8}
\[
X'\!=\!
\begin{bmatrix}
1 & 0 & 0 & 0 \\
0.50 & 0.50 & 0 & 0 \\
0.50 & 0.17 & 0.33 & 0 \\
0.50 & 0.50 & 0 & 0
\end{bmatrix},\
G'_{1}\!=\!
\begin{bmatrix}
- & 0.50 & 0.50 & 0.50 \\
- & - & 0.50 & 0.50 \\
- & - & - & 0 \\
- & - & - & -
\end{bmatrix},\
G'_{2}\!=\!
\begin{bmatrix}
- & 0 & 0 & 0 \\
- & - & 0.17 & 0.50 \\
- & - & - & 0 \\
- & - & - & -
\end{bmatrix},\
G'_{3}\!=\!G'_{4}\!=\!
\begin{bmatrix}
- & 0 & 0 & 0 \\
- & - & 0 & 0 \\
- & - & - & 0 \\
- & - & - & -
\end{bmatrix}.
\]
}

\subsection{Optimal solutions for Example~\ref{ex3}} \label{appdx:ex3}

For Instance~\ref{inst1}, an optimal solution of $\mathrm{LP(F^{S+}_{FP})}$ attaining the objective value $5$ is given by $(X^{+}_{FP},Z^{+})$, where
{
\small
\setlength{\arraycolsep}{1.8pt}
\renewcommand{\arraystretch}{0.8}
\[
X^{+}_{FP}\!=\!
\begin{bmatrix}
1    & 0    & 0 & 0 \\
0.62 & 0.38 & 0 & 0 \\
0.62 & 0.25 & 0.12 & 0 \\
0.50 & 0.25 & 0.12 & 0.12
\end{bmatrix},
\qquad
Z^{+}\!=\!
\begin{bmatrix}
- & 0.62 & 0.62 & 0.50 \\
- & - & 0.88 & 0.88 \\
- & - & - & 0.12 \\
- & - & - & -
\end{bmatrix}.
\]
}

\noindent For $\mathrm{LP(F^{S+}_{GW})}$, let
$G_k^{+}=(g_{ijk}^{+})$ for each $k\in\mathcal B$.
An optimal solution attaining the objective value $6.67$ is
$(X^{+}_{GW},G^{+})$, where
{\small
\setlength{\arraycolsep}{1.8pt}
\renewcommand{\arraystretch}{0.8} 
\[
X^{+}_{GW}\!=\!
\begin{bmatrix}
1&0&0&0\\
0.50&0.50&0&0\\
0.50&0.17&0.33&0\\
0.50&0.50&0&0
\end{bmatrix},\
G^{+}_{1}\!=\!
\begin{bmatrix}
-&0.50&0.50&0.50\\
-&-&0.50&0.50\\
-&-&-&0\\
-&-&-&-
\end{bmatrix},\
G^{+}_{2}\!=\!
\begin{bmatrix}
-&0&0&0\\
-&-&0.17&0.50\\
-&-&-&0.17\\
-&-&-&-
\end{bmatrix},\
G^{+}_{3}\!=\!G^{+}_{4}\!=\!
\begin{bmatrix}
-&0&0&0\\
-&-&0&0\\
-&-&-&0\\
-&-&-&-
\end{bmatrix}.
\]
}
\subsection{Optimal solutions for the numerical example of Lemma~\ref{lem:G_diag_lqp6}} \label{appdx:G1example}

For Instance~\ref{inst1}, let $(X,T)$ denote an optimal solution of $\mathrm{LP}(\mathrm{F^{S}_{G}})$.
Let $(X',T')$ denote an optimal solution of the LP relaxation obtained from $\mathrm{F^{S}_{G}}$ by replacing constraints~\eqref{G2} for $k=i$ with the strengthened inequalities
\eqref{G_diag_lqp6_valid}.
The corresponding optimal objective values are $3.6$ and $4.5$.
{
\small
\setlength{\arraycolsep}{1.8pt}
\renewcommand{\arraystretch}{0.8} 
\[
X\!=\!
\begin{bmatrix}
1    & 0    & 0 & 0 \\
0.25 & 0.75 & 0 & 0 \\
0.55 & 0.25 & 0.2 & 0 \\
0.05 & 0.75 & 0.2 & 0
\end{bmatrix},\,
T\!=\!
\begin{bmatrix}
4 & 0 & 0 & 0 \\
-5 & -11 & 0 & 0 \\
0 & 0 & 0 & 0
\end{bmatrix},\,
X'\!=\!
\begin{bmatrix}
1    & 0    & 0 & 0 \\
0.44 & 0.56 & 0 & 0 \\
0.81 & 0.19 & 0 & 0 \\
0.19 & 0.56 & 0 & 0.25
\end{bmatrix},\,
T'\!=\!
\begin{bmatrix}
7 & 0 & 0 & 0 \\
-8.75 & -8.25 & 0 & 0 \\
0 & 0 & 0 & 0
\end{bmatrix}.
\]
}
\subsection{Optimal solutions for Proposition~\ref{prop:RLTvsG}}
\label{appdx:RLTvsG}

We report the optimal solutions used in the proof of Proposition~\ref{prop:RLTvsG}.
For each $k\in\mathcal B$, let $\widehat{G}_k=(\widehat g_{ijk})$ and $\overline{G}_k=(\overline g_{ijk})$.
For Instance~\ref{inst1}, an optimal solution of $\mathrm{LP}(\mathrm{F^S_{RLT1}})$ attaining objective value $12$ is
$(\widehat X_R,\widehat G)$, where
{
\small
\setlength{\arraycolsep}{3pt}
\[
\widehat{X}_R =
\begin{bmatrix}
1 & 0 & 0 & 0 \\
0 & 1 & 0 & 0 \\
0 & 0 & 1 & 0 \\
0 & 1 & 0 & 0
\end{bmatrix},\;
\widehat{G}_{1} =
\begin{bmatrix}
1 & 0 & 0 & 0 \\
0 & 0 & 0 & 0 \\
0 & 0 & 0 & 0 \\
0 & 0 & 0 & 0
\end{bmatrix},\;
\widehat{G}_{2} =
\begin{bmatrix}
0 & 0 & 0 & 0 \\
0 & 1 & 0 & 1 \\
0 & 0 & 0 & 0 \\
0 & 1 & 0 & 1
\end{bmatrix},\;
\widehat{G}_{3} =
\begin{bmatrix}
0 & 0 & 0 & 0 \\
0 & 0 & 0 & 0 \\
0 & 0 & 1 & 0 \\
0 & 0 & 0 & 0
\end{bmatrix},\;
\widehat{G}_{4} = \mathbf{0}_{4\times 4}.
\]
}
An optimal solution of
$\mathrm{LP}(\mathrm{F^{S+}_{G}})$ attaining objective value $8.57$
is $(\widehat X_G,\widehat T)$, where
{
\small
\setlength{\arraycolsep}{3pt}
\[
\widehat{X}_G =
\begin{bmatrix}
1 & 0 & 0 & 0 \\
0.57 & 0.43 & 0 & 0 \\
0.86 & 0.14 & 0 & 0 \\
0 & 0.43 & 0 & 0.57
\end{bmatrix},
\qquad
\widehat{T} =
\begin{bmatrix}
4.57 & 0 & 0 & 0 \\
-6.86 & -5.14 & 0 & 0 \\
0 & 0 & 0 & 0
\end{bmatrix}.
\]
}
For Instance~\ref{inst2}, an optimal solution of $\mathrm{LP}(\mathrm{F^S_{RLT1}})$ attaining objective value $6.86$ is $(\overline X_R,\overline G)$, where
{
\small
\setlength{\arraycolsep}{1.4pt}
\[
\overline{X}_R =
\begin{bmatrix}
1 & 0 & 0 & 0 \\
0 & 1 & 0 & 0 \\
0 & 0.86 & 0.14 & 0 \\
0 & 0.86 & 0 & 0.14
\end{bmatrix},\;
\overline{G}_{1} =
\begin{bmatrix}
1 & 0 & 0 & 0 \\
0 & 0 & 0 & 0 \\
0 & 0 & 0 & 0 \\
0 & 0 & 0 & 0
\end{bmatrix},\;
\overline{G}_{2} =
\begin{bmatrix}
0 & 0 & 0 & 0 \\
0 & 1 & 0.86 & 0.86 \\
0 & 0.86 & 0.86 & 0.71 \\
0 & 0.86 & 0.71 & 0.86
\end{bmatrix},\;
\overline{G}_{3} =
\begin{bmatrix}
0 & 0 & 0 & 0 \\
0 & 0 & 0 & 0 \\
0 & 0 & 0.14 & 0 \\
0 & 0 & 0 & 0
\end{bmatrix},\;
\overline{G}_{4} =
\begin{bmatrix}
0 & 0 & 0 & 0 \\
0 & 0 & 0 & 0 \\
0 & 0 & 0 & 0 \\
0 & 0 & 0 & 0.14
\end{bmatrix}.
\]
}
Finally, an optimal solution of
$\mathrm{LP}(\mathrm{F^{S+}_{G}})$ attaining objective value $6.91$
is $(\overline X_G,\overline T)$, where
{
\small
\setlength{\arraycolsep}{3pt}
\[
\overline{X}_G =
\begin{bmatrix}
1 & 0 & 0 & 0 \\
0.45 & 0.55 & 0 & 0 \\
0.55 & 0.45 & 0 & 0 \\
0.09 & 0.55 & 0 & 0.36
\end{bmatrix},
\qquad
\overline{T} =
\begin{bmatrix}
3.64 & 0 & 0 & 0 \\
-5.45 & -6.55 & 0 & 0 \\
0 & 0 & 0 & 0
\end{bmatrix}.
\]
}
\section{Reduced formulations} \label{appendix:reduced_formulations}
\subsection{Reduced \texorpdfstring{${\mathrm{F^S_{FP}}}$}{F S FP} and \texorpdfstring{${\mathrm{F^{S+}_{FP}}}$}{F S FP+} formulations} \label{appdx:FP}
We reduce the size of formulation $\mathrm{F^S_{FP}}$ first by exploiting the symmetry-breaking constraints.
Since $x_{ik}=0$ for $i<k$ by~\eqref{SYM1}, these variables are not introduced.
Hence, constraints~\eqref{SYM1} are implicitly enforced.
Accordingly, constraints~\eqref{FP1} and~\eqref{FP2} need only be considered for $k\le i<j$.
Constraint~\eqref{FP3} is split into two cases.
For $k\le i<j$, the original inequality is retained.
For $i<k\le j$, we have $x_{ik}=0$, and hence~\eqref{FP3} reduces to $z_{ij}\le 1-x_{jk}$.
This yields constraint family~\eqref{FP_sign_neg3} below.

Finally, the formulation can be reduced further by exploiting the sign of the dissimilarity coefficients (also implemented in \cite{CLMI26}).
Define $\mathcal P := \{(i,j)\in\mathcal N\times\mathcal N:\ i<j,\ d_{ij}>0\}$ and $\mathcal M := \{(i,j)\in\mathcal N\times\mathcal N:\ i<j,\ d_{ij}<0\}$.
For $(i,j)\in\mathcal P$, minimization drives $z_{ij}$ toward its smallest feasible value, so only the lower-bounding constraints are required.
For $(i,j)\in\mathcal M$, the objective drives $z_{ij}$ toward its largest feasible value, so only the upper-bounding constraints are required.
If $d_{ij}=0$, the variable $z_{ij}$ and all associated linearization constraints can be eliminated.

The resulting reduced formulation is
\begin{align}
(\mathrm{F^{S}_{FP}}) \qquad
  \min \quad & \sum_{(i,j)\in\mathcal P\cup\mathcal M} d_{ij}z_{ij} + c\sum_{k=1}^{n}x_{kk} \label{FP_sign_obj}\\
  \text{s.t.}\quad
  & \sum_{k=1}^{i}x_{ik}=1  && i\in\mathcal N  \label{FP_sign_assign}\\
  & \sum_{i=k}^{n}w_i x_{ik}\le W x_{kk}  && k\in\mathcal B  \label{FP_sign_capacity}\\
  & x_{ik}\le x_{kk}  && i\in\mathcal N,\ k\in\mathcal B,\ k<i  \label{FP_sign_sym}\\
  & x_{ik}+x_{jk}-1\le z_{ij}  && (i,j)\in\mathcal P,\ k\le i  \label{FP_sign_pos}\\
  & z_{ij}\le 1-x_{ik}+x_{jk}  && (i,j)\in\mathcal M,\ k\le i  \label{FP_sign_neg1}\\
  & z_{ij}\le 1+x_{ik}-x_{jk}  && (i,j)\in\mathcal M,\ k<i  \label{FP_sign_neg2}\\
  & z_{ij}\le 1-x_{jk}  && (i,j)\in\mathcal M,\ i<k\le j  \label{FP_sign_neg3}\\
  & z_{ij}\ge0  && (i,j)\in\mathcal P\cup\mathcal M  \label{FP_sign_zdomain}\\
  & x_{ik}\in\{0,1\}  && i\in\mathcal N,\ k\le i.  \label{FP_sign_binary}
\end{align}
For the strengthened formulation $\mathrm{F^{S+}_{FP}}$, constraints~\eqref{FP_sign_pos} are replaced by the strengthened inequalities~\eqref{eq:strong-z}, restricted to the reduced index set:
\begin{equation}
x_{ik}+x_{jk}-x_{kk}\le z_{ij} \qquad (i,j)\in\mathcal P,\quad k\le i. \label{FP_lower_plus}
\end{equation}
All other constraints remain unchanged.

\subsection{Reduced \texorpdfstring{$\mathrm{F^S_{GW}}$}{F S GW} and \texorpdfstring{$\mathrm{F^{S+}_{GW}}$}{F S GW+} formulations}
\label{appdx:GW}
We first reduce the size of formulation $\mathrm{F^S_{GW}}$ by exploiting the symmetry-breaking constraints~\eqref{SYM1} and \eqref{SYM2}.
Since $x_{ik}=0$ for $i<k$ by~\eqref{SYM1}, these variables are not introduced.
Moreover, for $i<j$ and $k>i$, constraint~\eqref{GWL1} together with~\eqref{SYM1} implies
$g_{ijk}\le x_{ik}=0$.
Hence, the variables $g_{ijk}$ with $k>i$ can also be eliminated, and the remaining pairing variables are indexed by $k\le i<j$.

The model can be further reduced by exploiting the sign of the dissimilarity coefficients (which is also implemented in \cite{CLMI26}).
Using the sets $\mathcal P:=\{(i,j):1\le i<j\le n,\ d_{ij}>0\}$ and $\mathcal M:=\{(i,j):1\le i<j\le n,\ d_{ij}<0\}$,
for $(i,j)\in\mathcal P$, minimization drives $g_{ijk}$ toward its smallest feasible value, so only the lower-bounding constraints \eqref{GWL3}, together with nonnegativity constraints~\eqref{GWL4}, are required.
For $(i,j)\in\mathcal M$, the objective drives $g_{ijk}$ toward its largest feasible value, so only the upper-bounding constraints \eqref{GWL1} and \eqref{GWL2}, together with nonnegativity~\eqref{GWL4}, are required.
If $d_{ij}=0$, the corresponding $g_{ijk}$ variables and their associated linearization constraints can be eliminated.

The resulting reduced formulation is
\begin{align}
(\mathrm{F^S_{GW}})\qquad
\min\quad & \sum_{(i,j)\in\mathcal P\cup\mathcal M} \sum_{k=1}^{i} d_{ij}g_{ijk} +c\sum_{k=1}^{n}x_{kk} \label{GW_sign_base_obj}\\
\text{s.t.}\quad
& \sum_{k=1}^{i}x_{ik}=1 && i\in\mathcal N \label{GW_sign_base_assign}\\
& \sum_{i=k}^{n}w_i x_{ik}\le W x_{kk} && k\in\mathcal B \label{GW_sign_base_capacity}\\
& x_{ik}\le x_{kk} && k<i \label{GW_sign_base_sym}\\
& x_{ik}+x_{jk}-1\le g_{ijk} && (i,j)\in\mathcal P, k\le i \label{GW_sign_base_pos}\\
& g_{ijk}\le x_{ik} && (i,j)\in\mathcal M, k\le i \label{GW_sign_base_neg1}\\
& g_{ijk}\le x_{jk} && (i,j)\in\mathcal M, k\le i \label{GW_sign_base_neg2}\\
& g_{ijk}\ge0 && (i,j)\in\mathcal P\cup\mathcal M, k\le i \label{GW_sign_base_gdomain}\\
& x_{ik}\in\{0,1\} && k\le i. \label{GW_sign_base_binary}
\end{align}

For the strengthened formulation $\mathrm{F^{S+}_{GW}}$, the lower-bounding inequalities \eqref{GW_sign_base_pos} are replaced by the strengthened inequalities~\eqref{eq:strong-g} from Lemma~\ref{lemma:stronger} using the reduced index set.
This strengthening also permits an additional variable reduction: consider $k=i$ in $\mathrm{F^{S+}_{GW}}$.
Constraints~\eqref{GWL2} and~\eqref{eq:strong-g} give $g_{iji}\le x_{ji}$
and $x_{ii}+x_{ji}-x_{ii}\le g_{iji}$, respectively.
Hence, $g_{iji}=x_{ji}$ for $1\le i<j\le n$.
Therefore, the variables $g_{iji}$ can be eliminated without changing either the integer feasible region or the LP relaxation.
Consequently, their contribution to the objective function can be replaced by $\sum_{i=1}^{n-1}\sum_{j=i+1}^{n}d_{ij}x_{ji}$.
Therefore, the remaining pairing variables need only be defined for $k<i<j$.
The resulting reduced formulation is
\begin{align}
(\mathrm{F^{S+}_{GW}})\
\min\quad & \sum_{i=1}^{n-1}\sum_{j=i+1}^{n}d_{ij}x_{ji} + \sum_{(i,j)\in\mathcal P\cup\mathcal M} \sum_{k=1}^{i-1}d_{ij}g_{ijk} +c\sum_{k=1}^{n}x_{kk} \label{GW_sign_obj}\\
\text{s.t.}\quad
& \sum_{k=1}^{i}x_{ik}=1 && i\in\mathcal N \label{GW_sign_assign}\\
& \sum_{i=k}^{n}w_i x_{ik}\le W x_{kk} && k\in\mathcal B \label{GW_sign_capacity}\\
& x_{ik}\le x_{kk} && k<i \label{GW_sign_sym}\\
& x_{ik}+x_{jk}-x_{kk}\le g_{ijk} && (i,j)\in\mathcal P,\ k<i\label{GW_sign_pos}\\
& g_{ijk}\le x_{ik} && (i,j)\in\mathcal M,\ k<i \label{GW_sign_neg1}\\
& g_{ijk}\le x_{jk} && (i,j)\in\mathcal M,\ k<i \label{GW_sign_neg2}\\
& g_{ijk}\ge0 && (i,j)\in\mathcal P\cup\mathcal M,\ k<i \label{GW_sign_gdomain}\\
& x_{ik}\in\{0,1\} && k\le i. \label{GW_sign_binary}
\end{align}

\subsection{Reduced \texorpdfstring{${\mathrm{F^S_{RLT1}}}$}{F S RLT1} formulation} \label{appdx:RLT}

We reduce the size of formulation $\mathrm{F^S_{RLT1}}$ by eliminating variables that are fixed to zero by the symmetry-breaking constraints.
In particular, the variables $x_{ik}$ with $i<k$ are not introduced, since they are fixed to zero by~\eqref{SYM1}.
Similarly, the variables $\hat g_{ijk}$ with $\min\{i,j\}<k$ are not introduced.
Indeed, if $i<k$, then $\hat g_{ijk}\le x_{ik}=0$ by~\eqref{SYM1} and~\eqref{R1}.
Similarly, if $j<k$, the same conclusion follows from the symmetry relation $\hat g_{ijk}=\hat g_{jik}$ in~\eqref{R3}. Hence, $\hat g_{ijk}$ can be nonzero only when $k\le \min\{i,j\}$.

Moreover, constraints~\eqref{R2} are redundant and can be omitted, as shown in Lemma~\ref{lem:redundancyR2}.
After this omission, all remaining constraints are restricted to the reduced variable domains. Thus, constraints~\eqref{R1}, \eqref{R3}, and~\eqref{Connector} are imposed only for $k\le \min\{i,j\}$, $k\le i<j$, and $k\le i$, respectively.
The same restriction is applied to the RLT constraints. The constraints derived from the capacity constraints, namely~\eqref{RLTC1_1+} and~\eqref{RLTC2_1-}, are imposed only for $j>k$, with summations over $i\ge k$.
The RLT constraints derived from the symmetry-breaking inequalities, namely~\eqref{RLTS_1+} and~\eqref{RLTS_1-}, are imposed only for index pairs for which all variables are defined, namely $k<i$ and $k<j$.
The objective-function summations are adjusted accordingly.
\begin{align}
(\mathrm{F^S_{RLT1}}) \ \min \ & \frac{1}{2} \sum_{k \in \mathcal{B}}\sum_{i=k}^{n} \sum_{j=k}^{n} d_{ij}\hat{g}_{ijk}+c\sum_{k \in \mathcal{B}}x_{kk} \\ 
\text{s.t.} \quad
& \sum_{k=1}^{i} x_{ik} = 1 & i \in \mathcal{N} \label{RLT_R_assign}\\
& \sum_{i=k}^{n} w_i x_{ik} \le W x_{kk}  & k \in \mathcal{B}  \label{RLT_R_capacity} \\
& x_{ik} \le x_{kk}  &i\in \mathcal{N},k \in \mathcal{B}, k<i \label{RLT_R_sym} \\
&\hat{g}_{ijk} \leq x_{ik} & i,j\in \mathcal{N}, k\in \mathcal{B},k \le \min\{i,j\}\\
& \hat{g}_{ijk} = \hat{g}_{jik} & i,j\in \mathcal{N}, k\in \mathcal{B}, k\le i<j \\
&\sum_{i=k}^{n} w_i \hat{g}_{ijk}\leq W\, \hat{g}_{kjk} & j \in \mathcal{N}, k \in \mathcal{B}, k<j \\
&\sum_{i=k}^{n} w_i (x_{ik}-\hat{g}_{ijk}) \leq W (x_{kk} - \hat{g}_{kjk}) &   j \in \mathcal{N}, k \in \mathcal{B}, k< j \\
&\hat{g}_{ijk} \leq \hat{g}_{jkk} & i,j \in \mathcal{N}, k \in \mathcal{B}, k<\min\{i,j\} \\
&x_{ik}-x_{kk} \le \hat{g}_{ijk} - \hat{g}_{jkk} & i,j\in \mathcal{N},k \in \mathcal{B}, k<\min\{i,j\} \\
& x_{ik} = \hat{g}_{iik}  & i\in\mathcal{N}, k \in \mathcal{B}, k\le i \\
& \hat{g}_{ijk} \ge 0 & i,j\in \mathcal{N}, k\in \mathcal{B}, k \le \min\{i,j\}\\
& x_{ik} \in \{0,1\}  &i\in \mathcal{N},k \in \mathcal{B}, k\le i.  \label{RLT_R_binary}
\end{align}

\subsection{Reduced \texorpdfstring{$\mathrm{F^{S+}_{G}}$}{F S+ G} formulation}
\label{appdx:G}
Because of the symmetry-breaking constraints~\eqref{SYM1}, the variables $x_{ik}$ and $t_{ik}$ with $i<k$ can be omitted.
Accordingly, the constraints involving these variables are restricted to $k\le i$.
The refined bounds $\tilde l_{ik}$ and $\tilde u_{ik}$ apply to $t_{ik}$ for all $k\le i<n$.
For $k<i<n$, the constraints~\eqref{G2} use the refined deactivation coefficients $\tilde M_{ik}$.
For $k=i$, those constraints are replaced by the strengthened inequalities~\eqref{G_diag_lqp6_valid} from Lemma~\ref{lem:G_diag_lqp6}, so no coefficient $\tilde M_{ii}$ is needed.
The resulting reduced formulation is
\begin{align}
(\mathrm{F^{S+}_{G}}) \quad
\min \quad & \sum_{i=1}^{n-1} \sum_{k=1}^{i} t_{ik} + c \sum_{k\in\mathcal{B}} x_{kk} \\
\text{s.t.} \quad& \sum_{k=1}^{i} x_{ik} = 1 &  i \in \mathcal{N} \\
& \sum_{i=k}^{n} w_i x_{ik} \le W x_{kk} & k \in \mathcal{B}   \\
& x_{ik} \le x_{kk} & i\in \mathcal{N},k \in \mathcal{B}, k<i  \\[1.5mm]
& \tilde l_{ik} x_{ik} \le t_{ik} &i\in \mathcal{N}\setminus\{n\},k \in \mathcal{B}, k\le i  \\
& \sum_{j=i+1}^{n} d_{ij} x_{jk} - \tilde M_{ik} (1-x_{ik}) \le t_{ik} &i\in \mathcal{N}\setminus\{n\},k \in \mathcal{B}, k<i\\
& \sum_{j=i+1}^{n} d_{ij} x_{ji}\le t_{ii} &i\in \mathcal{N}\setminus\{n\} \\
& \tilde l_{ik}\le t_{ik}\le\tilde u_{ik} &i\in \mathcal{N}\setminus\{n\},k \in \mathcal{B}, k\le i \\
& x_{ik}\in\{0,1\} &i\in \mathcal{N},k \in \mathcal{B}, k\le i.
\end{align}

\subsection{Reduced item-partitioned hybrid formulation.} \label{appdx:PGW}

We reduce the size of formulation $\mathrm{F^{S+}_{HI}}$ by first exploiting the symmetry-breaking constraints~\eqref{PGW_SYM1} and \eqref{PGW_SYM2}.
Since $x_{ik}=0$ for $i<k$ by~\eqref{PGW_SYM1}, these variables are not introduced.
Moreover, for $i<j$ and $i<k$, constraints~\eqref{PGW_SYM1} and \eqref{PGW_gneg1} imply $g_{ijk}\le x_{ik}=0$.
Hence, the variables $g_{ijk}$ with $i<k$ can also be eliminated, and the remaining pairing variables are indexed by $k\in\bar{\mathcal A}$ with $k\le i<j$.

Constraints~\eqref{PGW_gpos} allow further variable reduction. For $i\in\bar{\mathcal{A}}$ and $i<j$, consider the $g$-variable corresponding to $k=i$.
Constraints~\eqref{PGW_gneg2} and~\eqref{PGW_gpos} give $g_{iji}\le x_{ji}$ and $x_{ii}+x_{ji}-x_{ii}\le g_{iji}$, respectively.
Hence, $g_{iji}=x_{ji}$ for $i\in\bar{\mathcal A},\quad i<j$.
Therefore, these variables can be eliminated.
Their contribution to the objective function can consequently be replaced by $\sum_{i\in\bar{\mathcal A}}\sum_{j=i+1}^{n} d_{ij}x_{ji}$,
and the remaining $g$-variables need only be defined for $k\in\bar{\mathcal{A}}$ with $k<i<j$.
Finally, as for $\mathrm{F^{S+}_{GW}}$, the model can be further reduced by exploiting the sign of the dissimilarity coefficients.
Let $\mathcal P:=\{(i,j):1\le i<j\le n,\ d_{ij}>0\}$, and $\mathcal M:=\{(i,j):1\le i<j\le n,\ d_{ij}<0\}$.
The resulting reduced formulation is
\begin{align}
\hspace{-2.2cm}
(\mathrm{F^{S+}_{HI}})\
\min\quad
& \sum_{p\in\mathcal P(\mathcal A)}\rho_p\lambda_p
+ \sum_{i\in\bar{\mathcal A}}\sum_{j=i+1}^{n}d_{ij}x_{ji}
+ \sum_{(i,j)\in\mathcal P\cup\mathcal M}
  \sum_{\substack{k\in\bar{\mathcal A}\\ k<i}}
  d_{ij}g_{ijk}
+ c\sum_{k\in\bar{\mathcal A}}x_{kk}
\label{PGW_red_obj}
\end{align}
\vspace{-1.2\baselineskip}
\begin{align}
\text{s.t.}\quad
& \sum_{p\in\mathcal P(\mathcal A)}\delta_{ip}\lambda_p
+ \sum_{\substack{k\in\bar{\mathcal A}\\ k\le i}}x_{ik} =1 && i\in\mathcal N \label{PGW_red_assign}\\
& x_{ik}\le x_{kk} && i\in\mathcal N,\ k\in\bar{\mathcal A},\ k<i \label{PGW_red_sym}\\
& \sum_{i=k}^{n}w_i x_{ik}\le W x_{kk} && k\in\bar{\mathcal A} \label{PGW_red_cap}\\
& x_{ik}+x_{jk}-x_{kk}\le g_{ijk} && (i,j)\in\mathcal P,\ k\in\bar{\mathcal A},\ k<i \label{PGW_red_pos}\\
& g_{ijk}\le x_{ik} && (i,j)\in\mathcal M,\ k\in\bar{\mathcal A},\ k<i \label{PGW_red_neg1}\\
& g_{ijk}\le x_{jk} && (i,j)\in\mathcal M,\ k\in\bar{\mathcal A},\ k<i \label{PGW_red_neg2}\\
& g_{ijk}\ge 0 && (i,j)\in\mathcal P\cup\mathcal M,\ k\in\bar{\mathcal A},\ k<i \label{PGW_red_gdomain}\\
& x_{ik}\in\{0,1\} && i\in\mathcal N,\ k\in\bar{\mathcal A},\ k\le i \label{PGW_red_xbin}\\
& \lambda_p\in\{0,1\} && p\in\mathcal P(\mathcal A). \label{PGW_red_lambin}
\end{align}

\subsection{The Fortet--Glover--Woolsey linearization approach  \texorpdfstring{(see \cite{CLMI26})}{}.}\label{appdx:EFGW}
\begin{align}
(\mathrm{F_{EFGW}})\,
\min\, & \sum_{k\in\mathcal B}  \sum_{i=1}^{n-1}\sum_{j=i+1}^{n}  d_{ij}g_{ijk} + c\sum_{k\in\mathcal B}x_{kk} \label{EFGW_obj}\\
\text{s.t.}\,
& \sum_{k\in\mathcal B}x_{ik}=1&& i\in\mathcal N \label{EFGW_assign}\\
& \sum_{i\in\mathcal N}w_i x_{ik}  \le W x_{kk}&& k\in\mathcal B\label{EFGW_capacity}\\
& x_{ik}=0&& i\in\mathcal N,\ k\in\mathcal B,\ i<k\label{EFGW_sym1}\\
& x_{ik}\le x_{kk}&& i\in\mathcal N,\ k\in\mathcal B,\ k<i\label{EFGW_sym2}\\
& \sum_{\substack{i\in\mathcal N\\i\neq j}}  w_i g_{ijk}  \le  Wg_{kjk}-w_jx_{jk}&& j\in\mathcal N,\ k\in\mathcal B\label{EFGW_RLT1}\\
& \sum_{\substack{i\in\mathcal N\\i\neq j}}  w_i\bigl(x_{ik}-g_{ijk}\bigr)  \le  W\bigl(x_{kk}-g_{kjk}\bigr)&& j\in\mathcal N,\ k\in\mathcal B\label{EFGW_RLT2}\\
& g_{ijk}\le x_{ik}&& i\in\mathcal N,\ k\in\mathcal B,\    j\in\mathcal N\setminus\{i\},\   d_{ij}<0\label{EFGW_neg}\\
& x_{ik}+x_{jk}-1\le g_{ijk}&& i\in\mathcal N,\ k\in\mathcal B,\    j\in\mathcal N\setminus\{i\},\
   d_{ij}>0\label{EFGW_pos}\\
& g_{ijk}=g_{jik}&& i,j\in\mathcal N,\ k\in\mathcal B,\ i<j\label{EFGW_gsym}\\
& g_{ijk}\ge0&& i\in\mathcal N,\ k\in\mathcal B,\    j\in\mathcal N\setminus\{i\},\   d_{ij}\neq0
\label{EFGW_gdomain}\\
& x_{ik}\in\{0,1\}&& i\in\mathcal N,\ k\in\mathcal B.
\label{EFGW_binary}
\end{align}

\end{document}